\documentclass[twoside]{article}
\usepackage{har2nat,natbib}
\PassOptionsToPackage{hyphens}{url}\usepackage[
    breaklinks=true,
    colorlinks,
    linkcolor={red!50!black},
    citecolor={blue!50!black},
    urlcolor={blue!80!black}
    ]{hyperref}
\usepackage{doi}
\usepackage{amsmath,amssymb} 
\usepackage{afterpage,rotating,pdflscape} 

\allowdisplaybreaks
\usepackage[UKenglish]{babel} 
	
\newcommand{\ignore}[1]{}
\newcommand{\citen}{\citeasnoun*}
\newcommand{\citeb}[1]{\citename*{#1}\ \citeyear*{#1}}

\usepackage[margin=1.19in]{geometry}
\usepackage{subfigure}
\usepackage{marginnote}
\newcommand{\anlchange}[1]{\textcolor{black}{#1}}
\newcommand{\ANL}[1]{\textcolor{black}{#1}}

\newcommand{\fo}[1]{{#1}} 

\newcommand*{\ie}{{\it i.e.}}

\newcommand*{\eg}{{\it e.g.}}

\numberwithin{figure}{section}
\numberwithin{table}{section}

\newcommand{\gfrac}[2]{{#1}/{#2}}
\newcommand{\gfracp}[2]{{#1}/{(#2)}}
\newcommand{\gpfrac}[2]{{(#1)}/{#2}}
\newcommand{\gpfracp}[2]{{(#1)}/{(#2)}}

\newcommand{\ffrac}{\displaystyle \frac}
 
\newcommand{\ds}{\displaystyle}
\newcommand{\D}{\mathrm{d}}

\newcommand{\gforall}{\text{for all}}
\newcommand{\Bigoh}{O}

\usepackage{color} 
\usepackage{graphics,graphicx} 
\graphicspath{ {./Figures/}}
\usepackage[ruled,vlined,linesnumbered]{algorithm2e}
\SetKwInOut{Input}{input}
\SetKwInOut{Output}{output}
\SetKw{Break}{break}
\SetAlgoNlRelativeSize{-5}
\SetKw{Continue}{continue}
\DontPrintSemicolon

\usepackage{epstopdf}

\usepackage{tikz}
\usetikzlibrary{positioning,shadows,arrows}

\newcommand{\T}{{\rm T}} 
\newcommand{\Reals}{\mathbb{R}} 
\newcommand{\Ints}{{\mathbb{Z}}} 

\newcommand{\eps}{\varepsilon} 			
\newcommand{\codes}[1]{{\rm \textsf{#1}}} 
\newcommand{\framework}[1]{{\rm #1}}
\DeclareMathOperator*{\argmin}{arg\,min}

\newcommand{\Proj}[1]{\mathrm{proj}( #1 )}
\newcommand{\Deltamax}{\Delta_{\rm max}}
\newcommand{\Deltamin}{\Delta_{\rm min}}
\newcommand{\kappaef}{\kappa_{{\rm ef}}}
\newcommand{\kappaeg}{\kappa_{{\rm eg}}}
\newcommand{\kappaeH}{\kappa_{{\rm eH}}}

\newcommand{\gammad}{\gamma_{\rm dec}}
\newcommand{\kappad}{\kappa_{\rm d}}
\newcommand{\gammai}{\gamma_{\rm inc}}

\newcommand{\defined}{=} 
\newcommand{\gS}{\gb^{\rm S}} 
\newcommand{\gK}{\gb^{\rm K}} 
\newcommand{\bigo}[1]{O ( #1 )}

\newcommand{\minimize}{\operatornamewithlimits{minimize}}

\renewcommand{\P}[1]{\operatorname{\mathbb{P}}[ #1]} 
\renewcommand{\dim}[1]{\operatorname{\textrm{dim}}( #1)} 
\newcommand{\qn}{(n+1)(n+2)/2} 
\newcommand{\cm}{\mathrm{cm}}

\newcommand{\B}[2]{\mathcal{B}(#1 ; #2 )} 
\newcommand{\Lip}[1]{L_{\rm #1}} 
\newcommand{\E}[2][]{\operatorname{\mathbb{E}}_{#1}\left[ #2\right]} 
\newcommand{\Egg}[2][]{\operatorname{\mathbb{E}}_{#1}\Biggl[ #2\Biggr]} 
\newcommand{\Ea}[2][]{\operatorname{\mathbb{E}}_{#1}[ #2]} 
\newcommand{\df}[1]{\delta_{\rm f}(#1)} 
\newcommand{\dc}[1]{\delta_{\rm c}(#1)} 
\newcommand{\trem}[1]{$^{\textrm #1}$} 

\newcommand{\tabref}[1]{\fo{Table~\ref{table:#1}}}
\newcommand{\algref}[1]{\fo{Algorithm~\ref{alg:#1}}}
\newcommand{\secref}[1]{Section~\ref{sec:#1}}
\newcommand{\secrefs}[2]{Sections~\ref{sec:#1} and \ref{sec:#2}}
\newcommand{\lineref}[1]{{line}~\ref{line:#1}}

\newcommand{\zerob}{\mathbf{0}} 
\usepackage{bm} 
\newcommand{\phib}{\bm{\phi}}
\newcommand{\Phib}{\bm{\Phi}}

\newcommand{\kappab}{\bm{\kappa}}
\newcommand{\thetab}{\bm{\theta}}
\newcommand{\ab}{\bm{a}}
\newcommand{\Ab}{\bm{A}}
\newcommand{\bb}{\bm{b}}
\newcommand{\Bb}{\bm{B}}
\newcommand{\cb}{\bm{c}}
\newcommand{\db}{\bm{d}}
\newcommand{\Db}{\bm{D}}
\newcommand{\eb}{\bm{e}}

\newcommand{\Fb}{\bm{F}}
\newcommand{\gb}{\bm{g}}

\newcommand{\Hb}{\bm{H}}
\newcommand{\Ib}{\bm{I}}
\newcommand{\lb}{\bm{l}}
\newcommand{\Mb}{\bm{M}}
\newcommand{\pb}{\bm{p}}
\newcommand{\Pb}{\bm{P}}
\newcommand{\Rx}{R_{\xb}} 
\newcommand{\Rlevel}{R_{\textrm{level}}} 
\newcommand{\sba}{\bm{s}} 
\newcommand{\Sb}{\bm{S}}
\newcommand{\ub}{\bm{u}}
\newcommand{\Ub}{\bm{U}}

\newcommand{\vb}{\bm{v}}

\newcommand{\xb}{\bm{x}}

\newcommand{\yb}{\bm{y}}
\newcommand{\Yb}{\bm{Y}}
\newcommand{\zb}{\bm{z}}
\newcommand{\Omegab}{\bm{\Omega}}

\newcommand{\xib}{\bm{\xi}}
\newcommand{\Xib}{\bm{\Xi}}
\newcommand{\lambdab}{\bm{\lambda}}
\newcommand{\Iset}{I}

\newcommand{\cC}{\mathcal{C}} 
\newcommand{\cLC}{\mathcal{LC}} 
 
\newcommand{\cH}{\mathcal{H}} 
\newcommand{\cK}{\mathcal{K}}
\newcommand{\cL}{\bm{\mathcal{L}}} 
\newcommand{\cN}{\mathcal{N}} 
\newcommand{\cP}{\mathcal{P}} 

\usepackage{lipsum}

\newcommand{\Nonord}{\mathbf{N}}
\newcommand{\epsb}{\bm{\epsilon}}

\usepackage{fancyhdr}
\begin{document}

\title{\textsc{Derivative-free optimization methods}}

\author
{Jeffrey Larson, Matt Menickelly and Stefan M.\ Wild
\\~\\
\textit{Mathematics and Computer Science Division,}\\ 
\textit{Argonne National Laboratory, Lemont, IL 60439, USA}\\
{\small
{\sf jmlarson@anl.gov}, {\sf mmenickelly@anl.gov}, {\sf wild@anl.gov}}
}

\maketitle

{\em Dedicated to the memory of Andrew R.\ Conn for his inspiring enthusiasm and his
 many contributions to the renaissance of derivative-free optimization methods.}
\vspace{30pt}

\begin{abstract}
In many optimization problems arising from scientific, engineering and
artificial intelligence applications,
objective and constraint functions are available only as the output of a black-box
or simulation oracle that does not provide derivative information.
Such settings necessitate the use of methods for derivative-free\ANL{, or zeroth-order,} optimization.
We provide a review and perspectives on
developments in these methods, with an emphasis on highlighting recent
developments and \ANL{on} unifying treatment of such problems in the non-linear
optimization and machine learning literature.
We categorize methods based on assumed properties of the black-box functions,
as well as features of the method\ANL{s}. We first overview the primary setting of
deterministic methods applied to unconstrained, non-convex optimization problems
where the objective function is defined by a deterministic black-box oracle.
We then discuss developments in randomized methods, methods that assume
some additional structure about the objective (including convexity,
separability and general non-smooth compositions),
methods for problems where the output of the black-box oracle is stochastic,
and methods for handling different types of constraints. 
\end{abstract}

\tableofcontents

\vspace{5pt}
\section{Introduction} %
\label{sec:intro}
The growth in computing for scientific, engineering and social applications 
has long been a driver of advances in methods for numerical
optimization. 
The development of derivative-free optimization methods -- those 
methods that do not require the availability of 
derivatives -- has especially been driven by the need to optimize increasingly 
complex and diverse problems.
One of the earliest calculations on 
MANIAC,\footnote{Mathematical Analyzer, Integrator, And Computer. Other 
lessons learned from this application are discussed by \citeasnoun{Anderson1986}.}
an early computer based on the von Neumann architecture, was the 
approximate solution of a six-dimensional non-linear 
least-squares problem using a derivative-free coordinate search \cite{Fermi1952}. Today, 
derivative-free methods are used 
routinely, for example by Google 
\cite{Golovin2017}, for the automation and tuning needed in the artificial 
intelligence era.

\begin{figure} %
\centering
\includegraphics[width=255pt,bb=0 0 370 300,clip]{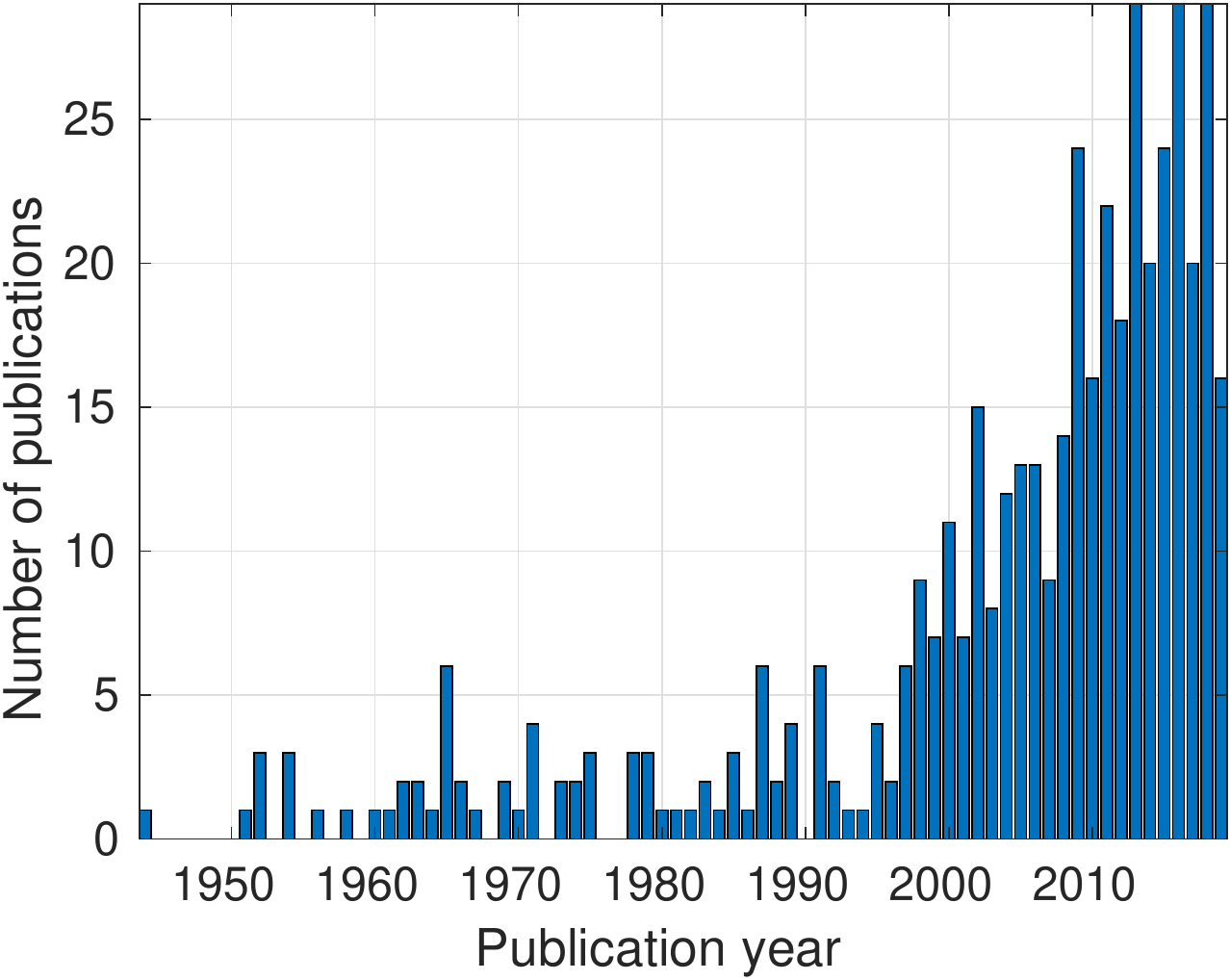}\\ %
\parbox{280pt}{\caption{Histogram of the references cited in the bibliography.\label{fig:references}}}
\end{figure}

In this paper we survey methods for derivative-free optimization and
key results for their analysis. 
Since the field -- also referred to 
as black-box optimization, gradient-free optimization, optimization 
without derivatives, simulation-based optimization and zeroth-order 
optimization -- is now far too expansive for a single survey, we focus 
on methods for local optimization of continuous-valued, single-objective 
problems. Although \secref{extensions} illustrates further
connections, here we 
mark the following notable omissions.
\begin{itemize}
\setlength\itemsep{5pt}
 \item We focus on methods that seek a 
local minimizer. Despite 
users understandably desiring the best possible solution, the problem of 
global optimization 
raises innumerably more mathematical and computational challenges than do the  
methods presented here. We instead point to the survey by \citeasnoun{Neumaier2004},
which importantly addresses general constraints,
and to
the textbook by
\citeasnoun{Forrester2008}, which lays a foundation for global surrogate modelling.

\item Multi-objective optimization and optimization in the 
presence of discrete variables are similarly popular tasks among users. Such 
problems 
possess fundamental challenges as well as differences from the methods 
presented here. 

\item In focusing on methods, we cannot do justice to the application problems 
that have driven the development of derivative-free 
methods and benefited from implementations of these methods. 
The recent textbook by \citeasnoun{AudetHare2017} contains a number 
of examples and references to applications; 
 \citeasnoun{Rios2013} and \citeasnoun{Auger2009} both reference a diverse set of 
implementations.
At the persistent page 
\[\mbox{\textsf{https://archive.org/services/purl/dfomethods}}\]
we intend to \anlchange{link} all works that cite the entries in our bibliography and those 
that cite this
survey; we hope this will provide a coarse, but dynamic, catalogue for the
reader interested in potential uses of these methods.
\end{itemize}
Given these limitations, we particularly note the intersection with the 
foundational books by 
\citeasnoun{Kelleybook} 
and 
\citeasnoun{Conn2009a}. 
Our intent is to highlight recent developments in, and the evolution of,
derivative-free optimization methods. \fo{Figure~\ref{fig:references}} summarizes 
our bias; over half of the references in this survey are from the past ten
years.

Many of the fundamental inspirations for the methods
discussed in this survey are
detailed to a lesser extent. We note in particular the activity in 
the United Kingdom in the 1960s (see \eg\ the works by
\citeb{Rosenbrock1960}, \citeb{Powell1964}, 
\citeb{NelderMead}, \citeb{Fletcherdfo65} and \citeb{Box1966}, and the later 
exposition and expansion by \citeb{brent1973algorithms})
and the Soviet Union (as evidenced by \citeb{Rastrigin1963}, 
\citeb{Matyas1965}, \anlchange{\citeb{Karmanov1974},} \citeb{polyakbook} and others). 
In addition to those mentioned later, we single out the work of
\citeasnoun{Powell1975a}, 
\citeasnoun{Wright95}, \citeasnoun{Davis2005Powell} and \citeasnoun{Leyffer2015optima} for 
insight into some of 
these early pioneers.

With our focus clear, we turn our attention to the deterministic optimization 
problem
\begin{equation}\label{eq:det_prob}
\begin{aligned}
 &\anlchange{\underset{\xb}{\minimize}} && f(\xb) \\
 &\mbox{\,subject to} && \xb \in \Omegab \subseteq \Reals^n
\end{aligned}
\tag{DET}\end{equation}
and the stochastic optimization problem
\begin{equation}\label{eq:stoch_prob}
\begin{aligned}
 &\anlchange{\underset{\xb}{\minimize}} && f(\xb) \defined 
\E[\xib]{\tilde{f}(\xb; \xib)} \\
 &\mbox{\,subject to} && \xb \in \Omegab.
\end{aligned}
\tag{STOCH}\end{equation}
Although important exceptions are noted throughout this survey, the majority of 
the methods discussed assume that the objective 
function $f$ in \eqref{eq:det_prob} and \eqref{eq:stoch_prob} is 
differentiable. 
This assumption may cause readers to pause (and some readers 
may never resume). The methods considered here do not 
necessarily address non-smooth optimization; instead they address problems where 
a (sub)gradient of the objective $f$ or a constraint function defining 
$\Omegab$ is not available to the optimization method.
Note that similar naming confusion has existed in 
non-smooth optimization, as evidenced by the introduction of \citeasnoun{IIASA1978}:
\begin{quote}
This workshop was held under the name \underline{Nondifferentiable 
Optimization}, but
it has been recognized that this is misleading, because it suggests 
`optimization without derivatives'.
\end{quote}

\subsection{Alternatives to derivative-free optimization methods}
\label{sec:intro:alts}
Derivative-free optimization methods are sometimes employed for 
convenience rather than by necessity. 
Since the decision to use a derivative-free method typically 
limits the performance -- in terms of accuracy, expense or problem 
size -- relative to what 
one might expect from gradient-based optimization methods, we 
first mention alternatives to using derivative-free methods. 

The design of derivative-free optimization methods is informed by the 
alternatives of algorithmic and numerical differentiation. 
For the former, the purpose 
seems clear: since the methods use only function values, they apply even 
in cases when one cannot produce a computer code for the function's derivative. 
Similarly, derivative-free optimization methods should be designed in order to 
outperform (typically measured in terms of the number of function evaluations) 
gradient-based optimization methods that employ numerical differentiation.

\subsubsection{Algorithmic differentiation}
\label{sec:intro:ad}

Algorithmic differentiation\footnote{\emph{Algorithmic differentiation} is 
sometimes 
referred to as \emph{automatic differentiation}, but we follow the preferred 
convention of \citeasnoun{Griewank2003}.} (AD) is a means of generating derivatives 
of mathematical functions \anlchange{that are expressed in} 
 computer code 
\cite{Griewank2003,griewank2008edp}. The forward mode of AD may be viewed as 
performing differentiation of elementary mathematical operations in each line
of source code 
by means of the chain rule, while the reverse mode may be seen as traversing 
the resulting computational graph in reverse order. 

Algorithmic differentiation has the benefit of automatically exploiting 
function structure, such as partial separability or other sparsity, 
and the 
corresponding ability of producing a derivative code whose computational cost is 
comparable to the cost of evaluating the function code itself.

AD has seen significant adoption and advances 
in the past decade \cite{AD2012advances}. Tools for algorithmic 
differentiation cover a growing set of compiled and interpreted languages, with 
an evolving list summarized on the 
community portal at 
\[\mbox{\textsf{http://www.autodiff.org}.}\] 
Progress has also been made on 
algorithmic differentiation of piecewise smooth functions, such as those 
with breakpoints resulting 
from absolute values or 
conditionals in a code; see, for example, \citeasnoun{Griewank2015}. The machine 
learning renaissance has also fuelled demand and interest in AD, driven in 
large part by the success of algorithmic differentiation in backpropagation 
\cite{2015arXiv150205767G}.

\subsubsection{Numerical differentiation}
\label{sec:intro:fd}

Another alternative to derivative-free methods is to estimate the derivative of
$f$ by numerical differentiation and then to 
use the estimates in a derivative-based 
method. This approach has the benefit that only zeroth-order information (\ie\
the
function value) is needed; however, depending on the derivative-based method 
used, the quality of the derivative estimate may be a limiting factor. 
Here we remark that for the finite-precision (or even fixed-precision) 
functions encountered in scientific applications, finite-difference estimates of 
derivatives may be sufficient for many purposes; see
\secref{orphaned-fd}.

When numerical derivative estimates are
used, the optimization method must 
tolerate inexactness in the derivatives. Such methods have been 
classically studied for both non-linear equations and unconstrained 
optimization; see, for example, the works of \citeasnoun{Powell1965}, 
\citeasnoun{Brown1971} and 
\citeasnoun{Mifflin1975} and the references therein. Numerical 
derivatives 
continue to be employed by recent methods (see \eg\ the works of 
\citeb{Cartis2012} and \citeb{BBN2018}).
Use in practice is typically determined by whether the limit on the derivative 
accuracy and the 
expense in terms of function evaluations are acceptable.

\subsection{Organization of the paper}
\label{sec:intro:organization} 
This paper is organized principally by problem class: unconstrained domain 
(Sections~\ref{sec:det_det} and \ref{sec:rand_det}), convex objective (\secref{convex}), structured objective
(\secref{structured}), stochastic optimization (\secref{stoch}) and 
constrained 
domain 
(\secref{constrained}). 

\secref{det_det} presents deterministic methods for solving
\eqref{eq:det_prob} when $\Omegab = \Reals^n$. The section is split between 
direct-search methods and model-based methods, although the lines between these 
are increasingly blurred; see, for example, 
\citeasnoun{Conn2013}, \citeasnoun{CRV2008}, \citeasnoun{Gramacy2015} and \citeasnoun{Gratton2016}. 
Direct-search methods are 
summarized in far greater detail by \citeasnoun{Kolda2003} and \citeasnoun{Kelleybook},
and in the more recent survey by \citeasnoun{Audet2014a}. Model-based methods that 
employ trust regions are given full treatment by \citeasnoun{Conn2009a}, and those 
that
employ stencils are detailed by \citeasnoun{Kel2011}.

In \secref{rand_det} we review randomized methods for solving
\eqref{eq:det_prob} when $\Omegab = \Reals^n$. These methods are 
often variants of the deterministic methods in \secref{det_det} but require 
additional notation to capture the resulting stochasticity; the analysis of 
these methods can also deviate significantly from their deterministic 
counterparts.

In \secref{convex} we 
discuss derivative-free methods intended 
primarily for convex optimization. We make this \anlchange{delineation} 
because such methods have distinct lines of analysis and can often 
solve considerably higher-dimensional problems than can general 
methods for non-convex derivative-free optimization. 

In \secref{structured} we survey methods that address particular structure 
in the objective $f$ in \eqref{eq:det_prob}. Examples of such structure 
include non-linear least-squares objectives, composite non-smooth objectives 
and partially separable objectives.

In \secref{stoch} we address derivative-free stochastic optimization, 
that is, when methods have access only to a stochastic realization of a
function in pursuit of solving \eqref{eq:stoch_prob}. This topic is 
increasingly intertwined with 
simulation optimization and Monte Carlo-based optimization; for these areas 
we refer to the surveys by 
\citeasnoun{HomemdeMello201456}, \citeasnoun{Fu2005}, \citeasnoun{Amaran2015} and 
\citeasnoun{Kim2014}.

\secref{constrained} presents methods for deterministic 
optimization problems with constraints (\ie\ $\Omegab \subset \Reals^n$).
Although many of these methods rely on 
the foundations laid in \secrefs{det_det}{rand_det}, we highlight particular 
difficulties associated with constrained derivative-free optimization. 

In \secref{extensions} we briefly highlight related problem areas (including 
global and multi-objective derivative-free optimization), 
methods and 
other implementation considerations.

\vspace{5pt}
\section{Deterministic methods for deterministic objectives} %
\label{sec:det_det}
We now address deterministic methods for solving \eqref{eq:det_prob}. We
discuss direct-search methods in \secref{det_det_DS}, model-based methods in
\secref{det_det_model} and other methods in \secref{det_det_others}. 
At a coarse level, direct-search methods use comparisons of function values to
directly determine \anlchange{candidate points}, whereas model-based methods use a surrogate of $f$
to determine \anlchange{candidate points}. Naturally, some hybrid methods incorporate
ideas from both model-based and direct-search methods and may not be so easily
categorized. 
An early survey of direct-search and
model-based methods is given in \citeasnoun{Powell1998a}.

\subsection{Direct-search methods} 
\label{sec:det_det_DS}

Although \citeasnoun{Hooke1961} are credited with originating the term `direct
search', there is no agreed-upon definition of what constitutes a
direct-search method. 
We follow the convention of \citeasnoun{Wright95}, wherein a direct-search method is a method that uses only function values and `does not
``in its heart'' develop an approximate gradient'.

We first discuss simplex methods, including the Nelder--Mead method -- perhaps
the most widely used direct-search method. We follow this discussion with a presentation
of directional direct-search methods; hybrid direct-search methods are discussed
in \secref{det_det_others}.
(The global direct-search method \codes{DIRECT} is discussed in \secref{global}.)

\subsubsection{Simplex methods}
\label{sec:det_det_NM}

Simplex methods (not to be confused with Dantzig's simplex method for linear
programming) move and manipulate a collection of $n+1$ affinely independent
points (\ie\ the vertices of a simplex in $\Reals^n$) when solving \eqref{eq:det_prob}. 
The method of \citeasnoun{Spendley1962} involves either taking the point in the simplex with
the largest function value and reflecting it through the hyperplane defined by
the remaining $n$ points or moving the $n$ worst points toward the best vertex
of the simplex. In this manner, the geometry of all simplices remains the same
as that of the starting simplex. (That is, all simplices are similar in the geometric sense.) 

\begin{figure} %
\centering
\begin{tikzpicture}[line join = round, line cap = round]
\pgfmathsetmacro{\factori}{sqrt(2)};
\pgfmathsetmacro{\factor}{1/sqrt(2)};
\pgfmathsetmacro{\shift}{2.4};

\foreach \xg in {0,1,2,3,4}{
 \coordinate [] (A4b) at ( 0+\xg*\shift, 0, 0*\factor);
 \coordinate [] (B4b) at (-2+\xg*\shift, 0, 0*\factor);
 \coordinate [] (C4b) at (-1+\xg*\shift, 0, 2*\factor);
 \coordinate [] (D4b) at (-2+\xg*\shift,-2,-.33*\factori);
 \draw[-,dotted,blue,opacity=.95] (B4b)--(D4b)--(A4b);
 \draw[-,dotted,blue,opacity=.95] (C4b)--(B4b)--(D4b)--cycle;
 \draw[-,dotted,blue,opacity=.95] (D4b)--(A4b)--(C4b)--cycle;
 \draw[-,dotted,blue,opacity=.95] (A4b)--(C4b)--(B4b)--cycle; 
}

\coordinate 
  (A1) at ( 0+0*\shift, 0, 0*\factor);
\coordinate [label=above:$\xb_{(1)}$] (B1) at (-2+0*\shift, 0, 0*\factor);
\coordinate 
  (C1) at (-1+0*\shift, 0, 2*\factor);
\coordinate [label=-20:$\xb_{(n+1)}$] (D1) at (-2+0*\shift,-2,-.33*\factori);
\draw[-, fill=red!30, opacity=.5  ] (B1)--(D1)--(A1);
\draw[-, fill=red!30, opacity=.5  ] (C1)--(B1)--(D1)--cycle;
\draw[-, fill=purple!70, opacity=.5] (D1)--(A1)--(C1)--cycle;
\draw[-, fill=green!30, opacity=.5] (A1)--(C1)--(B1)--cycle; 

\coordinate 
  (A2) at ( 0+1*\shift, 0, 0*\factor);
\coordinate 
  (B2) at (-2+1*\shift, 0, 0*\factor);
\coordinate 
  (C2) at (-1+1*\shift, 0, 2*\factor);
\coordinate [label=above:$\xb^{\rm new}$] 
  (D2) at ( 0+1*\shift, 2, 1.67*\factor);
\draw[dashed] (A2)--(B2);
\draw[-, fill=red!30, opacity=.5  ] (B2)--(D2)--(A2);
\draw[-, fill=green!30, opacity=.5] (A2)--(C2)--(B2); 
\draw[-, fill=red!30, opacity=.5  ] (C2)--(B2)--(D2)--cycle;
\draw[-,fill=purple!70, opacity=.5] (D2)--(A2)--(C2)--cycle;

\coordinate (A5) at ( 0+2*\shift, 0, 0*\factor);
\coordinate (B5) at (-2+2*\shift, 0, 0*\factor);
\coordinate (C5) at (-1+2*\shift, 0, 2*\factor);
\coordinate [label=above:$\xb^{\rm new}$] (D5) at ( 0.5+2*\shift, 3, 1.67*\factor);
\draw[dashed] (A5)--(B5);
\draw[-, fill=red!30, opacity=.5  ] (B5)--(D5)--(A5);
\draw[-, fill=green!30, opacity=.5] (A5)--(C5)--(B5); 
\draw[-, fill=red!30, opacity=.5  ] (C5)--(B5)--(D5)--cycle;
\draw[-,fill=purple!70, opacity=.5] (D5)--(A5)--(C5)--cycle;

\coordinate 
  (A3) at ( 0  +3*\shift, 0,  0*\factor);
\coordinate 
  (B3) at (-2 +3*\shift, 0,  0*\factor);
\coordinate 
  (C3) at (-1  +3*\shift, 0,  2*\factor);
\coordinate (D3) at (-1.5+3*\shift,-1,  
.167*\factori);
\coordinate [label={[label distance=-0.1cm]right:$\xb^{\rm new}$}] (D3a) at (-1.0+2.95*\shift,-1.4, .167*\factori);
\draw[line width = 0.7,->,>=stealth,color=black] (D3a) -- (D3);

\draw[-, fill=red!30, opacity=.5  ] (B3)--(D3)--(A3);
\draw[-, fill=red!30, opacity=.5  ] (C3)--(B3)--(D3)--cycle;
\draw[-,fill=purple!70, opacity=.5] (D3)--(A3)--(C3)--cycle;
\draw[-, fill=green!30, opacity=.5] (A3)--(C3)--(B3)--cycle; 

\coordinate (A4) at ( -1+4*\shift, 0, 
0*\factor);
\coordinate [label={[label distance=-0.25cm]above right:$\xb_2^{\rm new}$}] (A4a) at ( -1+4.1*\shift, 0.4, 0*\factor);
\coordinate 
(B4) at (-2+4*\shift, 0, 0*\factor);
\coordinate (C4) at (-1.5+4*\shift, 0,
1*\factor);
\coordinate [label={[label distance=-0.1cm]right:$\xb_1^{\rm new}$}] (C4a) at (-1.5+4.35*\shift, -0.9, 1*\factor);
\coordinate (D4) at 
(-2+4*\shift,-1,-.166*\factori);
\coordinate [label={[label distance=-0.1cm]right:$\xb_3^{\rm new}$}] (D4a) at (-2+4.35*\shift,-1.8,-.166*\factori);
\draw[-, fill=red!30, opacity=.5  ] (B4)--(D4)--(A4);
\draw[-, fill=red!30, opacity=.5  ] (C4)--(B4)--(D4)--cycle;
\draw[-,fill=purple!70, opacity=.5] (D4)--(A4)--(C4)--cycle;
\draw[-, fill=green!30, opacity=.5] (A4)--(C4)--(B4)--cycle; 
\draw[line width = 0.7,->, >=stealth, color=black] (A4a) -- (A4);
\draw[line width = 0.7,->, >=stealth, color=black] (C4a) -- (C4);
\draw[line width = 0.7,->, >=stealth, color=black] (D4a) -- (D4);

%

\end{tikzpicture} %
 \caption{Primary Nelder--Mead simplex operations: original simplex, 
reflection, expansion, inner contraction, and shrink. 
\label{fig:nelder-mead}}
\end{figure}
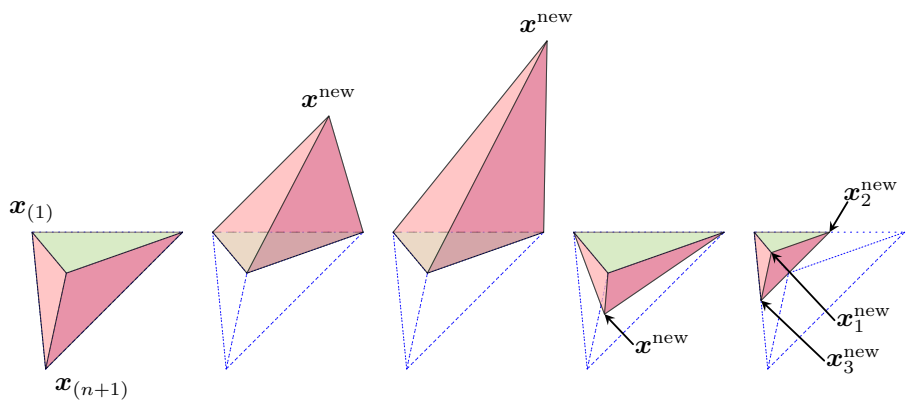

\citeasnoun{NelderMead} extend the possible simplex operations, as shown in
\fo{Figure~\ref{fig:nelder-mead}} by allowing the `expansion' and `contraction'
operations in addition to the `reflection' and `shrink' operations of
\citeasnoun{Spendley1962}. These operations enable the Nelder--Mead simplex method
to distort the simplex in order to account for possible curvature
present in the objective function.

\citeasnoun{NelderMead} propose stopping further function evaluations when the
standard error of the function values at the simplex vertices is small. Others,
\citeasnoun{WoodsPhD} for example, propose stopping when the size of the
simplex's longest side incident to the \anlchange{best simplex vertex} is small.

Nelder--Mead is an incredibly popular method, in no small part due to its
inclusion in {\em Numerical Recipes} \cite{pres07},
which has been cited over \anlchange{$125\,000$} times and no doubt used many times more.
The method (as implemented by \citeb{Lagarias2012})
is also the algorithm underlying \codes{fminsearch} in \framework{MATLAB}. 
Benchmarking studies highlight Nelder--Mead performance in practice
\cite{JJMSMW09,Rios2013}.

The method's
popularity from its inception was not diminished by the lack of theoretical
results proving its ability to identify stationary points. \citeasnoun{WoodsPhD}
presents a non-convex, two-dimensional function where Nelder--Mead converges to a
non-stationary point (where the function's Hessian is singular). 
Furthermore,
\citeasnoun{McKinnon1998} presents a class of thrice-continuously differentiable,
strictly convex functions on $\Reals^2$ where the Nelder--Mead simplex fails to
converge to the lone stationary point. The only operation that Nelder--Mead
performs on this relatively routine function is repeated `inner contraction' of
the initial simplex.

Researchers have continued to develop convergence results for modified or
limited versions of Nelder--Mead.
\citeasnoun{Kelley1999} addresses Nelder--Mead's theoretical deficiencies by
restarting the method when the objective decrease on consecutive iterations is not larger
than a multiple of the 
simplex gradient norm. Such restarts do not ensure that Nelder--Mead will 
converge: \citeasnoun{Kelley1999} shows an example of such behaviour.
\citeasnoun{Price2002} embed Nelder--Mead in a different (convergent) algorithm using
positive spanning sets. \citeasnoun{nazareth2002gilding} propose a clever, though
perhaps superfluous, variant that connects Nelder--Mead to golden-section
search.

\citeasnoun{JCL98} show that Nelder--Mead (with appropriately chosen reflection and
expansion coefficients) converges to the global minimizer of strictly
convex functions when $n=1$. 
\citeasnoun{Gao2010} show that the contraction and expansion steps of Nelder--Mead
satisfy a descent condition on uniformly convex functions. 
\citeasnoun{Lagarias2012} show that a restricted version of the Nelder--Mead
method -- one that does not allow an expansion step -- can converge to
minimizers of any twice-continuously differentiable function with a
positive-definite Hessian and bounded level sets. (Note that the class of
functions from \citeasnoun{McKinnon1998} have singular Hessians at only one point -- their
minimizers -- and not at the point to which the \anlchange{simplex vertices} are converging.) 

The simplex method of \citeasnoun{Rykov1980} includes ideas from model-based
methods. Rykov varies the number of reflected vertices
 from iteration to iteration, following one of three rules that depend on
the function value at the simplex centroid $\xb_c$. Rykov considers both
evaluating $f$ at the centroid and approximating $f$ at the centroid using the
values of $f$ at the vertex. The non-reflected vertices are also moved in
parallel with the reflected subset of vertices. In general, the number of
reflected vertices is chosen so
that $\xb_c$ moves in a direction closest to $-\nabla f(\xb_c)$.
This, along with a test of sufficient decrease in $f$, ensures convergence of the
modified simplex method to a minimizer of convex, continuously differentiable
functions with bounded level sets and Lipschitz-bounded gradients. (The
sufficient-decrease condition is also shown to be efficient for the classical
Nelder--Mead algorithm.)

\citeasnoun{Tseng99} proposes a modified simplex method that keeps the $b_k$ best
simplex vertices on a given iteration $k$ and uses them to reflect the
remaining vertices. Their method prescribes that `the rays emanating from the
reflected vertices toward the $b_k$ best vertices should contain, in their
convex hull, the rays emanating from a weighted centroid of the $b_k$ best
vertices toward the to-be-reflected vertices'. Their method also includes a
\emph{fortified descent} condition that is stronger than common sufficient-decrease
conditions. If $f$ is continuously differentiable and bounded below and $b_k$
is fixed for all iterations, \citeasnoun{Tseng99} prove that every cluster point of
the sequence of \anlchange{candidate points generated by} their method is a stationary point.

\anlchange{\citeasnoun{Brmen2006} propose a convergent version of a simplex method  that does not
require a sufficient descent condition be satisfied. Instead, they ensure that
evaluated points lie on a grid of points, and they show that this grid will
be refined as the method proceeds.}

\subsubsection{Directional direct-search methods} \label{sec:det_det_DDS}
Broadly speaking, each iteration of a directional direct-search (\framework{DDS}) method
generates a finite set of points near the current \anlchange{point $\xb_k$}; these \emph{poll
points} are generated by taking \anlchange{$\xb_k$}  and adding terms of
the form $\alpha_k \db$, where $\alpha_k$ is a positive step size and $\db$ is an
element from a finite set of directions $\Db_k$. \citeasnoun{Kolda2003} propose the term
\emph{generating set search methods} to encapsulate this class of
methods.\footnote{The term generating set arises from a need to generate a cone
from the nearly active constraint normals when $\Omegab$ is defined by linear constraints.} 
The objective function $f$ is then evaluated at all or some of the poll points,
and \anlchange{$\xb_{k+1}$} is selected to be some poll point that produces
a (sufficient) decrease in the objective and the step size is possibly increased. If
no poll point provides a sufficient decrease, \anlchange{$\xb_{k+1}$ is set to $\xb_k$} and the
step size is decreased. In either case, the set of directions $\Db_k$ can (but
need not) be modified \anlchange{to obtain $\Db_{k+1}$}. 

\begin{algorithm}[tb]
 \SetKw{break}{break}
 Initialize $\xb^+\gets \xb$\\
 \For{$\pb_i \in \Pb $}
 {
 Evaluate $f(\pb_i)$\\
 \If{$f(\pb_i) - f(\xb)$ acceptable \label{line:descent}} 
 {
 $\xb^+ \gets \pb_i$\\
 optional \break \label{line:optional_break}
 }
 }
 \caption{$\xb^+ = \mathtt{test\_descent}(f,\xb,\Pb)$ \label{alg:descent}}
\end{algorithm}

\begin{algorithm}[tb]
 Set parameters $0 < \gammad < 1 \le \gammai$\\
 Choose initial point $\xb_0$ and step size $\alpha_0 >0$\\
 \For{$k=0,1,2,\ldots$}
 {
 Choose and order a finite set $\Yb_k\subset \Reals^n$ \label{line:search_step} \tcp*{(search step)}
 $\xb_{k}^+ \gets \mathtt{test\_descent}(f,\xb_k,\Yb_k)$\\
 \If{$\xb_k^+ = \xb_k $}
 		{
 Choose and order poll directions $\Db_k\subset\Reals^n$ \label{line:poll_step} \tcp*{(poll step)} 
 $\xb_k^+ \gets \mathtt{test\_descent}(f, \xb_k, \{\xb_k + \alpha_k \db_i: \db_i\in \Db_k\})$ 
 		}
 \eIf{$\xb_k^+ = \xb_k $}
 		{
 $\alpha_{k+1}\gets \gammai\alpha_k$\\
 		}
 		{
 $\alpha_{k+1} \gets \gammad \alpha_k$\\
 		}
 $\xb_{k+1}\gets \xb_k^+$
 	}
 \caption{Directional direct-search method\label{alg:ds}}
\end{algorithm}

A general \framework{DDS} method is provided in \algref{ds}, which includes a
\emph{search step} where $f$ is evaluated at any finite set of points 
$\Yb_k$, including
$\Yb_k=\emptyset$. The search step allows one to (potentially) improve 
the performance of \algref{ds}. For example, \anlchange{points} could be randomly
sampled during the search step from the domain in the hope of finding a better
local minimum, or a person running the algorithm may have problem-specific
knowledge that can generate candidate points given the \anlchange{observed history of evaluated points and their} function values. 
While the search step allows for this insertion of such heuristics,
rigorous convergence results are driven by the more disciplined poll step.
When testing for objective decrease in \algref{descent}, one can stop
evaluating points in $\Pb$ (\lineref{optional_break}) as soon as the first
point is identified where there is (sufficient) decrease in $f$. In this case, the polling
(or search) step is considered \emph{opportunistic}.

\framework{DDS} methods are largely distinguished by how they generate the set of poll
directions $\Db_k$ at \lineref{poll_step} of \algref{ds}. 
Perhaps the first approach is \emph{coordinate search}, in which the poll
directions are defined as $\Db_k = \{\pm \eb_i:i=1,2,\ldots,n\}$,
where $\eb_i$ denotes \anlchange{the $i$th elementary basis vector (\ie\ } column $i$ of the identity matrix in $n$ dimensions\anlchange{)}.
The first known description of coordinate search appears in the work of \citeasnoun{Fermi1952}
where the smallest positive integer
$l$ is sought such that $f(\xb_k + l\alpha \eb_1/2) > f(\xb_k + (l-1)\alpha
\eb_1/2)$. If an increase in $f$ is observed at $\eb_1/2$ then
$-\eb_1/2$ is considered. After such an integer $l$ is identified for
the first coordinate direction, $\xb_k$ is updated to $\xb_k \pm l \eb_1/2$ and
the second coordinate direction is
considered. If $\xb_k$ is unchanged after cycling through all coordinate
directions, then the method is repeated but with $\pm \anlchange{\eb}_i/2$ replaced with $\pm \anlchange{\eb}_i/16$, 
terminating when no improvement is observed for this smaller $\alpha$. In terms
of \algref{ds} the search set $\Yb_k=\emptyset$ at \lineref{search_step}, and the
descent test at \lineref{descent} of \algref{descent} merely
tests for simple decrease, that is, $f(\pb_i)-f(\xb)<0$. Other versions of
acceptability in \lineref{descent} of \algref{descent} are employed by methods
discussed later.

Proofs that \framework{DDS} methods converge first appeared in the works of \citeasnoun{Cea1971} and
\citeasnoun{Yu1979positive}, \anlchange{al}though both require the sequence of step-size parameters
to be non-increasing. 
\citeasnoun{Lewis2000} attribute the first global
convergence proof for coordinate search to \citeasnoun[p.~43]{Polak1971computational}. In turn, Polak cites the `method of local variation'
of \citeasnoun{Banichuk1966}; although \citeasnoun{Banichuk1966} do develop parts of a
convergence proof, they state in Remark~1 that `the question of the strict
formulation of the general sufficient conditions for convergence of the
algorithm to a minimum remains open'. 

Typical convergence results for \framework{DDS} require that the set $\Db_k$ is a 
\emph{positive spanning set} (PSS) for the domain $\Omegab$; that is, any point 
$\xb \in
\Omegab$ can be written as 
\[
\xb = \sum_{i = 1}^{| \Db_k |} \lambda_i \db_i,
\]
where $\db_i \in \Db_k$ and $\lambda_i \ge 0$ for all $i$. Some of the first
discussions of properties of positive spanning sets were presented by \citeasnoun{Davis1954} and
\citeasnoun{McKinney1962}, but recent treatments have also appeared in
\citeasnoun{regis2016properties}.
In addition to requiring positive spanning sets during the poll step, earlier \framework{DDS}
convergence results depended on $f$ being continuously differentiable. 
When $f$ is non-smooth, no descent direction is guaranteed for these early
\framework{DDS} methods, even when the step size is arbitrarily small.
See, for example, the modification of the Dennis--Woods \cite{Dennis1987}
function by \citeasnoun[Figure~6.2]{Kolda2003} and a discussion of why
coordinate-search methods (for example) will not move when started at a point of
non-differentiability; moreover, when started at differentiable points, coordinate-search
methods tend to converge to a point that is not (Clarke) stationary.

The pattern-search method of \citeasnoun{Torczon1991} revived 
interest in direct-search methods. The method therein contains ideas from 
both
DDS and simplex methods. 
Given a simplex defined by $\xb_k, \yb_1, \ldots, \yb_n$ (where $\xb_k$ is the 
simplex vertex with smallest function value), the polling directions are given 
by $\Db_k = \{
\yb_i - \xb_k : i = 1,\ldots,n \}$.
If a decrease is observed at the best poll point in $\xb_k+\Db_k$, the 
simplex is set to either 
$\xb_k \bigcup \xb_k + \Db_k $ or some expansion thereof. If no improvement is 
found
during the poll step, the simplex is contracted. 
\citeasnoun{Torczon1991} shows that if $f$ is continuous on the level set of $\xb_0$ and this
level set is compact, then a subsequence of \anlchange{$\{\xb_k\}$} converges to
a stationary point of $f$, a point where $f$ is non-differentiable, or a point
where $f$ is not continuously differentiable. 

A generalization of pattern-search methods is the class of \emph{generalized
pattern-search \framework{(GPS)} methods}. Early \framework{GPS} methods did not
allow for a search step; the search-poll paradigm was introduced by
\citeasnoun{Booker1999}. 
\framework{GPS} methods are characterized by
fixing a positive spanning set $\Db$ and selecting $\Db_k\subseteq \Db$ during
the poll step at \lineref{poll_step} on each iteration of \algref{ds}.
\citeasnoun{Torczon1997} assumes that the test for decrease in \lineref{descent} in
\algref{descent} is simple decrease, that is, that $f(\pb_i) < f(\xb)$. Early
analysis of \framework{GPS} methods using simple decrease required 
the step size $\alpha_k$ to remain rational \cite{CAJD02,Torczon1997}. \citeasnoun{Audet2004b} shows that
such an assumption is necessary by constructing
small-dimensional examples where \framework{GPS} methods do not converge if $\alpha_k$ is
irrational. Works below show that
if a sufficient (instead of simple) decrease is ensured, $\alpha_k$
can take irrational values.

A refinement of the analysis of \framework{GPS} methods was made by \citeasnoun{Dolan2003}, which
shows that when $\nabla f$ is Lipschitz-continuous, the step-size parameter
$\alpha_k$ scales linearly with $\|\nabla f(\xb_k)\|$. Therefore $\alpha_k$ can
be considered a reliable measure of first-order stationarity and justifies the
traditional approach of stopping a \framework{GPS} method when $\alpha_k$ is small.
Second-order convergence analyses of \framework{GPS} methods have also been considered.
\citeasnoun{Abramson2005} shows that, when applied to \anlchange{a} twice-continuously
differentiable $f$,
\anlchange{a \framework{GPS} method that infinitely often has $\Db_k$ 
include a fixed orthonormal basis and its negative 
will have a limit point} satisfying a `pseudo-second-order' stationarity
condition. Building off the use of curvature information in \citeasnoun{FriSte07}, 
\citeasnoun{Abramson2013} show that a modification of the \framework{GPS}
framework that constructs approximate Hessians of $f$ will converge to points
that are second-order stationary provided that certain conditions on the
Hessian approximation hold (and a fixed orthonormal basis and its negative are
in $\Db_k$ infinitely often). 

In general, first-order convergence results (there exists a limit point $\xb_*$
of \anlchange{$\{\xb_k\}$ generated by a} \framework{GPS} method such that $\nabla f(\xb_*) = \zerob$) for 
\framework{GPS} methods can be demonstrated
when $f$ is continuously differentiable.
For general Lipschitz-continuous (but non-smooth) functions $f$, however,
one can only demonstrate that on a particular subsequence $\cK$,
satisfying
$\{\xb_k\}_{k\in \cK}\to \xb_*$, 
for each $\db$ that appears infinitely many times in $\{\Db_k\}_{k\in \cK}$,
it holds that
$f'(\xb_*;\db)\geq 0$; that is, the directional derivative at $\xb_*$ in the
direction $\db$ is non-negative. 

The flexibility of \framework{GPS} methods inspired various extensions.
\citeasnoun{Abramson2004gps} consider adapting \framework{GPS} to utilize derivative 
information when it is available in order to reduce the number of points
evaluated during the poll step.
\citeasnoun{MAAbramson_etal_2008} and \citeasnoun{Frimannslund2011} re-use previous function
evaluations in order to determine the next set of directions. 
\citeasnoun{ACLV2007} consider re-using previous function evaluations 
to compute simplex gradients; they show that the information obtained from
simplex gradients can be used to reorder the poll points $\Pb$ in
\algref{descent}. A similar use of simplex gradients in the non-smooth
setting is considered by \citeasnoun{ACLVJD2008}. 
\citeasnoun{Torczon2003b} discuss modifications to \algref{ds} that allow for
increased efficiency when concurrent, asynchronous evaluations of $f$ are
possible; an implementation of the method of \citeasnoun{Torczon2003b} is presented
by \citeasnoun{Gray2006}. 

The early analysis of \citeasnoun[Section~7]{Torczon1991} of pattern-search 
methods when $f$ is non-smooth carries over to \framework{GPS}
methods as well; such methods may converge to a non-stationary point. This
motivated a further generalization of \framework{GPS} methods, \emph{mesh adaptive direct
search} (\codes{MADS}) methods \cite{Audet06mads,Abramson2006}. Inspired by
\citeasnoun{Coope2000b}, \codes{MADS} methods augment \framework{GPS} methods by incorporating a mesh
parametrized by a
\emph{mesh parameter} $\beta^m_k>0$. 
 In the $k$th iteration,
 given the fixed PSS $\Db$ and the mesh parameter $\beta^m_k$, the \codes{MADS} mesh around
 the current point $\xb_k$ is
 \[
\mathcal{M}_k = \bigcup_{\xb \in \Sb_k} 
\Biggl\{\xb + \beta^m_k \sum_{j=1}^{| \Db |} \lambda_j \db_j : \db_j \in \Db,
\lambda_j \in \mathbb{N}\bigcup\{ 0 \}\Biggr\},
\]
 where $\Sb_k$ is the set of points at which $f$ has been evaluated prior to 
the $k$th iteration of the method.

\codes{MADS} methods 
additionally define
a frame
\[
\mathcal{F}_k = \{\xb_k + \beta^m_k \db^f: \db^f\in\Db^f_k\},
\]
where $\Db^f_k$ is a finite set of directions, each of which is expressible as
\[
\db^f = \sum_{j=1}^{| \Db |} \lambda_j\db_j,
\]
with each $\lambda_j\in\mathbb{N}\bigcup\{0\}$ and $\db_j \in \Db$.
Additionally, \codes{MADS} methods define
a \emph{frame parameter} $\beta^f_k$
and require that each $\db^f\in\Db^f_k$ satisfies
$\beta^m_k\|\db^f\|\leq\beta^f_k\max\{\|\db\|:\db\in\Db\}$. 
Observe that in each iteration, $\mathcal{F}_k\subsetneq \mathcal{M}_k$.
Note that the mesh is never explicitly constructed nor
stored over the domain.
Rather, points are evaluated only at
what \emph{would be} nodes of some implicitly defined mesh
via the frame. 

In the poll step of \algref{ds}, the set of poll directions $\Db_k$ 
is chosen as $\{\yb-\xb_k:\yb\in\mathcal{F}_k\}$. 
The role of the step-size parameter $\alpha_k$ in \algref{ds}
is completely replaced by the behaviour of $\beta^f_k,\beta^m_k$. 
If there is no improvement
at a candidate solution during the poll step, $\beta^m_k$ is
decreased, resulting in a finer mesh; likewise $\beta^f_k$ is decreased,
resulting in a finer local mesh around $\xb_k$. 
\codes{MADS} intentionally allows the parameters
$\beta^m_k$ and $\beta^f_k$ to be decreased at different rates; 
roughly speaking, 
by driving $\beta^m_k$ to zero faster than $\beta^f_k$ is driven to zero, 
and by choosing the sequence $\{\Db^f_k\}$ to satisfy certain conditions, 
the directions in $\mathcal{F}_k$ become asymptotically dense around limit points of $\xb_k$. 
That is, it is possible to decrease $\beta^m_k,\beta^f_k$ at rates such that 
poll
directions will be arbitrarily close to any direction. 
This ensures that the Clarke directional derivative is non-negative in all
directions around any limit point of the sequence of \anlchange{$\xb_k$ generated by \codes{MADS}}; that is,
\begin{equation}
f_C'(\xb_*;\db)\geq 0 \quad \mbox{for all directions } \db,
\label{eq:uncon_Clarke_stationary}
\end{equation}
with an analogous result also holding for constrained problems, with
\eqref{eq:uncon_Clarke_stationary} reduced to all \emph{feasible} directions
$\db$. (\framework{DDS} methods for constrained optimization will be discussed in \secref{constrained}.)
This powerful result highlights the ability of directional direct-search
methods to address non-differentiable functions $f$.

\codes{MADS} does not prescribe any one approach for adjusting $\beta^m_k,\beta^f_k$ so that the poll
directions are dense, but \citeasnoun{Audet06mads} demonstrate an approach where randomized
directions are completed to be a PSS and $\beta^f_k$ either is $n\sqrt{\beta^m_k}$ or
$\sqrt{\beta^m_k}$ results in a asymptotically dense poll directions for any
\anlchange{convergent subsequence of $\{\xb_k\}$}.
\codes{MADS} does not require a sufficient-decrease condition.

Recent advances to \codes{MADS}-based algorithms have focused on reducing the number of 
function evaluations required in practice by adaptively reducing the number of 
poll points queried; see, for example, \citeasnoun{Audet2014} and \citeasnoun{Alarie2018}. Smoothing-based extensions to noisy deterministic problems include \citeasnoun{Audet2018}.
\citeasnoun{LVAC12} show that \codes{MADS} methods converge to local minima even for
a limited class of \emph{discontinuous} functions that satisfy some assumptions
concerning the behaviour of the disconnected regions of the epigraph at limit points.

\paragraph{Worst-case complexity analysis.}
Throughout this survey, when discussing classes of methods, we will refer to their
worst-case complexity (WCC). 
Generally speaking, WCC refers to an upper bound on the number of function evaluations 
$N_{\epsilon}$ 
required to attain
an $\epsilon$-accurate solution to a problem drawn from a problem class.
Correspondingly, the definition of $\epsilon$-accurate varies between different problem classes. 
For instance, and of particular immediate importance, 
if an objective function is assumed Lipschitz-continuously differentiable (which we denote by $f\in\cLC^1$),
then an appropriate notion of first-order $\epsilon$-accuracy is
\begin{equation}
\label{eq:nonconvex_stationarity}
\|\nabla f(\xb_k)\|\leq\epsilon.
\end{equation}
That is, the WCC of a method applied to the class $\cLC^1$ is characterized by $N_{\epsilon}$,
an upper bound on the number of function evaluations the method requires before  \anlchange{\eqref{eq:nonconvex_stationarity} is satisfied} for \emph{any} $f\in\cLC^1$. 
 Similarly, we can define a notion of second-order $\epsilon$-accuracy as
 \begin{equation}
 \label{eq:2ndnonconvex_stationarity}
 \max\{\|\nabla f(\xb_k)\|,-\lambda_k\}\leq\epsilon,
 \end{equation}
where $\lambda_k$ denotes the minimum eigenvalue of $\nabla^2 f(\xb_k)$. 

Note that WCCs can only be derived for methods for which convergence results have been established. 
Indeed, in the problem class $\cLC^1$, first-order convergence results canonically have the form
\begin{equation}
\label{eq:first_order_convergence}
\lim_{k\to\infty} \|\nabla f(\xb_k)\| = 0.
\end{equation}
The convergence in \eqref{eq:first_order_convergence} automatically implies the weaker lim-inf-type result
\begin{equation}
\label{eq:liminf_first_order_convergence}
\liminf_{k\to\infty} \|\nabla f(\xb_k)\| = 0,
\end{equation}
from which it is clear that for any $\epsilon>0$, there must exist finite $N_{\epsilon}$ so that 
\eqref{eq:nonconvex_stationarity} holds. 
In fact, in many works, demonstrating a result of the form \eqref{eq:liminf_first_order_convergence} is a stepping
stone to proving a result of the form \eqref{eq:first_order_convergence}.
Likewise, demonstrating a second-order WCC of the form \eqref{eq:2ndnonconvex_stationarity} depends on showing
\begin{equation}
\label{eq:second_order_convergence}
\lim_{k\to\infty}\max\{\|\nabla f(\xb_k)\|,-\lambda_k\} = 0,
\end{equation}
which guarantees the weaker lim-inf-type result
\begin{equation}
\label{eq:liminf_second_order_convergence}
\liminf_{k\to\infty}\max\{\|\nabla f(\xb_k)\|,-\lambda_k\} = 0. 
\end{equation}

Proofs of convergence for \framework{DDS} methods applied to functions $f\in\cLC^1$ 
often rely on 
a (sub)sequence of positive spanning sets $\{\Db_k\}$ satisfying
\begin{equation}\label{eq:cosine}
 \cm(\Db_k) \defined
 \min_{\vb\in\Reals^n\setminus\{\zerob\}}\max_{\db\in\Db_k}\ffrac{\db^\top 
\vb}{\|\db\|\|\vb\|} 
 \geq \kappa > 0,
\end{equation}
where $\cm(\cdot)$ is the \emph{cosine measure} of a set. 
Under Assumption~\eqref{eq:cosine},
\citeasnoun{Vicente2013} obtains a WCC of type \eqref{eq:nonconvex_stationarity} 
for a method in the 
\algref{ds} framework. In that work, it is assumed that $\Yb_k=\emptyset$ at every search step.
Moreover, \emph{sufficient} decrease is tested at \lineref{descent} of \algref{descent};
in particular, \citeasnoun{Vicente2013} checks in this line whether 
$f(\pb_i) < f(\xb) - c\alpha_k^2$ for some $c>0$,
where $\alpha_k$ is the current step size in \algref{ds}.
Under these assumptions, \citeasnoun{Vicente2013} demonstrates a WCC in 
$\bigo{\epsilon^{-2}}$. 
Throughout this survey, we will refer to Table~\ref{table:rates} for more details 
concerning specific WCCs.
In general, though, we will often summarize WCCs in terms of their $\epsilon$-dependence,
as this provides an asymptotic characterization of a method's complexity in terms of the accuracy
to which one wishes to solve a problem. 

When $f\in\cLC^2$, 
work by 
\citeasnoun{Gratton2016} essentially
augments the \framework{DDS} method analysed by \citeasnoun{Vicente2013},
but forms an
approximate Hessian via central differences from function evaluations obtained (for free) 
by using a particular choice of $\Db_k$. 
\citeasnoun{Gratton2016} then demonstrate that this augmentation of \algref{ds} has
a subsequence that converges to a second-order stationary point. That is,
they prove a convergence result of the form \eqref{eq:liminf_second_order_convergence} 
and demonstrate
a WCC result of type \eqref{eq:2ndnonconvex_stationarity} 
in $\bigo{\epsilon^{-3}}$ (see Table~\ref{table:rates}).

We are unaware of WCC results for \codes{MADS} methods; this situation may be unsurprising since
\codes{MADS} methods are motivated by non-smooth problems, which depend on the generation of 
a countably infinite number of poll directions.
However, WCC results are not necessarily impossible to obtain in \emph{structured} non-smooth cases,
which we discuss in \secref{structured}. 
We will discuss a special case where \emph{smoothing functions} of
a non-smooth function are assumed to be available in \secref{cno_noncon}.

\subsection{Model-based methods}
\label{sec:det_det_model}

In the context of derivative-free optimization, model-based methods are methods 
whose updates are based primarily on the predictions of a model that serves as 
a surrogate of the objective function or of a related merit function. 
We begin with basic properties and construction of popular models; readers 
interested in algorithmic frameworks such as trust-region methods and 
implicit filtering can proceed to \secref{det_det_TRM}. Throughout 
this section, we assume that models are intended as a surrogate for the 
function $f$; in future sections, these models will be extended to capture 
functions arising, for example, as constraints or separable components.
The methods in this section assume some smoothness in $f$ and therefore operate 
with smooth models; in \secref{structured}, we examine model-based methods that exploit 
knowledge of non-smoothness.

\subsubsection{Quality of smooth model approximation} 
\label{sec:det_det_modelquality}

A natural first indicator of the quality of a model used for optimization is 
the degree to which the model locally approximates the function $f$ and its 
derivatives. To say anything about the quality of such approximation, one must make an 
assumption about the smoothness of both the model and function. For the 
moment, we leave this assumption implicit, but it will be formalized in 
subsequent sections.

A function $m: \Reals^n\to \Reals$ is said to be a 
\anlchange{$\kappab$}-fully linear
 model of $f$ on $\B{\xb}{\Delta}=\{\yb: \|\xb-\yb\|\leq \Delta\}$ 
if
\begin{subequations}
\label{eq:fullylinear}
 \begin{alignat}{3}
 \label{eq:fl:f}
 |f(\xb+\sba) - m(\xb+\sba)| & \leq \kappaef \Delta^2,
 &	 \quad &\gforall\ \sba \in \B{\zerob}{\Delta},\\*
 \label{eq:fl:g}
 \|\nabla f(\xb+\sba) - \nabla m(\xb+\sba)\| & \leq \kappaeg\Delta,
 &	 \quad &\gforall\ \sba \in \B{\zerob}{\Delta},
 \end{alignat}
\end{subequations} 
for $\kappab = (\kappaef, \kappaeg)$.
Similarly, for $\kappab = (\kappaef, \kappaeg, \kappaeH)$, $m$ is said to 
be a \anlchange{$\kappab$}-fully quadratic model of $f$ on $\B{\xb}{\Delta}$ if 
\begin{subequations}
\label{eq:fullyquadratic}
 \begin{alignat}{3}
 \label{eq:fq:f}
 |f(\xb+\sba) - m(\xb+\sba)| & \leq \kappaef \Delta^3,
	 & \quad &\gforall\ \sba \in \B{\zerob}{\Delta},\\*
 \label{eq:fq:g}
 \|\nabla f(\xb+\sba) - \nabla m(\xb+\sba)\| & \leq 
	 \kappaeg\Delta^2,
	 &\quad &\gforall\ \sba \in \B{\zerob}{\Delta},\\*
 \label{eq:fq:H}
 \|\nabla^2 f(\xb+\sba) - \nabla^2 m(\xb+\sba)\| & \leq 
	 \kappaeH\Delta,
	 &\quad &\gforall\ \sba \in \B{\zerob}{\Delta}.	 
 \end{alignat}
\end{subequations} 
Extensions to higher-degree approximations follow a similar form, 
but the computational expense associated with achieving higher-order guarantees 
is not a strategy pursued by derivative-free methods that we are aware of.

Models satisfying \eqref{eq:fullylinear} or \eqref{eq:fullyquadratic} are 
called Taylor-like models. To understand why, consider the second-order 
Taylor model
\begin{equation}
 \label{eq:taylormodel}
 m(\xb+\sba) = f(\xb) + \nabla f(\xb)^\T \sba + \ffrac{1}{2}\sba^\T\nabla^2 
 f(\xb)\sba.
\end{equation}
This model is a $\kappab$-fully quadratic model of $f$, with
\[
(\kappaef, 
\kappaeg, \kappaeH) = (\Lip{H}/6, \Lip{H}/2, \Lip{H}),
\]
on 
any $\B{\xb}{\Delta}$, where $f$ has a Lipschitz-continuous second derivative 
with Lipschitz constant $\Lip{H}$.

As illustrated in the next section, one also can guarantee that 
models that do not employ derivative information satisfy these approximation 
bounds in \eqref{eq:fullylinear} or \eqref{eq:fullyquadratic}. This 
approximation quality is used by derivative-free algorithms to ensure that 
a sufficient reduction predicted by the model $m$ yields an attainable 
reduction in the function $f$ as $\Delta$ becomes smaller.

\subsubsection{Polynomial models}
\label{sec:polymodels}

Polynomial models are the most commonly used models for derivative\anlchange{-free} local 
optimization. We let $\cP^{d,n}$ denote the space of polynomials of $n$ 
variables of degree $d$ and $\phib: \Reals^n \to 
\Reals^{\dim{\cP^{d,n}}}$ define a basis for this space. 
For example, quadratic models can be obtained by using the monomial basis
\begin{equation}
 \label{eq:quad_monomials}
 \phib(\xb) = [1, \, x_1, \, \ldots , \, x_n, \, x_1^2, \, \ldots 
x_n^2, \, x_1x_2 , \, \ldots, x_{n-1}x_n]^\T,
\end{equation}
for which $\dim{\cP^{2,n}}= \qn$; linear models can be obtained by using the 
first $\dim{\cP^{1,n}}=n+1$ components of \eqref{eq:quad_monomials}; quadratic 
models with diagonal Hessians, which \anlchange{are} considered by \citeasnoun{Powell2003a}, can 
be obtained by using the first $2n+1$ components of \eqref{eq:quad_monomials}.

Any polynomial model $m\in \cP^{d,n}$ is defined by $\phib$ and coefficients 
$\ab\in \Reals^{\dim{\cP^{d,n}}}$ through 
\begin{equation}
 \label{eq:basic_poly_model}
 m(\xb) = \sum_{i=1}^{\dim{\cP^{d,n}}} a_i \phi_i(\xb).
\end{equation}
Given a set of $p$ points $\Yb = \{\yb_1, \ldots, \yb_p\}$, a model that 
interpolates $f$ on $\Yb$ is defined by the solution $\ab$ to
\begin{equation}
 \label{eq:interpolation_system}
 \Phib(\Yb) \ab
 =
 \begin{bmatrix} 
 \phib(\yb_1) & \cdots & \phib(\yb_p)
\end{bmatrix}^\T 
 \ab
 =
 \begin{bmatrix}
 f(\yb_1) \\ \vdots \\ f(\yb_p)
 \end{bmatrix}\!.
\end{equation}

The existence, uniqueness and conditioning of a solution to 
\eqref{eq:interpolation_system} depend on the location of the sample points 
$\Yb$ through the matrix $\Phib(\Yb)$. We note that when $n>1$, 
$|\Yb|=\dim{\cP^{d,n}}$
is insufficient for guaranteeing that $\Phib(\Yb)$ is non-singular 
\cite{Wendland}. Instead, additional conditions, effectively on the geometry of 
the sample points $\Yb$, must be satisfied.

\paragraph{Simplex gradients and linear interpolation models.} 
\label{sec:det_det_simplexgradients}

The geometry conditions needed to uniquely define a linear model are relatively 
straightforward: the sample points $\Yb$ must be affinely independent; that is, 
the columns of 
\begin{equation}
 \label{eq:unscaled_linear}
 \Yb_{-1} = \begin{bmatrix}
 \yb_2-\yb_1 & \cdots & \yb_{n+1}-\yb_1
 \end{bmatrix}
\end{equation}
must be linearly 
independent. Such sample points define what is referred to as a simplex gradient 
$\gb$ through $\gb = [a_2, \ldots, a_{n+1}]^\T$, when the monomial 
basis $\phi$ is used in \eqref{eq:interpolation_system}.

Simplex gradients can be viewed as a generalization of first-order
finite-difference estimates (\eg\ the forward differences based on
evaluations at the points $\{\yb_1, \yb_1+\Delta \eb_1, \ldots ,
\yb_1+\Delta \eb_n\}$); their use in optimization algorithms dates at
least back to the work of \citeasnoun{Spendley1962} that inspired
\citeasnoun{NelderMead}. Other example usage includes pattern search
\cite{ACLV2007,ACLVJD2008} and noisy optimization
\cite{Kelleybook,Bortz1998}; the study of simplex gradients continues
with recent works such as those of \citeasnoun{Regis2015calculus} and
\citeasnoun{Coope2017}.

Provided that \eqref{eq:unscaled_linear} is non-singular, it is 
straightforward to show that linear interpolation models are 
$\kappab$-fully linear model of $f$ in a neighbourhood of $\yb_1$. In particular, 
if $\Yb\subset \B{\yb_1}{\Delta}$ and $f$ has an $\Lip{g}$-Lipschitz-continuous 
first derivative on an open domain containing $\B{\yb_1}{\Delta}$, then 
\eqref{eq:fullylinear} holds on $\B{\yb_1}{\Delta}$ with 
\begin{equation}
 \label{eq:linearmodelerrors}
\kappaeg= \Lip{g}(1+\sqrt{n}\Delta\|\Yb_{-1}^{-1}\|/2)
\quad \mbox{and} \quad 
\kappaef = \Lip{g}/2+\kappaeg. 
\end{equation}

The expressions in \eqref{eq:linearmodelerrors} also provide a recipe for 
obtaining a model with a potentially tighter error bound over 
$\B{\yb_1}{\Delta}$: modify $\Yb\subset \B{\yb_1}{\Delta}$ to decrease 
$\|\Yb_{-1}^{-1}\|$. 
 We note that when $\Yb_{-1}$ 
contains orthonormal directions scaled by $\Delta$, one recovers
$\kappaeg= \Lip{g}(1+\sqrt{n}/2)$ and $\kappaef = 
\Lip{g}(3+\sqrt{n})/2$, which is the least value one can obtain from
 \eqref{eq:linearmodelerrors} given the restriction that $\Yb\subset 
\B{\yb_1}{\Delta}$.
Hence, by performing LU or QR factorization with pivoting, one can obtain 
directions (which are then scaled by $\Delta$) in order to improve the 
conditioning of $\Yb_{-1}^{-1}$ and hence the approximation bound. Such an 
approach \anlchange{is} performed by \citeasnoun{Conn2006a} for linear models and by 
\citeasnoun{SMWCAS11} for fully linear radial basis function models.

The geometric conditions on $\Yb$, induced by the approximation bounds in 
\eqref{eq:fullylinear} or \eqref{eq:fullyquadratic}, can be viewed as playing a 
similar role to the geometric conditions (\eg\ positive spanning) imposed on 
$\Db$ in directional direct-search methods. Naturally, the choice of basis 
function used for any model affects the quantitative measure of that model's 
quality.

Note that many practical methods employ interpolation sets 
contained within a constant multiple of the trust-region 
radius (\ie\ $\Yb\subset \B{\yb_1}{c_1\Delta}$ for a 
constant $c_1\in[1,\infty)$). 

\paragraph{Quadratic interpolation models.}
\label{sec:det_det_quadmodels}

Quadratic interpolation models have been used for derivative-free optimization 
for at least fifty years \cite{Winf69,Winf73} and were employed by a series of methods 
that revitalized interest in model-based methods;
see, for example, \citeasnoun{Conn1996}, \citeasnoun{Conn1997a}, \citeasnoun{Conn1997b}
and \citename{Powell1998b} \citeyear{Powell1998b,Powell2002a}.

Of course, the quality of an interpolation model (quadratic or otherwise) in a
region of interest is determined by the position of the underlying points being
interpolated. For example, if a model $m$ interpolates a function $f$ at points
far away from a certain region of interest, the model value may differ greatly
from the value of $f$ in that region. $\Lambda$-poisedness is a concept to
measure how well a set of points is dispersed through a region of interest, and
ultimately how well a model will estimate the function in that region.

The most commonly used metric for quantifying how well points are
positioned in a region of interest is based on Lagrange polynomials. Given a
set of $p$ points $\Yb = \{ \yb_1,\ldots,\yb_p \}$, a basis of Lagrange
polynomials satisfies
\begin{equation} \label{eq:Lagrange_poly}
 \ell_j(\yb_i) = 
\begin{cases}
 1 & \mbox{if } i = j,\\
 0 & \mbox{if } i \neq j.
 \end{cases}
\end{equation}

We now define $\Lambda$-poisedness. A set of points $\Yb$ is said to be
$\Lambda$-poised on a set \anlchange{$\Bb$} if $\Yb$ is
linearly independent and the Lagrange
polynomials $\{ \ell_1, \ldots, \ell_p \}$ associated with $\Yb$
satisfy
\begin{equation}\label{eq:Lambda-poised}
 \Lambda \ge \max_{1 \le i \le p} \max_{\xb \in B} | \ell_i(\xb) |.
\end{equation}
(For an equivalent definition of $\Lambda$-poisedness, see
\citeasnoun[Definition~3.6]{Conn2009a}.)
Note that the definition of $\Lambda$-poisedness is independent
of the function being modelled. Also, the points $\Yb$ need not necessarily be
elements of the set \anlchange{$\Bb$}. Also, note that if a model is poised on a set \anlchange{$\Bb$}, it
is poised on any subset of \anlchange{$\Bb$}. One is usually interested in the least value of
$\Lambda$ so that \eqref{eq:Lambda-poised} holds. 

Powell's unconstrained optimization by quadratic approximation method (\codes{UOBYQA}) follows such an approach in maximizing the Lagrange 
polynomials. In \citeasnoun{Powell1998b}, \citeasnoun{Powell2001a} and \citeasnoun{Powell2002a}, significant care 
is given to the linear algebra expense associated with this maximization and 
the associated change of basis as the methods change their interpolation sets. 
For example, in \citeasnoun{Powell1998b}, particular sparsity in the Hessian 
approximation is employed with the aim of capturing curvature while keeping 
linear algebraic expenses low. 

Maintaining, and the question of to what extent it is necessary to maintain, 
this 
geometry for quadratic models has been intensely studied; see, for example, 
\citeasnoun{GFJMJN09}, \citeasnoun{Marazzi2002a}, \citeasnoun{DAmbrosio2017} and \citeasnoun{Scheinberg10}.

\paragraph{Underdetermined quadratic interpolation models.}
\label{sec:det_det_powellmodels}

A fact not to be overlooked in the context of derivative-free optimization is 
that employing an interpolation set $\Yb$ requires availability of the 
$|\Yb|$ function values $\{f(\yb_i): \yb_i\in\Yb\}$. When the function $f$ is 
computationally expensive to evaluate, the $\qn$ points required by fully 
quadratic models can be a burden, potentially with little benefit, to obtain 
repeatedly in an optimization algorithm.

Beginning with \citeasnoun{Powell2003a}, Powell investigated quadratic models 
constructed from fewer than $\qn$ points. The most successful of these 
strategies was detailed in \citeasnoun{Powell2004a} and \citeasnoun{Powell2004b} and
resolved the $\qn-|\Yb|$ remaining degrees of 
freedom by solving problems of the form
\begin{equation}\label{eq:minchange}
\begin{aligned}
 &\minimize_{m\in \cP^{2,n}} && \|\nabla^2 m(\check{\xb}) - \Hb\|^2_F \\
 &\mbox{\,subject to} && m(\yb_i) = f(\yb_i), \quad \gforall\ \yb_i \in \Yb
\end{aligned}
\end{equation}
to obtain a model $m$ about a point of interest $\check{\xb}$. Solutions to 
\eqref{eq:minchange} are models with a Hessian closest in Frobenius norm to a specified 
$\Hb=\Hb^\T$ among all models that interpolate $f$ on $\Yb$. A popular 
implementation of this strategy is the \codes{NEWUOA} solver \cite{MJDP06}.

By using the basis 
\begin{align}
 \label{eq:quad_powell}
 \phib(\check{\xb}+\xb) &= \bigl[\phib_{\rm fg}(\check{\xb}+\xb)^\T && | \, 
\phib_{\rm H}(\check{\xb}+\xb)^\T\bigr]^\T \\*
&= \biggl[1, \, x_1, \, \ldots , \, x_n && \bigg| \,
\ffrac{1}{2} x_1^2, \, \ldots, \, \ffrac{1}{2} x_n^2, \, \ffrac{1}{\sqrt{2}} 
x_1x_2 , 
\, \ldots, \ffrac{1}{\sqrt{2}} x_{n-1}x_n\biggr]^\T, \notag
\end{align}
the problem \eqref{eq:minchange} is equivalent to the problem
\begin{align}\label{eq:seminorm}
 &\minimize_{\ab_{\rm fg}, \ab_{\rm H}} && \|\ab_{\rm 
H}\|^2_2 \\*
 &\mbox{\,subject to} && \ab_{\rm fg}^\T\phib_{\rm fg}(\yb_i) 
 + \ab_{\rm H}^\T \phib_{\rm H}(\yb_i) 
 = f(\yb_i) - \ffrac{1}{2} \yb_i^\T \Hb \yb_i, \quad \gforall\ \yb_i \in \Yb.
 \notag
\end{align}
Existence and uniqueness of solutions to \eqref{eq:seminorm} again depend on 
the positioning of the points in $\Yb$. Notably, a necessary condition for 
there to be a unique minimizer of the seminorm is that at least 
$n+1$ of the points in $\Yb$ be affinely independent. Lagrange polynomials can 
be defined for this case; \citeasnoun{Conn2006b} \anlchange{establish} conditions for 
$\Lambda$-poisedness (and hence a fully linear, or better, approximation 
quality) of such models. 

\citename{Powell2004d} \citeyear{Powell2004d,Powell2007a,MJDP07b} \anlchange{develops} efficient solution methodologies for \eqref{eq:seminorm} when 
$\Hb$ and $m$ are constructed from interpolation sets that \anlchange{differ} by at most 
one point, and \anlchange{employ} these updates in 
\codes{NEWUOA} and subsequent solvers. \anlchange{\citeasnoun{SW08} and \citeasnoun{CRV2008} 
use} $\Hb=\zerob$ in order to obtain tighter fully linear error bounds of 
models resulting from \eqref{eq:seminorm}. A strategy of using even fewer 
interpolation points (including those in a proper subspace of $\Reals^n$) \anlchange{is} 
developed by \citeasnoun{Powell2013beyond} and \citeasnoun{Zhang2014}.
In \secref{sepspa}, we summarize approaches that exploit knowledge of sparsity of the 
derivatives of $f$ in building quadratic models that interpolate fewer than $\qn$ points.

\begin{figure}
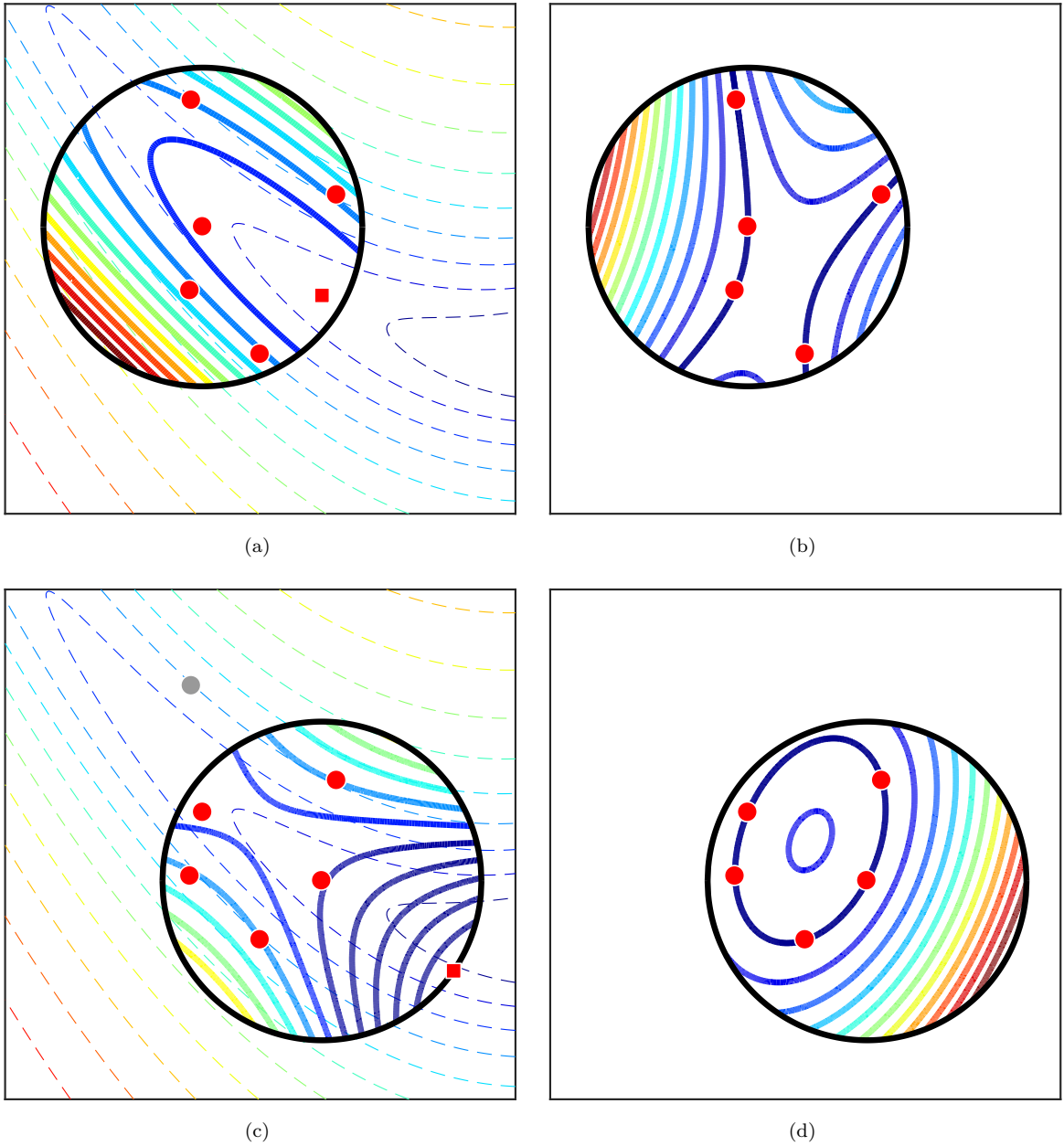
 %
\centering
\subfigure[\label{fig:iter_1}]{\includegraphics[width=.49\linewidth]{{{fig2_2a}}}}\hfill %
\subfigure[\label{fig:iter_1_pts}]{\includegraphics[width=.49\linewidth]{{{fig2_2b}}}}\\ %
\subfigure[\label{fig:iter_2}]{\includegraphics[width=.49\linewidth]{{{fig2_2c}}}}\hfill %
\subfigure[\label{fig:iter_2_pts}]{\includegraphics[width=.49\linewidth]{{{fig2_2d}}}}\\ %
\caption{
 (a) Minimum-norm-Hessian model through five points in $\B{\xb_k}{\Delta_k}$ 
and 
its minimizer.
 (b) Absolute value of a sixth Lagrange polynomial for the five points.
 (c) Minimum-norm-Hessian model through five points in 
$\B{\xb_{k+1}}{\Delta_{k+1}}$ and its minimizer.
 (d) Absolute value of a sixth Lagrange polynomial for the five points.
\label{fig:dfotr}}
\end{figure}

\fo{Figure~\ref{fig:dfotr}} shows quadratic models in two dimensions that interpolate
$\qn-1=5$ points as well as the associated magnitude of the remaining Lagrange 
polynomial (note that this polynomial vanishes at the five interpolated points).

\paragraph{Regression models.}
\label{sec:det_det_regmodels}

Just as one can establish approximation bounds and geometry conditions 
when $\Yb$ is linearly independent, the same can be done for 
overdetermined regression models \cite{Conn2006b,Conn2009a}. 
This can be accomplished by extending the definition of Lagrange polynomials
from \eqref{eq:Lagrange_poly} to the regression case. That is, given a
basis $\phib: \Reals^n \to \Reals^{\dim{\cP^{d,n}}}$ and points $\Yb =
\{ \yb_1,\ldots,\yb_p \}\subset \Reals^n$ with $p>\dim{\cP^{d,n}}$,
the set of polynomials satisfies 
\begin{equation} \label{eq:Lagrange_poly_reg}
 \ell_j(\yb_i) \overset{\rm l.s.}{=}
 \begin{cases}
 1 & \mbox{if } i = j,\\
 0 & \mbox{if } i \neq j,
 \end{cases}
\end{equation}
where $\overset{\rm l.s.}{=}$ denotes the least-squares solution.
The regression model can be recovered finding the least-squares solution (now
overdetermined) system from \eqref{eq:interpolation_system}, and the definition
of $\Lambda$-poisedness (in the regression sense) is equivalent to
\eqref{eq:Lambda-poised}.
Ultimately, given a linear regression model through a set of $\Lambda$-poised points
$\Yb\subset \B{\yb_1}{\Delta}$, and if $f$ has an $\Lip{g}$-Lipschitz-continuous 
first derivative on an open domain containing $\B{\yb_1}{\Delta}$, then 
\eqref{eq:fullylinear} holds on $\B{\yb_1}{\Delta}$ with 
\begin{equation}
 \label{eq:linearmodelerrors_reg}
 \kappaeg= \ffrac{5}{2}\sqrt{p} \Lip{g}\Lambda
\quad \mbox{and} \quad 
\kappaef = \ffrac{1}{2}\Lip{g}+\kappaeg. 
\end{equation}
\citeasnoun{Conn2006b} note the fact that the extension of Lagrange
polynomials does not apply to the 1-norm or infinity-norm case.
\citeasnoun{BLG2013} show that the definition of Lagrange polynomials can be extended
to the weighted regression case.
\citeasnoun{Verdrio2017} show that \eqref{eq:fullylinear} can also be recovered 
for support vector regression models.

Efficiently minimizing the model (regardless of type) over a trust region is
integral to the usefulness of such models within an optimization algorithm. 
In fact, this necessity is a primary reason for the use of 
low-degree polynomial models by the majority of derivative-free trust-region 
methods. For quadratic models, the resulting subproblem remains one of 
the most difficult non-convex optimization problems solvable in 
polynomial time, as illustrated by \citeasnoun{More553}. As exemplified by
\citeasnoun{Powell1997a}, 
the implementation of subproblem solvers is a key concern in methods seeking to 
perform as few algebraic operations between function evaluations as possible.

\subsubsection{Radial basis function interpolation models}
\label{sec:rbfmodels}

An additional way to model non-linearity with potentially less 
restrictive geometric conditions is by using radial basis 
functions (RBFs). Such models take the form
\begin{equation}
 \label{eq:rbf_model}
 m(\xb) = \sum_{i=1}^{|\Yb|} b_i \psi(\|\xb-\yb_i\|) + 
\ab^\T\phib(\xb),
\end{equation}
where
$\psi:\Reals_+ \to \Reals$ is a conditionally positive-definite 
univariate function and $\ab^\T\phib(\xb)$ represents a 
(typically low-order) polynomial as before; see, for example, \citeasnoun{Buhmann2000}. 
Given a sample set $\Yb$, RBF model coefficients $(\ab,\bb)$ can be 
obtained by solving the augmented interpolation equations
\begin{equation}
\label{eq:rbfsystem}
\begin{bmatrix} 
 \psi(\|\yb_1-\yb_1\|) & \cdots & 
\psi(\|\yb_1-\yb_{|\Yb|}\|) & \phib(\yb_1)^\T \\
 \vdots & & \vdots & \vdots \\
 \psi(\|\yb_{|\Yb|}-\yb_1\|) & \cdots & 
\psi(\|\yb_{|\Yb|}-\yb_{|\Yb|}\|) & \phib(\yb_{|\Yb|})^\T \\
 \phib(\yb_1) & \cdots & \phib(\yb_{|\Yb|}) & \zerob
 \end{bmatrix}
\begin{bmatrix}
 \bb \\ 
\ab
 \end{bmatrix}
=
\begin{bmatrix}
 f(\yb_1)\\
 \vdots \\ 
 f(\yb_{|\Yb|})\\ 
 \zerob
 \end{bmatrix}\!.
\end{equation}

That RBFs are conditionally positive-definite ensures that \eqref{eq:rbfsystem} is non-singular 
provided that the degree $d$ of the polynomial $\phib$ is sufficiently large 
and that $\Yb$ is poised for degree-$d$ polynomial interpolation. For example, 
cubic ($\psi(r)=r^3$) RBFs require a linear polynomial; 
multiquadric ($\psi(r)=-(\gamma^2+r^2)^{1/2}$) RBFs require a constant 
polynomial; and inverse multiquadric ($\psi(r)=(\gamma^2+r^2)^{-1/2}$) 
and Gaussian ($\psi(r)=\exp(-\gamma^{-2} r^2)$) RBFs do not 
require a polynomial. Consequently, RBFs have relatively unrestrictive 
geometric requirements on the interpolation points $\Yb$ while allowing for 
modelling a wide range of non-linear behaviour.

This feature is typically exploited in global optimization
(see \eg\ \citeb{Bjorkman2000a}, \citeb{Gutmann2001b} and \citeb{RGRCAS2007}), whereby an RBF 
surrogate model is employed to globally approximate $f$. However, 
works such as \citeasnoun{Oeuvray2008}, \citeasnoun{Oeuvray2005a}, \citeasnoun{WildPhD} and \citeasnoun{SMWCAS13} \anlchange{establish and 
use} local approximation properties of these models. This approach is typically 
performed by relying on a linear polynomial $\ab^\T\phib(\xb)$, which can 
be used to establish that the RBF model in \eqref{eq:rbf_model} can be a fully 
linear local approximation of smooth $f$.

\subsubsection{Trust-region methods}
\label{sec:det_det_TRM}

Having discussed issues of model construction, we are now ready to present 
a general statement of a model-based trust-region method in \algref{tr}.

\begin{algorithm}
 Set parameters $\epsilon > 0$, $0 < \gammad < 1 \le \gammai$, $0 < \eta_0 \le \eta_1 < 1$, $\Deltamax$\\
 Choose initial point $\xb_0$, trust-region radius $0 < \Delta_0 \le 
\Deltamax$, and set of previously evaluated points $\Yb_k$\\
 \For{$k = 0,1,2\ldots$}
 {
 Select a subset of $\Yb_k$ (or augment $\Yb_k$ and evaluate $f$ at new points) for model building\\
 Build a model $m_k$ using points in $\Yb_k$ and their function values\\
 \While{$\| \nabla m_k(\xb_k) \| < \epsilon $}
 {
 \If{$m_k$ is accurate on $\B{\xb_k}{\Delta_k}$}
 {
 $\Delta_k \gets \gammad \Delta_k$
 }
 \Else{By updating $\Yb_k$, make $m_k$ accurate on $\B{\xb_k}{\Delta_k}$}
 }
 Generate a direction $\sba_k\in \B{\zerob}{\Delta_k}$ so that $\xb_k + \sba_k$
 approximately minimizes $m_k$ on $\B{\xb_k}{\Delta_k}$ \label{line:trsp}\\[3pt]

 Evaluate $f(\xb_k+\sba_k)$ and $\rho_k \gets \ffrac{f(\xb_k) - f(\xb_k + 
\sba_k)}{m_k(\xb_k)
 - m_k(\xb_k+\sba_k)}$\\[3pt]

 \If{$\rho_k < \eta_1$ and $m_k$ is inaccurate on $\B{\xb_k}{\Delta_k}$}
 {
 Add model improving point(s) to $\Yb_k$ 
 }

 \uIf{$\rho_k \ge \eta_1$}
 {
 $\Delta_{k+1} \gets \min\{\gammai \Delta_k, \Deltamax\}$
 }
 \uElseIf{$m_k$ is accurate on $ \B{\xb_k}{\Delta_k}$}
 {
 $\Delta_{k+1} \gets \gammad \Delta_k$
 }
 \Else{$\Delta_{k+1} \gets \Delta_k$}

 \leIf{$\rho_k \ge \eta_0$}
 {
 $\xb_{k+1} \gets \xb_k + \sba_k$
 }
 {
 $\xb_{k+1} \gets \xb_k$
 }

 $\Yb_{k+1} \gets \Yb_k$
 }

 \caption{Derivative-free model-based trust-region method\label{alg:tr}}
\end{algorithm}

A distinguishing characteristic of derivative-free model-based trust-region methods
is how they manage $\Yb_k$, the set of points used to construct the model
$m_k$. Some methods ensure that $\Yb_k$ contains a scaled
stencil of points around $\xb_k$; such an approach can be attractive since the
objective at such points can be evaluated in parallel. A fixed stencil can also
ensure that all models sufficiently approximate the objective. Other methods
construct $\Yb$ by using previously evaluated points near $\xb_k$, for example,
those points within $\B{\xb_k}{c_1\Delta_k}$ for some constant $c_1 \in
[1,\infty)$. Depending on the set of previously evaluated points, 
such methods may need to add points to $\Yb_k$ that most improve the model quality.
Determining which additional points to add to $\Yb_k$
can be computationally expensive, but the method should be willing to do so in the hope of needing fewer 
evaluations of the objective function at new points in $\Yb_k$. Most methods 
do not ensure that models are valid on every iteration but rather
make a single step toward improving the model. Such an approach can ensure a
high-quality model in a finite number of improvement steps. (Exceptional methods
that ensure model quality before $\sba_k$ is calculated are the methods of
Powell and manifold sampling of \citeasnoun{KLW18}.)
The \codes{ORBIT} method \cite{SWRRCS07} places a limit on the size of $\Yb_k$ 
(\eg\ in order to
limit the amount of linear algebra or to prevent overfitting).
In the end, such restrictions on $\Yb_k$ may determine whether $m_k$ is an 
interpolation or regression model.

Derivative-free trust-region methods share many similarities with traditional
trust-region methods, for example, 
the use of a $\rho$-test to determine whether a step
is taken or rejected. 
As in a traditional trust-region method, the $\rho$-test
measures the ratio of actual decrease observed in the
objective versus the decrease predicted by the model.

On the other hand, the management of the trust-region radius parameter $\Delta_k$
in \algref{tr} differs 
remarkably from traditional trust-region
methods.
Derivative-free variants require an additional test of model quality, the failure
of which results in shrinking $\Delta_k$. 
When derivatives are available, Taylor's theorem ensures
model accuracy for small $\Delta_k$. In the derivative-free case, such a condition
must be explicitly checked in order to ensure that $\Delta_k$ does not go to zero
merely because the model is poor, hence the inclusion of tests of model quality.
As a direct result of these considerations, $\Delta_k \to 0$ as
\algref{tr} converges; this is generally not the case in traditional trust-region methods.

As in derivative-based trust-region methods, 
the solution to the trust-region subproblem in \lineref{trsp} of \algref{tr}
must satisfy a \emph{Cauchy decrease condition}.
Given the model $m_k$ used in \algref{tr}, 
we define the optimal step length in the direction $-\nabla m_k(\xb_k)$ by
\[
t_k^C = \argmin_{t\geq 0:\xb_k-t\nabla m_k(\xb_k)\in\B{\xb_k}{\Delta_k}} m_k(\xb_k-t\nabla m_k(\xb_k)),
\]
and the corresponding \emph{Cauchy step} 
\[
\sba_k^C = -t_k^C\nabla m_k(\xb_k).
\]
It is straightforward to show (see \eg\ \citeb{Conn2009a}, Theorem~10.1) that
\begin{equation}
\label{eq:cauchy_decrease}
m_k(\xb_k)-m_k(\xb_k+\sba_k^C) \geq \ffrac{1}{2}\|\nabla m_k(\xb_k)\|\min\biggl\{\ffrac{\|\nabla m_k(\xb_k)\|}{\|\nabla^2 m_k(\xb_k)\|},\Delta_k\biggr\}.
\end{equation}
That is, \eqref{eq:cauchy_decrease} states that,
provided that both $\Delta_k \approx \|\nabla m_k(\xb_k)\|$ and
a uniform bound exists on the norm of the model Hessian, 
the model decrease attained by the Cauchy step $\sba_k^C$ is of the order of
$\Delta_k^2$. 
In order to prove convergence, it is desirable
to ensure that each step $\sba_k$ generated in \lineref{trsp} of \algref{tr}
decreases the model $m_k$ by no less than $\sba_k^{C}$ does, or at least some
fixed positive fraction of the decrease achieved by $\sba_k^{C}$.
Because successful iterations ensure that the actual decrease attained in an iteration is at least
a constant fraction of the model decrease, the sequence of decreases of \algref{tr} are square-summable,
provided that $\Delta_k\to 0$. (This is indeed the case for derivative-free trust-region methods.) 
Hence, in most theoretical treatments of these methods, it is commonly
stated as an assumption that the subproblem solution $\sba_k$ obtained in 
\lineref{trsp}
of \algref{tr} satisfies 
\begin{equation}
\label{eq:cauchy_assumption}
m_k(\xb_k) - m_k(\xb_k+\sba_k) \geq \kappa_{\rm 
fcd}(m_k(\xb_k)-m_k(\xb_k+\sba_k^C)),
\end{equation}
where $\kappa_{\rm fcd}\in(0,1]$ is the
fraction of the Cauchy decrease. 
In practice, when $m_k$ is a quadratic model, subproblem solvers have been 
well studied and often come with guarantees concerning the satisfaction of 
\eqref{eq:cauchy_assumption} \cite{Trmbook}.
\citeasnoun[Figure~4.3]{SWRRCS07} demonstrate the satisfaction of an assumption like
\eqref{eq:cauchy_assumption} when the model $m_k$ is a radial basis function. 

Under reasonable smoothness assumptions, 
most importantly $f\in\cLC^1$,
algorithms in the \algref{tr} framework have been shown to be first-order
convergent (\ie\ \eqref{eq:first_order_convergence}) 
and second-order convergent (\ie\ \eqref{eq:second_order_convergence}),
 with the (arguably) most well-known proof given
by \citeasnoun{Conn2009}.
In more recent work, \citeasnoun{Garmanjani2016} provide a WCC bound
of the form \eqref{eq:nonconvex_stationarity}
for \algref{tr}, recovering essentially the same upper bound 
on the number of function evaluations required by \framework{DDS} methods found in \citeasnoun{Vicente2013},
that is, a WCC bound in $\bigo{\epsilon^{-2}}$ (see \tabref{rates}).
When $f\in\cLC^2$, \citeasnoun{Gratton2017b} demonstrate a second-order WCC bound of the form
\eqref{eq:2ndnonconvex_stationarity} in $\bigo{\epsilon^{-3}}$; in
order to achieve this result, fully quadratic models $m_k$ are required.
In \secref{rand_det_tr}, a similar result is achieved by using randomized variants that 
do not require a fully quadratic model in every iteration.

Early analysis of Powell's \codes{UOBYQA} method shows that, with minor
modifications, %
 the 
algorithm can converge superlinearly in neighbourhoods of strict 
convexity \cite{Han2004}. A key distinction between Powell's methods and other model-based trust-region methods is the use of separate neighbourhoods for model quality and trust-region steps, with each of these neighbourhoods changing dynamically. Convergence of such methods \anlchange{is} addressed by 
\citename{Powell2010convergence}\ \citeyear{Powell2010convergence,Powell2012convergence}. 

The literature on derivative-free trust-region methods is extensive. 
We mention in passing several additional classes of trust-region methods
that have not fallen neatly into our discussion thus far.
\emph{Wedge methods} \cite{Marazzi2002a} explicitly enforce geometric properties
($\Lambda$-poisedness) of the sample set between iterations by adding additional constraints to
the trust-region subproblem.
\citeasnoun{Alexandrov1998} consider a trust-region method utilizing a hierarchy of
model approximations.
In particular, if derivatives can be obtained
but are expensive, then the method of \citeasnoun{Alexandrov1998} uses
a model that interpolates not only zeroth-order 
information
but also first-order (gradient) information. 
For problems with deterministic noise, 
\citeasnoun{Elster95} propose a method that projects the solutions of a trust-region subproblem
onto a dynamically refined grid, encouraging better practical behaviour. 
Similarly, for problems with deterministic noise, 
\citeasnoun{Maggiar2018}
propose a model-based trust-region method that implicitly convolves the objective function 
with a Gaussian kernel, again yielding better practical behaviour.

\subsection{Hybrid methods and miscellanea}
\label{sec:det_det_others}

While the majority of work in derivative-free methods for deterministic
problems can be classified as direct-search or model-based methods, some work
defies this simple classification. In fact, several works
\cite{Conn2013,CRV2008,Dennis1997managing,Frimannslund2011} propose methods that seem to
hybridize these two classes, existing somewhere in the intersection. For
example, \citeasnoun{Custodio2005} and \citeasnoun{CRV2008} develop the \codes{SID-PSM}
method, which extends \algref{ds} so that the search step consists of
minimizing an approximate quadratic model of the objective (obtained either by
minimum-Frobenius norm interpolation or by regression) over a trust region. Here, we highlight
methods that do not neatly belong to the two aforementioned classes of methods.

\subsubsection{Finite differences}
\label{sec:orphaned-fd}

As noted in \secref{intro:fd}, many of the earliest derivative-free methods 
employed finite-difference-based estimates of derivatives. The most popular 
first-order directional derivative estimates include the 
forward/reverse difference
\begin{equation}
 \df{f;\xb;\db;h} \defined \ffrac{f(\xb+h\db)-f(\xb)}{h}
 \label{eq:fdiff}
\end{equation}
and central difference
\begin{equation}
 \dc{f;\xb;\db;h} \defined \ffrac{f(\xb+h\db)-f(\xb-h\db)}{2h},
 \label{eq:cdiff}
\end{equation}
where $h\neq 0$ is the difference parameter and the non-trivial $\db\in 
\Reals^n$ defines the direction. Several recent methods, including 
the methods described in
Sections~\ref{sec:imfil}, \ref{sec:arc} and \ref{sec:rand_det_rs}, use such estimates and employ 
difference parameters or directions that dynamically change.

As an example of a potentially dynamic choice of difference parameter, 
we consider the usual case of roundoff errors. We denote 
by $f_{\infty}'(\xb;\db)$ the directional derivative at $\xb$ of the 
infinite-precision 
(\ie\ based on real arithmetic) objective function
$f_{\infty}$ in the unit direction $\db$ (\ie\ $\|\db\|=1$). We then have the 
following error for forward or reverse finite-difference estimates based on the 
function $f$ available through computation:
\begin{equation}
\label{eq:roundoffbound}
 |\df{f;\xb;\db;h} - f_{\infty}'(\xb;\db) |
 \leq \ffrac{1}{2} \Lip{g}(\xb) |h| + 2 
\ffrac{\epsilon_{\infty}(\xb)}{|h|},
\end{equation}
provided that $|f''_{\infty}(\cdot;\db)|\leq \Lip{g}(\xb)$ and 
$|f_{\infty}(\cdot)-f(\cdot)|\leq \epsilon_{\infty}(\xb)$ on the interval 
$[\xb, 
\xb+h\db]$. In \citeasnoun{gill1981po} and \citeasnoun{gill1983cfd}, the recommended 
difference 
parameter is $h= 2\sqrt{\epsilon_{\infty}(\xb)/\Lip{g}(\xb)}$, 
which yields the minimum value 
$2\sqrt{\epsilon_{\infty}(\xb) \Lip{g}(\xb)}$ 
of the upper bound in \eqref{eq:roundoffbound}; 
when $\epsilon_{\infty}$ is a bound on the roundoff 
error and $\Lip{g}$ is of order one, then the familiar 
$h\in\bigo{\sqrt{\epsilon_{\infty}}}$ is obtained.

Similarly, if one 
models the error between $f_{\infty}$ and $f$ as a stationary stochastic 
process 
(through the ansatz denoted by $f_{\xib}$)
with variance $\epsilon_{\rm f}(\xb)^2$, minimizing the upper bound on the 
mean-squared error,
\begin{equation}
\label{eq:msebound}
 \E[\xib]{(\df{f_{\xib};\xb;\db;h} - f_{\infty}'(\xb;\db) 
)^2}
\leq \ffrac{1}{4} \Lip{g}(\xb)^2 h^2 + 2 
\ffrac{\epsilon_{\rm f}(\xb)^2}{h^2},
\end{equation}
yields the choice 
$h = ({\sqrt{8} \epsilon_{\rm f}(\xb)/\Lip{g}(\xb)})^{1/2}$
with an associated root-mean-squared error of 
$({\sqrt{2}\epsilon_{\infty}(\xb) \Lip{g}(\xb)})^{1/2}$;
see, for example,
\citename{more2011edn} \citeyear{more2011edn,more2014nd}. 
A rough procedure for computing 
$\epsilon_{\rm f}$ is provided in \citeasnoun{more2011ecn} and used in recent 
methods such as that of \citeasnoun{BBN2018}.

In both cases \eqref{eq:roundoffbound} and \eqref{eq:msebound}, the first-order 
error is $c \sqrt{\epsilon (\xb) \Lip{g}(\xb)}$ (for a constant 
$c\leq 2$), which can be used to guide the decision on whether the derivatives 
estimates are of sufficient accuracy.

\subsubsection{Implicit filtering}
\label{sec:imfil}

Implicit filtering is a hybrid of a grid-search algorithm (evaluating all
points on a lattice) and a Newton-like local optimization method. The gradient 
(and possible Hessian) estimates for local optimization are approximated by 
the central 
differences $\{\dc{f;\xb_k;\eb_i;\Delta_k}: i=1,\ldots,n\}$. The difference 
parameter $\Delta_k$ decreases 
when implicit
filtering encounters a \emph{stencil failure} \anlchange{at $\xb_k$}, that is,
\begin{equation}\label{eq:stencil_failure}
 f(\xb_k) \le f(\xb_k \pm \Delta_k \eb_i), 
\end{equation}
where $\eb_i$ is the $i$th \anlchange{elementary basis vector}. This is similar to direct-search
methods, but notice that implicit filtering is not polling opportunistically:
all polling points are evaluated on each iteration.
\begin{algorithm}
 Set parameters $\mathtt{feval\_max} > 0$, $\Deltamin > 0$, $\gammad \in (0,1)$ and $\tau > 0$\\
 Choose initial point $\xb_0$ and step size $\Delta_0 \geq \Deltamin$\\
 $k \gets 0$; evaluate $f(\xb_0)$ and set $\mathtt{fevals} \gets 1$ \\
 \While{$\mathtt{fevals} \le \mathtt{feval\_max}$ and $\Delta_k \ge \Deltamin$}
 {
 Evaluate $f(\xb_k \pm \Delta_k \eb_i)$ for $i \in \{1,\ldots,n\}$ and 
approximate $\nabla f(\xb_k)$ via $\{\dc{f;\xb_k;\eb_i;\Delta_k}: 
i=1,\ldots,n\}$\label{line:stencil}\\
 \eIf{equation \eqref{eq:stencil_failure} is satisfied or $\| \nabla f(\xb_k) \| \le \tau \Delta_k$}
 {
 $\Delta_{k+1} \gets \gammad \Delta_k$\\
 $\xb_{k+1} \gets \xb_k$
 }
 {
 Update Hessian estimate $\Hb_k$ (or set $\Hb_k \gets \Ib$)\\
 $\sba_k \gets -\Hb_k^{-1}\nabla f(\xb_k)$\\
 Perform a line search in the direction $\sba_k$ to generate $\xb_{k+1}$\\
 $\Delta_{k+1} \gets \Delta_k$
 }
 $k \gets k+1$
 }
 \caption{Implicit-filtering method \label{alg:imfil}}
\end{algorithm}
The basic version of implicit filtering from \citeasnoun{Kel2011} is
outlined in \algref{imfil}. Note that most implementations of
implicit filtering require a bound-constrained domain.

Considerable effort has been devoted to extensions of \algref{imfil} when
$f$ is `noisy'.
\citeasnoun{PGCTK95} show that implicit filtering converges to local minima of
\eqref{eq:det_prob} when the objective $f$ is the sum of a smooth function
$f_s$ and a high-frequency, low-amplitude function $f_n$, with $f_n \to 0$
quickly in a neighbourhood of all minimizers of $f_s$.
Under similar assumptions,
\citeasnoun{Kelley2000a} show that \algref{imfil} converges superlinearly if the
step sizes $\Delta_k$ are defined as a power of the norm of the previous
iteration's gradient approximation.

\subsubsection{Adaptive regularized methods}
\label{sec:arc}

\citeasnoun{Cartis2012} perform an analysis of adaptive regularized 
cubic (\codes{ARC}) methods and propose
a derivative-free method, \codes{ARC-DFO}. 
\codes{ARC-DFO} is an extension of \codes{ARC} whereby gradients are replaced with central 
finite differences of the form \eqref{eq:cdiff}, with the difference parameter 
monotonically decreasing within a single iteration of the method.
\codes{ARC-DFO} is an intrinsically model-based method akin to 
\algref{tr}, but the objective within each subproblem regularizes third-order behaviour of the model.
Thus, like a trust-region method, \codes{ARC-DFO} employs trial steps and model gradients.
During the main loop of \codes{ARC-DFO},
 if the difference parameter exceeds a constant factor
of the minimum of the trial step norm or the model gradient norm, then the
difference parameter is shrunk by a constant factor, and the iteration restarts 
to obtain a new trial step. 
This mechanism is structurally similar to a derivative-free trust-region method's
checks on model quality. 
\citeasnoun{Cartis2012} show that \codes{ARC-DFO} demonstrates a WCC result of type
\eqref{eq:nonconvex_stationarity}
in $\bigo{\epsilon^{-3/2}}$, the same asymptotic result (in terms of
$\epsilon$-dependence)
that the authors demonstrate for derivative-based variants of \codes{ARC} methods. 
In terms of dependence on $\epsilon$, this result is a strict improvement over the WCC results of
the same type demonstrated for \framework{DDS} and trust-region methods, 
although this result is proved under the stronger assumption that $f\in\cLC^2$. 

In a different approach, \citeasnoun{Hare2013c} show convergence of a
derivative-free method that penalizes large steps via a proximal regularizer,
thereby removing the necessity for a trust region. \citeasnoun{Lazar2016}
regularize their line-search with a term seeking to minimize a weighted change
of the model's third derivatives.

\subsubsection{Line-search-based methods}
\label{sec:linesearch}
Several line-search-based methods for derivative-free optimization have been 
developed.
\citeasnoun{Grippo1988} and \citeasnoun{DeLeone1984} (two of the few papers appearing in the
1980s concerning derivative-free optimization) both analyse conditions on the
step sizes used in a derivative-free line-search algorithm, and provide
methods for constructing such steps.
\citeasnoun{Luicidi2002} present methods that combine pattern-search and line-search
approaches in a convergent framework.
The \codes{VXQR} method of \citeasnoun{Neumaier2011} performs a line search on a
direction computed from a QR factorization of previously evaluated points.
\citeasnoun{Neumaier2011} apply \codes{VXQR} to problems with $n = 1000$, 
a large problem dimension among the methods considered here. 

Consideration has also been given to non-monotone 
line-search-based derivative-free methods. 
Since gradients are not available 
in derivative-free optimization, the search direction in a line-search method may not be a descent direction. Non-monotone methods allow one to still employ such directions in a globally convergent framework.
\citeasnoun{Grippo2007} extend line-search strategies based on coordinate search and
the method of \citeasnoun{BARZILAI1988} to develop a globally convergent non-monotone
derivative-free method. \citeasnoun{Grippo2014} extend such non-monotone strategies to broader
classes of algorithms that employ simplex gradients, hence further unifying
direct-search and model-based methods.
Another non-monotone line-search method is proposed by \citeasnoun{DinizEhrhardt2008}, who 
encapsulate early examples of randomized \framework{DDS} methods 
(\secref{rand_det_ds}).

\subsubsection{Methods for non-smooth optimization}
\label{sec:nsmooth}

In \secref{det_det_DDS}, we discuss how \codes{MADS} handles non-differentiable
objective functions by densely sampling directions on a mesh, thereby ensuring
that all Clarke directional derivatives are non-negative (\ie\ 
\eqref{eq:uncon_Clarke_stationary}). 
Another early analysis of a DDS method on a class of non-smooth objectives 
was performed by \citeasnoun{Garcia-Palomares2002}.

Gradient sampling methods are a developing class of algorithms for general
non-smooth non-convex optimization; see the recent survey by \citeasnoun{BCLOS18}. 
These methods attempt to estimate the \emph{$\epsilon$-subdifferential} at
a point $\xb$ by evaluating a random sample of gradients in the neighbourhood of
$\xb$ and constructing the convex hull of these gradients. In a
derivative-free setting, the approximation of these gradients is not as
immediately obvious in the presence of non-smoothness, but 
there exist
gradient-sampling methods that use 
finite-difference estimates with specific smoothing techniques \cite{Kiwiel10}. 

In another distinct line of research, \citeasnoun{Bagirov2007} analyse 
a derivative-free variant of subgradient descent, where subgradients are
approximated via so-called discrete gradients. 
In \secref{nso}, we will further discuss
methods for minimizing \emph{composite non-smooth objective functions}
of the form $f = h\circ F$, where $h$ is non-smooth but a closed-form expression is known
and $F$ is assumed smooth. These methods are characterized
by their exploitation of the knowledge of $h$, making them less general than the methods
for non-smooth optimization discussed so far.

\vspace{5pt}
\section{Randomized methods for deterministic objectives} %
\label{sec:rand_det}
We now summarize randomized methods for solving \eqref{eq:det_prob}.
Such methods often have promising theoretical properties, although some
practitioners may dislike the non-deterministic behaviour of these methods. 
We discuss
randomization within direct-search methods in \secref{rand_det_ds} and within
trust-region methods in \secref{rand_det_tr}, but we first begin with a discussion of 
random search as applied to deterministic objectives.

In any theoretical treatment of randomized methods, one must be careful 
to distinguish between random variables and their realizations.
For the sake of terseness in this survey, we will intentionally conflate 
variables
with realizations and refer to respective papers for more careful statements of
theoretical results.

\subsection{Random search}
\label{sec:RS}

We highlight two randomized methods for minimizing a deterministic
objective: pure random search and Nesterov random search. 

\subsubsection{Pure random search}
\label{sec:pureRS}

Pure random search is a natural method to start with for randomized derivative-free optimization.
Pure random search is popular for multiple reasons; in particular, it 
is easy to implement (with few or no user-defined tolerances), and (if the 
points generated are independent of one another) 
it exhibits perfect scaling in terms of evaluating $f$ at many points 
simultaneously.

A pure random-search method is given in \algref{pure_random}, where
points are generated randomly from $\Omegab$. For example, if $\Omegab = 
\{ \xb : \cb(\xb)
\le \zerob, \anlchange{\lb} \le \xb \le \ub \}$, \lineref{generate} of
\algref{pure_random} may involve drawing points uniformly at random from $[\anlchange{\lb}, \ub]$ and
checking whether they satisfy $\cb(\xb) \le \zerob$.
\begin{algorithm}[tb]
 Choose initial point $\hat{\xb} \in \Omegab$, termination test, and point generation scheme
 
 \While{Termination test is not satisfied}
 {
 Generate $\xb \in \Omegab$\label{line:generate}\\

 \If {$f(\xb) < f(\hat{\xb})$}
 {
 $\hat{\xb} \gets \xb$
 }
 }
 \caption{Pure random search \label{alg:pure_random}}
\end{algorithm}
If the procedure for generating points in \lineref{generate} of
\algref{pure_random} is independent of the function values observed, then the
entire set of points used within pure random search can be generated
beforehand: we intentionally omit the index $k$ in the statement of \algref{pure_random}.

Nevertheless, for the sake of analysis, it is useful to consider an ordering of
the sequence of random points generated by \algref{pure_random}. 
With such a sequence $\{\xb_k\}$, one can analyse the best points after $N$ evaluations,
\[
\hat{\xb}_N \in \argmin_{k=1,\ldots,N} f(\xb_k).
\]
If $f_*$ is the global minimum value, then 
\[
\P{f(\hat{\xb}_N) \le f_* + \epsilon} = 1 - 
\prod_{k=1}^{N}(1-\P{\xb_k \in
\cL_{f_* + \epsilon}(f)}),
\]
where $\epsilon \ge 0$ and $\cL_{\alpha}(f) = \{ \xb : f(\xb) \le \alpha\}$.
Provided that the procedure used to generate points at \lineref{generate} of
\algref{pure_random} satisfies 
\[
\lim_{N \to \infty} \prod_{k=1}^{N}(1-\P{\xb_k \in \cL_{f_* +
\epsilon}(f)}) = 0
\]
for all $\epsilon > 0$, then $f(\hat{\xb}_k)$ converges in probability to
$f_*$. 
For example, if each $\xb_k$ is drawn independently and uniformly over 
$\Omegab$, then one 
can
calculate the number of evaluations required to ensure that the $\hat{\xb}_k$ 
returned by \algref{pure_random}
satisfies $\hat{\xb}_k \in \cL_{f_* +
\epsilon}$ with probability $p \in (0,1)$,
that is, 
\[N \ge \ffrac{\log(p)}{\log\biggl(1 - \ffrac{\mu(\cL_{f_* + 
\epsilon} \bigcap \Omegab )}{\mu(\Omegab)}\biggr)},\]
provided $\mu (\cL_{f(\xb_*) + \epsilon} \cap \Omegab ) > 0$ and
$\Omegab$ is measurable.

Random-search methods typically make few assumptions about 
$f$; see \citeasnoun{Zhigljavsky1991} for further
discussion about the convergence of pure random search.
Naturally, a method that assumes only that 
$f$ is measurable on $\Omegab$ is likely to produce \anlchange{function values} that converge
more slowly to \anlchange{$f_*$} when applied to an $f \in \cC^0$ than does a method that
exploits the continuity of $f$.
Heuristic modifications of random search have sought to
improve empirical performance on certain classes of
problems, while still maintaining random search's global 
optimization property;
see, for example, the work of \citeasnoun{Zabinsky1992} and \citeasnoun{Patel1989}.

\subsubsection{Nesterov random search}
\label{sec:rand_det_rs}
We refer to the method discussed in this section as Nesterov random
search because of the seminal article by \citeasnoun{Nesterov2015}, but the idea 
driving this method is much older. 
A similar method, for instance, is discussed in \citeasnoun[Chapter~3.4]{polyakbook}. 

The method of Nesterov random search is largely motivated by Gaussian smoothing.
In particular, given a covariance (\ie\ symmetric positive-definite) 
matrix $\Bb$ and a
smoothing parameter $\mu>0$, consider the Gaussian smoothed function 
\[
f_{\mu}(\xb) =
\sqrt{\ffrac{\textrm{det}(\Bb)}{(2\pi)^n}} \int_{\Reals^n} f(\xb +
\mu\ub)\exp\biggl(-\ffrac{1}{2}\ub^\T \Bb \ub\biggr) \D \ub.
\]
This smoothing
has many desirable properties; for instance, if $f$ is
Lipschitz-continuous with constant $\Lip{f}$, then $f_{\mu}$ is
Lipschitz-continuous with a constant no worse than $\Lip{f}$ for all $\mu>0$. 
Likewise,
if $f$ has Lipschitz-continuous gradients with constant $\Lip{g}$, then $f_\mu$ 
has Lipschitz-continuous gradients with a constant no worse than $\Lip{g}$
for all
$\mu>0$. 
If $f$ 
is convex, then $f_\mu$ is convex. 

One can show that
\[
\nabla f_{\mu}(\xb) = \ffrac{1}{\mu}\sqrt{\ffrac{\textrm{det}(\Bb)}{(2\pi)^n}}
\int_{\Reals^n} (f(\xb + \mu\ub)-f(\xb))\exp\biggl(-\ffrac{1}{2}\ub^\T \Bb \ub\biggr)\Bb\ub \D \ub.
\]
In other words, $\nabla f_{\mu}(\xb)$, 
which can be understood as an approximation of $\nabla f(\xb)$ in the smooth case,
can be computed via an expectation over $\ub\in\Reals^n$ weighted by the 
finite difference $f(\xb + \mu\ub)-f(\xb)$ and inversely weighted by a radial distance from $\xb$. 
With this interpretation in mind, \citeasnoun{Nesterov2015} propose a collection of 
\emph{random gradient-free oracles}, where one first generates a Gaussian 
random vector $\ub\in\cN (\zerob,\Bb^{-1})$ 
and then uses one of 
\begin{equation}
\label{gradient_free_oracle}
\begin{aligned}
\gb_{\mu}(\xb;\ub) & \defined 
\df{f;\xb;\ub;\mu}
\Bb\ub,
\quad \mbox{or}\\
\hat{\gb}_{\mu}(\xb;\ub) & \defined 
\dc{f;\xb;\ub;\mu}
\Bb\ub,
\end{aligned}
\end{equation}
for a difference parameter $\mu > 0$.
\citename{Nesterov2015} also propose a third oracle, 
$\gb_0(\xb;\ub) \defined f'(\xb;\ub) \Bb\ub$, 
intended for the optimization of non-smooth functions;
this oracle assumes the ability to compute directional
derivatives $f'(\xb;\ub)$.
For this reason, \citename{Nesterov2015} refer to all oracles as gradient-free 
instead of derivative-free. 
Given the scope of this survey,
we focus on 
the derivative-free oracles
$\gb_\mu$ and $\hat{\gb}_\mu$
displayed in \eqref{gradient_free_oracle}. 

With an oracle $\gb$ chosen as either oracle in \eqref{gradient_free_oracle},
Nesterov random-search methods are straightforward to
define, and we do so in \algref{rs-nesterov}. In \algref{rs-nesterov},
$\Proj{\cdot;\Omegab}$ denotes projection onto a 
domain $\Omegab$. 

\begin{algorithm}[tb]
 Choose initial point $\xb_0\in \Omegab$, sequence of step sizes
 $\{\alpha_k\}_{k=0}^{\infty}$, oracle $\gb$ from
 \eqref{gradient_free_oracle}, smoothing parameter $\mu>0$ and covariance
 matrix $\Bb$ \\

 \For{$k=0,1,2,\ldots$}
 {
 $\hat{\xb}_k \gets \argmin_{j\in \{0,1,\ldots,k\}} f(\xb_j)$\\
 Generate $\ub_k\in\cN (\zerob,\Bb^{-1})$;
 compute $\gb(\xb_k;\ub_k)$\\
 $\xb_{k+1} \gets \Proj{\xb_k-\alpha_k\Bb^{-1}\gb(\xb_k;\ub_k),\Omegab}$
 }
 \caption{Nesterov random search \label{alg:rs-nesterov}}
\end{algorithm}

A particularly striking result proved in \citeasnoun{Nesterov2015} was perhaps the 
first WCC result
for an algorithm (\algref{rs-nesterov}) in the case where
$f\in\cLC^0$ -- that is, $f$ may be both non-smooth and non-convex. 
Because of the randomized nature of iteratively sampling from 
a Gaussian distribution in \algref{rs-nesterov},
complexity results are given as expectations.
That is, letting $\Ub_k = \{\ub_0,\ub_1,\ldots,\ub_k\}$ denote the random variables associated with
the first $k$ iterations of \algref{rs-nesterov}, complexity results are stated in terms of expectations
with respect to the filtration defined by these variables. 
A WCC is given as an upper bound on the number of $f$ evaluations needed to 
attain the approximate ($\epsilon >0$) optimality condition
\begin{equation}
 \label{eq:gauss_nonconvex_complexity}
 \E[\Ub_{k-1}]{\|\nabla f_{\check{\mu}}(\hat{\xb}_k)\|}\leq \epsilon, 
\end{equation}
where $\hat{\xb}_k = \argmin_{j=0,1,\ldots,k-1} f(\xb_j)$.
By fixing a particular choice of $\check{\mu}$ (dependent on $\epsilon$, $n$ 
and Lipschitz constants), 
\citeasnoun{Nesterov2015} demonstrate that the number of $f$ evaluations needed to 
attain 
\eqref{eq:gauss_nonconvex_complexity} is in $\bigo{\epsilon^{-3}}$; see \tabref{rates}. 
For $f\in\cLC^1$ (but still non-convex), \citeasnoun{Nesterov2015} prove a WCC result of type
\eqref{eq:gauss_nonconvex_complexity} in $\bigo{\epsilon^{-2}}$ for the same method.
WCC results of \algref{rs-nesterov} under a variety of stronger assumptions on
the convexity and differentiability of $f$ 
are also shown in \tabref{rates} and discussed in \secref{convex}.
We further note that some randomized methods of the form \algref{rs-nesterov}
have also been developed for \eqref{eq:stoch_prob}, which we discuss in
\secref{stoch}.

We remark on an undesirable feature of the convergence analysis for variants of
\algref{rs-nesterov}: the analysis of these methods supposes that the sequence $\{\alpha_k\}$ is chosen as a constant
that \anlchange{depends on parameters, including $\Lip{f}$, that may not be available to the method}. 
Similar assumptions concerning the preselection of $\{\alpha_k\}$ also appear in 
the convex cases discussed in \secref{convex}, and we highlight these
dependencies in \tabref{rates}.

\subsection{Randomized direct-search methods}
\label{sec:rand_det_ds}

Randomization has also been used in the
\framework{DDS} framework discussed in \secref{det_det_DS}
in the hope of more efficiently using evaluations of $f$.
Polling every point in a PSS requires at least $n+1$ function evaluations;
if $f$ is expensive to evaluate and $n$ is relatively large, this can be
wasteful. A deterministic strategy for performing fewer 
evaluations on many iterations is opportunistic polling.
Work in randomized direct-search methods attempts to address, formalize and analyse
the situation where polling directions are \emph{randomly sampled} from some distribution in each iteration.
The ultimate goal is to replace the $\bigo{n}$ per-iteration function evaluation cost 
with an $\bigo{1}$ per-iteration cost,\footnote{We note that a $\bigo{1}$ cost 
can naturally be achieved by deterministic \framework{DDS} methods when derivatives are available 
\cite{Abramson2004gps}.} while still guaranteeing some form of 
global convergence.

In \secref{det_det_DS}, we mentioned \codes{MADS} methods 
that consider the random generation of polling directions in each
iteration 
(in order to satisfy the asymptotic density required of 
search directions
for the minimization of non-smooth, but Lipschitz-continuous, $f$).
Examples include \citeasnoun{Audet06mads} and \citeasnoun{VanDyke2013}, which implement
\codes{LTMADS} and \codes{QRMADS}, respectively. 
While this direction of research is within the scope of randomized methods, 
the purpose of randomization in \codes{MADS} methods is to overcome
particular difficulties encountered when optimizing general non-smooth objectives.
 This particular randomization does not fall within the scope of this section,
 where randomization is intended to decrease a method's dependence on $n$. 
 In the remainder of this section, we focus on a body of work that seems to exist entirely for the 
unconstrained case
where $f$ is assumed sufficiently smooth. 

\citeasnoun{Gratton2015}
extend the direct-search framework (\algref{ds})
by assuming that the set of polling directions $\Db_k$ includes only a descent
direction with probability $p$ 
(as opposed to assuming $\Db_k$ always includes a descent direction, which comes for 
free when $f\in\cLC^1$ provided $\Db_k$ is, for example, a PSS).
To formalize,
given $p\in(0,1)$, 
a random sequence of polling directions $\{\Db_k\}$ is said to be 
$p$-probabilistically $\kappad$-descent provided that,
given a deterministic starting point $\xb_0$, 
\begin{equation}\label{eq:descent_part1}
 \P{\cm ( [\Db_0,-\nabla f(\xb_0) ] ) \geq \kappad} \geq p,
\end{equation}
and for all $k\geq 1$,
\begin{equation}\label{eq:descent_part2}
 \P{\cm ( [\Db_k,-\nabla f(\xb_k) ] ) \geq \kappad \mid 
 \Db_0,\ldots,\Db_{k-1}} \geq p,
\end{equation}
where $\cm(\cdot)$ is the cosine measure in \eqref{eq:cosine}. 
A collection of polling directions $\Db_k$ satisfying 
 \eqref{eq:descent_part1} and \eqref{eq:descent_part2} 
can be obtained by drawing directions uniformly on the unit ball.

As with the other methods in this section, $\xb_k$ in \eqref{eq:descent_part2} 
is in fact a random variable due to the random sequence 
$\{\Db_k\}$ generated by the algorithm, and hence it 
makes sense to view \eqref{eq:descent_part2} as a
probabilistic statement. 
In words, \eqref{eq:descent_part2} states that with probability at
least $p$, 
the set of polling directions used in iteration $k$ has a positive cosine measure with
the steepest descent direction $-\nabla f(\xb_k)$, regardless of the past history of the algorithm.
\citeasnoun{Gratton2015} use Chernoff bounds in order to bound the worst-case complexity 
\emph{with high probability}.
Roughly, they show that if $f\in\cLC^1$, and if
\eqref{eq:descent_part1} and $\eqref{eq:descent_part2}$ hold with $p > 1/2$ (this constant changing when $\gammad \neq 1/\gammai$),
then
$\|\nabla f(\xb_k)\| \leq \epsilon$ holds within $\bigo{\epsilon^{-2}}$ function evaluations with 
a probability
that increases exponentially to 1 as $\epsilon\to 0$. 
The WCC result of \citeasnoun{Gratton2015} demonstrates that as $p\to 1$ (\ie\ 
$\Db_k$ almost always includes a descent direction),
the known WCC results for \algref{ds} discussed in \secref{det_det_DDS} are recovered. 
A more precise statement of this WCC result is included in \tabref{rates}.

\anlchange{\citeasnoun{Bibi2019}
propose a randomized direct-search method in which the two poll directions in each iteration
are $\Db_k=\{\eb_i,-\eb_i\}$, where $\eb_i$ is the $i$th elementary basis vector.
In the $k$th iteration, $\eb_i$ is selected from $\{\eb_1,\dots,\eb_n\}$ 
with a probability proportional to the Lipschitz constant of the $i$th partial
derivative of $f$.
\citeasnoun{Bibi2019} perform WCC analysis of this method assuming 
a known upper bound on Lipschitz constants of partial derivatives;
this assumption leads to improved constant factors, but they
essentially prove an upper bound on the number of iterations needed to attain
$\E[]{\|\nabla f(\xb_k)\|}\leq \epsilon$ in $\bigo{\epsilon^{-2}}$, where the expectation
is with respect to the random sequence of $\Db_k$. 
\citeasnoun{Bibi2019} prove additional WCC results
in cases where $f$ is convex or $c$-strongly convex.} 

An early randomized \framework{DDS} derivative-free method that 
only occasionally employs a descent direction is developed by 
\citeasnoun{DinizEhrhardt2008}. 
There, a non-monotone line-search strategy is used to accommodate 
search directions along which descent may not be initially apparent.
\citeasnoun{Belitz2013} develop a randomized \framework{DDS} method that 
employs surrogates and an adaptive lattice.

\subsection{Randomized trust-region methods}
\label{sec:rand_det_tr}

Whereas the theoretical convergence of 
a \framework{DDS} method depends on the set of polling directions satisfying some spanning property
(\eg\ a cosine measure bounded away from zero),
the theoretical convergence of a trust-region method (\eg\ \algref{tr}) 
depends on the use of fully linear models. 
Analogous to how randomized \framework{DDS} methods relax the requirement of the use of a positive spanning
set in \emph{every} iteration, \eqref{eq:descent_part2},
it is reasonable to ask whether one can relax the requirement of being fully 
linear in every iteration of a trust-region method.
Practically speaking, in the unconstrained case it may not be necessary to ensure that every model is built by using a
$\Lambda$-poised set of points (therefore ensuring that the model is fully linear) on every iteration,
since ensuring $\Lambda$-poised sets entails additional function evaluations. 

\citeasnoun{Bandeira2014} 
consider a sequence of random models $\{m_k\}$
and a random sequence of  \anlchange{trust-region centres and} radii $\{\xb_k,\Delta_k\}$.\ 
They say that the sequence of random models is $p$-probabilistically 
\anlchange{$\kappab$}-fully linear provided
\begin{equation}\label{eq:tr_prob}
 \P{ m_k \text{ is a } \anlchange{\kappab}\text{-fully linear model of $f$ on } 
\B{\xb_k}{\Delta_k} \anlchange{ \; \big| \; \cH_{k-1}} } \geq p,
\end{equation}
where $\anlchange{\cH_{k-1}}$ is the filtration of the random process prior to the current iteration. That is, $\anlchange{\cH_{k-1}}$ is the $\sigma$-algebra generated by the algorithm's history. 
Under additional standard assumptions concerning \algref{tr} (\eg\ $\gammad = 1/\gammai$), the authors show that if \eqref{eq:tr_prob} holds with $p > 1/2$, then 
$\lim_{k\rightarrow \infty} \|\nabla f(\xb_k)\|\anlchange{=}0$ almost surely (\ie\ with probability one).
\citeasnoun{Gratton2017complexity} build on this result; they demonstrate that, up to
constants, the same (with high probability) WCC bound that was proved for \framework{DDS} methods in
\citeasnoun{Gratton2015} holds for the randomized trust-region method proposed by
\citeasnoun{Bandeira2014}. 
\anlchange{Higher-order versions of \eqref{eq:tr_prob} also exist; in 
\secref{sepspa} we discuss settings for which \citeasnoun{Bandeira2014} obtain 
probabilistically 
$\kappab$-fully quadratic models by interpolating $f$ on a set of fewer than 
$\qn$
points.}

\vspace{5pt}
\section{Methods for convex objectives} %
\label{sec:convex}
As is true of derivative-based optimization, convexity in the objective of 
\eqref{eq:det_prob} or \eqref{eq:stoch_prob} 
can be exploited either when designing new methods or when analysing 
existing methods. 
Currently, this split falls neatly into two categories.
The majority of work considering \eqref{eq:det_prob} when $f$ is convex
sharpens the WCCs for frameworks already discussed
in this survey.
On the other hand, the influence of machine learning, particularly large-scale
empirical risk minimization,
has led to entirely new derivative-free methods for solving \eqref{eq:stoch_prob}
when $f$ is convex.

\subsection{Methods for deterministic convex optimization}
\label{sec:convex_det}

We first turn our attention to the solution of \eqref{eq:det_prob}. 
In convex optimization,
one can prove WCC bounds on the difference between a method's estimate of the global minimum of $f$ and the value of $f$
at a global minimizer $\xb_*\in\Omegab$ (\ie\
a point satisfying $f(\xb_*)\leq f(\xb)$ for all $\xb\in\Omegab$). 
This differs from the local WCC bounds on the objective gradient, namely
\eqref{eq:first_order_convergence}, that are commonly shown when $f$
is not assumed to be convex. 
For convex $f$, under appropriate additional assumptions,
one typically demonstrates that a method satisfies
\begin{equation}
\label{eq:zeroth_order_convergence}
\lim_{k\to\infty} f(\xb_k)-f(\xb_*) = 0,
\end{equation}
where $\xb_k$ is the $k$th point of the method. 
Hence, 
an appropriate measure of $\epsilon$-optimality
when minimizing an unconstrained convex objective is the satisfaction of
\begin{equation}
\label{eq:convex_stationarity}
f(\xb_k)-f(\xb_*)\leq\epsilon.
\end{equation} 

For the \framework{DDS} methods discussed in \secref{det_det_DDS},
\citeasnoun{Dodangeh2016} and \citeasnoun{Dodangeh2015} analyse the worst-case complexity
of \algref{ds} when there is no search step and when 
$f$ is smooth and convex, 
and has a bounded level set. 
By imposing an appropriate upper bound on the step sizes $\alpha_k$,
and (for $c>0$) using a test of the form 
\begin{equation}
\label{ds_suff_dec}
f(\pb_i)\leq f(\xb_k) - c\alpha_k^2,
\end{equation}
in \lineref{descent} of
\algref{ds},
\citeasnoun{Dodangeh2016} show that 
the worst-case number of $f$ evaluations to achieve 
\eqref{eq:nonconvex_stationarity}
is in $\bigo{n^2\Lip{g}^2\epsilon^{-1}}$. 
\citeasnoun{Dodangeh2015} show that this $n^2$-dependence is optimal 
(in the sense that it cannot be improved) 
within the class of deterministic methods that employ positive spanning sets. 
Recall from \secref{rand_det_ds} that randomized methods allow one to reduce this 
dependence to be linear in $n$.
Under additional assumptions, which are satisfied, for example, when $f$ is strongly convex,
\citeasnoun{Dodangeh2016} prove $R$-linear convergence of \algref{ds}, yielding a WCC
of type \eqref{eq:convex_stationarity} with the dependence on $\epsilon$ 
reduced to $\log(\epsilon^{-1})$. We note that, in the convex setting, $R$-linear convergence had been previously established for a \framework{DDS} method by
 \citeasnoun{Dolan2003}.
It is notable that the method analysed in \citeasnoun{Dolan2003}
requires only strict decrease, whereas, to the authors' knowledge,
the \framework{DDS} methods for which
WCC results have been established all require sufficient decrease.

\citeasnoun{KP2014} propose a \framework{DDS} method
that does not allow for increases in the step size $\alpha_k$
and analyse the method on strongly convex, convex and non-convex
objectives. Although \citeasnoun{KP2014} demonstrate WCC bounds with
the same dependence on $n$ and $\epsilon$ as do \citeasnoun{Dodangeh2016}, they
additionally assume that one has explicit knowledge of a Lipschitz gradient constant $\Lip{g}$ 
and can thus replace the test \eqref{ds_suff_dec} explicitly with
\begin{equation}
\label{ds_suff_dec2}
f(\pb_i)\leq f(\xb_k) - \ffrac{\Lip{g}}{2}\alpha_k.
\end{equation}
Exploiting this additional knowledge of $\Lip{g}$, the WCC result in \citeasnoun{KP2014}
exhibits a strictly better dependence on $\Lip{g}$ than does the WCC result in \citeasnoun{Dodangeh2016},
with a WCC of type \eqref{eq:convex_stationarity} in $\bigo{n^2\Lip{g}\epsilon^{-1}}$. 
Additionally assuming $f$ is
$c$-strongly convex, \citeasnoun{KP2014} provide a result showing a WCC
of type \eqref{eq:convex_stationarity} in $\bigo{\log(\epsilon^{-1})}$.

In the non-convex
case, \citeasnoun{KP2014} recover the same WCC
of type \eqref{eq:nonconvex_stationarity} from \citeasnoun{Vicente2013}; see the
 discussion in \secref{det_det_DS}.
Once again, however, the result by \citeasnoun{KP2014} assumes knowledge of
$\Lip{g}$ and again recovers a 
strictly better dependence on $\Lip{g}$ by using the test \eqref{ds_suff_dec2} in \lineref{descent}
of \algref{ds}. 

Recalling the Nesterov random search methods discussed in \secref{rand_det_rs}, 
we remark that \citeasnoun{Nesterov2015} explicitly
give results for deterministic, convex $f$. 
In particular, \citename{Nesterov2015} prove WCCs of a specific type.
Because of
the randomized nature of \algref{rs-nesterov}, WCCs are given as
expectations of the form
\begin{equation} 
 \label{eq:gauss_convex_stationarity}
 \E[\Ub_{k-1}]{f(\hat{\xb}_k)}-f(\xb_*) \leq\epsilon, 
\end{equation} 
where $\hat{\xb}_k = \argmin_{j \in \{0,1,\ldots,k-1\}} f(\xb_j)$ and where $\Ub_{k-1} =
\{\ub_0,\ub_1,\ldots,\ub_{k-1}\}$ is the filtration of
Gaussian samples. 
The form of $\epsilon$-optimality represented by \eqref{eq:gauss_convex_stationarity} 
can be interpreted as a probabilistic variant of \eqref{eq:convex_stationarity}.
 The WCC results of \citename{Nesterov2015} show that the worst-case number of $f$
 evaluations to achieve \eqref{eq:gauss_convex_stationarity} is in
 $\bigo{\epsilon^{-1}}$ when $f\in\cLC^1$.
 Additionally assuming that $f$ is $c$-strongly convex yields an improved result;
 the WCC of type \eqref{eq:gauss_convex_stationarity}
 is now in $\bigo{\log(\epsilon^{-1})}$. 
 Moreover, by mimicking the method of accelerated gradient descent 
 (see \eg\ \ANL{\citeb{Nesterov2013introductory}, Chapter~2.2}), 
 \citename{Nesterov2015}
 present a variant of \algref{rs-nesterov} with a WCC of type
 \eqref{eq:gauss_convex_stationarity} in $\bigo{\epsilon^{-1/2}}$. 
 When $f\in\cLC^0$, \citename{Nesterov2015} provide a 
 WCC of type \eqref{eq:gauss_convex_stationarity} in $\bigo{\epsilon^{-2}}$, but this result
 assumes that \algref{rs-nesterov} uses an oracle with access to exact directional derivatives of $f$.
 Thus, the method achieves the $\bigo{\epsilon^{-2}}$ result when $f\in\cLC^0$ 
 is not a derivative-free method.
 
As remarked in \secref{rand_det}, these convergence results
depend on preselecting a
sequence of step sizes $\{\alpha_k\}$ for \algref{rs-nesterov}; 
in the convex case, the necessary $\{\alpha_k\}$ depends not only on
the Lipschitz constants but also on a bound $\anlchange{\Rx}$ on the
distance between the initial point and the global minimizer (\ie\
$\|\xb_0-\xb_*\|\leq \anlchange{\Rx}$). 
The aforementioned WCC results will hold only if one chooses $\{\alpha_k\}$ and $\mu$ 
(the difference parameter of the oracle used in \algref{rs-nesterov})
that scale with $\Lip{g}$ 
and $\anlchange{\Rx}$ appropriately.
When additionally assuming $f$ is $c$-strongly convex, 
$\{\alpha_k\}$ and $\mu$ also depend on $c$.
\citeasnoun{Stich2011}\footnote{A careful reader may be caught off guard by the fact 
that
\citeasnoun{Stich2011} was published before \citeasnoun{Nesterov2015}. This
is not a typo; 
\citeasnoun{Stich2011} build on the results from
an early preprint of \citeasnoun{Nesterov2015}.}
 extend the framework of 
\algref{rs-nesterov} with an approximate line search
that avoids the need for predetermined sequences of step sizes.

In general, the existing WCC results for derivative-free methods match
the WCC results for their derivative-based counterparts in 
 $\epsilon$-dependence.
The WCC results for derivative-free methods 
tend to involve an additional factor of $n$ when compared with their derivative-based counterparts.
This observation mirrors a common intuition in derivative-free optimization:\ 
since a number of $f$ evaluations in $\bigo{n}$ can
guarantee a suitable approximation of the gradient, then for any class of
gradient-based method for which we replace gradients with approximate gradients or model gradients, one
should expect to recover the WCC of that method, but
with an extra factor of $n$. 
WCC results such as those discussed thus far in this survey
add credence to this
intuition. As we will see in \secref{stoch_convex_bandit}, however, this
optimistic trend does not always hold:\ we will see problem classes
for which derivative-free methods are provably worse than their
derivative-based counterparts by a factor of $\epsilon^{-1}$. 

\citeasnoun{Bauschke2014} offer an alternative approach to \algref{rs-nesterov} for
the solution of \eqref{eq:det_prob}
when $f$ is assumed convex and, additionally, lower-$\cC^2$.
(Such an assumption on $f$ is
obviously stronger than plain convexity but contains, for example, functions
that are defined as the
pointwise maximum over a collection of convex functions.)
\citeasnoun{Bauschke2014} show that linear interpolation through function values is
sufficient for obtaining approximate subgradients of convex, lower-$\cC^2$ function; these approximate subgradients are used in
lieu of subgradients in a mirror-descent algorithm 
(see \eg\ \citeb{Srebro2011}) 
similar to 
\algref{2ptmirror} in \secref{two-point-feedback}. \citeasnoun{Bauschke2014} establish convergence of their method
in the sense of \eqref{eq:zeroth_order_convergence},
 and they demonstrate the performance of their method when applied to pointwise
maxima of collections of convex quadratics. 

\subsection{Methods for convex stochastic optimization} 
\label{sec:stoch_convex_bandit}
We now turn our attention to the solution of the problem \eqref{eq:stoch_prob}
when $\tilde{f}(\xb;\xib)$ is assumed convex in $\xb$ for each realization $\xib$.
Up to this point in the survey, we have typically assumed
that \eqref{eq:stoch_prob} is unconstrained, that is,
 $\Omegab=\Reals^n$. 
In this section, however, it will become more frequent that
$\Omegab$
is a compact set.

In the machine learning community, zeroth-order information
(\ie\ evaluations of $\tilde{f}$ only)
 is frequently referred
to as \emph{bandit feedback},\footnote{A `one-armed bandit' is an American colloquialism for a casino slot machine,
the arm of which must be pulled to reveal a player's losses or rewards.}
due to the concept of multi-armed bandits from
reinforcement learning. 
Multi-armed bandit problems are sequential allocation
problems defined by a prescribed set of actions. 
\citeasnoun{Robbins1952} formulates
a multi-armed bandit problem as a gambler's desire to minimize
the total losses accumulated
 from pulling discrete sequence (of length $T$) of $A<\infty$ slot machine
 arms. 
 The gambler does not have to decide the full length-$T$ sequence up front. Rather,
 at time $k\leq T$,
 the losses associated with the first $k-1$ pulls are known to the gambler
 when deciding which of the $A$ arms to pull next. 
 The gambler's decision of the $k$th arm to pull is represented 
 by the scalar variable $x_k\in\Omegab$.
 Given additional environmental
 variables outside the gambler's control, $\xib_k\in\Xib$, 
 which represents the stochastic nature of the slot machine,\footnote{Depending
 on the problem set-up, the environment of a bandit problem
 may be either stochastic or adversarial. Because this section is discussing stochastic
 convex optimization, we will assume that losses are stochastic; that is,
 the $\xib_k$ are i.i.d.\ and
 independent of the gambler's decisions $x_k$.
 See \citeasnoun{BubeckCesaBianchi2012}
 for a survey of bandit problems more general than those discussed
 in this survey.}
 the environment makes a decision $\xib$
 simultaneously with the gambler's decision $x$.
 The gambler's loss is then $\tilde{f}(x;\xib)$. 

 Within this multi-armed bandit setting,
 the typical metric of the gambler's performance is
 the \emph{cumulative regret}.
 Provided that the expectation 
 \anlchange{$\Ea[\xib]{\tilde{f}(x;\xib)} = f(x)$} exists for each
 $x\in\{1,\dots,A\}$,
 then the gambler's best long-run strategy in terms of minimizing the
expected total losses 
 is to constantly play
 $x_*\in \argmin_{x\in\{1,\ldots,A\}} f(x)$.
 If, over the course of $T$ pulls, the gambler 
 makes a sequence of decisions $x_1,\dots,x_T$ and
 the environment makes a sequence of decisions $\xib_1,\dots,\xib_T$
 resulting in a sequence of realized losses
 $\tilde{f}(x_1;\xib_1),\dots, \tilde{f}(x_T;\xib_T)$,
 then the 
 cumulative regret $\anlchange{r_T}$ associated with the gambler's sequence of decisions is
the difference between the cumulative loss
 incurred by the gambler's strategy
$(x_1,\ldots,x_T)$ and the loss incurred by the best possible long-run strategy
$(x_*,\ldots,x_*)$. 
Analysis of methods for bandit problems in this set-up is generally concerned with \emph{expected cumulative regret} 
\begin{equation}
 \label{expected_regret}
 \E[\xib]{\anlchange{r_T}(x_1,\dots,x_T)} \defined \Egg[\xib]{\sum_{k=1}^{T} \tilde{f}(x_k;\xib_k)}-Tf(x_*),
\end{equation}
where the expectation is computed over $\xib\in \Xib$, 
since we assume here that the sequence of $\xib_k$ is 
independent and identically distributed.

 This particular treatment of the bandit problem
 with a discrete space of actions $x\in\{1,\ldots,A\}$ 
 was the one considered by \citeasnoun{Robbins1952}
 and has been given extensive treatment
 \cite{ACFS1995,LaiRobbins1985,Agarwal1995,ACF2002}.
 
 Extending multi-armed bandit methods to \emph{infinite-armed} 
 bandits makes the connections to derivative-free optimization -- particularly 
 derivative-free convex optimization -- readily apparent. 
 \citeasnoun{Auer2002} extends the multi-armed bandit problem to allow
 for a compact (as opposed to discrete) set of actions for the gambler
 $\xb\in\Omegab\subset\Reals^n$ as well as a compact set of vectors
 for the environment, 
 $\xib\in\Xib\subset\Reals^n$. 
 The vectors $\xib$ in this set-up define linear functions;
 that is, if the gambler chooses $\xb_k$ in their $k$th pull,
 and the environment chooses $\xib_k\in\Xib$,
 then the gambler incurs loss $\tilde{f}(\xb_k;\xib_k)=\xib_k^{\T}\xb_k$.
 In this linear regime, expected regret takes a 
 form remarkably similar to \eqref{expected_regret}, that is, 
 \begin{equation} \label{expected_regret_convex}
 \begin{aligned}
 \E[\xib]{\anlchange{r_T}(\xb_1,\dots,\xb_T)} & = \Egg[\xib]{ \sum_{k=1}^T 
 \tilde{f}(\xb_k;\xib_k)}- \min_{\xb\in\Omegab}\Egg[\xib]{ \sum_{k=1}^T 
 \tilde{f}(\xb_k;\xib_k)}\\
 & = \Biggl( \sum_{k=1}^T f(\xb_k)\Biggr) - Tf(\xb_*),
 \end{aligned}
 \end{equation}
 with $\xb_*\in \argmin_{\xb\in\Omegab} f(\xb)$. 
As in \eqref{expected_regret}, \eqref{expected_regret_convex} defines the 
expected regret with 
respect to the best long-run strategy $\xb\in\Omegab$ that the gambler could 
have played for the $T$ 
rounds (\ie\ the strategy that would minimize the expected cumulative 
losses).\footnote{We remark that many of the references we provide here may also refer to 
methods minimizing \eqref{expected_regret_convex} as methods of \emph{online optimization};
one reference in particular \cite{pmlr-v49-bach16} even suggests a taxonomy identifying
bandit learning as a restrictive case of online optimization.}

By using a bandit method known as Thompson sampling, 
one can show (under appropriate additional assumptions) that if
$\tilde{f}(\xb;\xib) = \xib^{\T}\xb$, then
\eqref{expected_regret_convex} can be bounded as
\anlchange{$\Ea[\xib]{\anlchange{r_T}(\xb_1,\dots,\xb_T)}\in\bigo{n\sqrt{T\log(T)}}$}  
\cite{RussoVanRoy2014}.
Analysis of bounds on \eqref{expected_regret_convex} in the linear case raises an interesting question:
to what extent can 
similar analysis be performed
for classes of functions $\tilde{f}(\xb;\xib)$ that are
non-linear in $\xb$?

In much of the bandit literature, the additional structure defining a class of non-linear functions is convexity;
here, by convexity, we mean that $\tilde{f}(\xb;\xib)$ is convex in $\xb$ for 
each realization
$\xib\in\Xib$. 

Regret bounds on \eqref{expected_regret_convex} automatically imply WCC results
for stochastic convex optimization.
To see this, define an \emph{average point}
\[
\bar{\xb}_k = \ffrac{1}{k} \sum_{t=1}^k \xb_t.
\]
Because of the convexity of $f$,
\begin{equation}
\label{eq:regret_implies_opt}
f(\bar{\xb}_T) - f(\xb_*) \leq \ffrac{1}{T} \sum_{k=1}^T f(\xb_k)-f(\xb_*)
=\ffrac{\E[\xib]{\anlchange{r_T}(\xb_1,\dots,\xb_T)}}{T},
\end{equation}
where the inequality follows from an application of Jensen's inequality.
We see in \eqref{eq:regret_implies_opt} that, given an upper bound of 
$\bar{\anlchange{r}}(T)$ on the
expected regret \anlchange{$\Ea[\xib]{\anlchange{r_T}}$}, we can automatically
derive, provided $\bar{\anlchange{r}}(T)/T \leq\epsilon$, a WCC result for stochastic convex optimization of the form
\begin{equation}
\label{convex_stoch_stationarity}
\E[\xib]{\tilde{f}(\bar{\xb}_{T};\xib) - \tilde{f}(\xb_*;\xib)}
 = f(\bar{\xb}_{T})-f(\xb_*)\leq \epsilon.
\end{equation}
Equation \eqref{convex_stoch_stationarity} corresponds to a stochastic version of a 
WCC result
of type \eqref{eq:convex_stationarity} 
where $N_{\epsilon}=T$.

The converse implication, however, is not generally true: a
small optimization error does not imply a small regret.
That is, WCC results of type \eqref{convex_stoch_stationarity} 
do not generally imply bounds on expected regret \eqref{expected_regret_convex}.
This is particularly highlighted by \citeasnoun{Shamir13}, who considers a class of 
strongly convex quadratic objectives $f$.
 For such objectives, \citeasnoun{Shamir13} establishes
a strict gap between
the upper bound on the optimization error (the left-hand side of \eqref{convex_stoch_stationarity})
and the lower bound on the expected regret \eqref{expected_regret_convex}
attainable by any method in the bandit setting.
Such a gap, between expected regret and optimization error,
has also been proved for bandit methods applied
to problems of logistic regression \cite{Hazan2014}. 
 
Results such as those of \citeasnoun{Shamir13}, \citeasnoun{Jamieson2012} and 
\citeasnoun{Hazan2014}
have led researchers to %
 consider
methods of derivative-free (stochastic) optimization 
within the paradigm of convex bandit learning.
In this survey, we group these
methods into two categories
of assumptions on the type of bandit feedback (\ie\ the observed realizations of $\tilde{f}$)
available to the method:
one-point bandit feedback or two-point (multi-point) bandit
feedback. Although many of the works we cite have also analysed regret bounds for these methods,
we focus on optimization error. 
 
 \subsubsection{One-point bandit feedback}
 \label{sec:one-point-feedback}
In one-point bandit feedback, a method is assumed to have 
access to an oracle $\tilde{f}$ that returns unbiased estimates of $f$. 
In particular, given two
points $\xb_1$, $\xb_2\in\Omegab$, two separate calls to the oracle
will return $\tilde{f}(\xb_1;\xib_1)$ and $\tilde{f}(\xb_2;\xib_2)$; methods do not have control over
the selection of the random variables $\xib_1$ and $\xib_2$.
Many one-point bandit methods do not fall neatly into the frameworks discussed
in this survey \cite{Agarwal2011,Belloni15,Bubeck2017}. See
\tabref{one-pt_rates} for a summary of best known WCC results
of \anlchange{type} \eqref{convex_stoch_stationarity}
for one-point
bandit feedback methods.

One example of a method using one-point bandit feedback, whose development falls
naturally into our discussion thus far, is given by
\citeasnoun{Flaxman2005};
they analyse a method resembling
\algref{rs-nesterov}, but the gradient-free oracle is chosen as
\begin{equation} \label{eq:gradient_free_oracle_stoch_onept}
 \gb_{\mu}(\xb;\ub;\xib) \defined \ffrac{\tilde{f}(\xb + \mu\ub;\xib)}{\mu}\ub,
\end{equation}
where $\ub$ is drawn uniformly from the unit $n$-dimensional sphere. 
That is, given a realization $\xib\in\Xib$, a stochastic gradient
estimator based on a stochastic evaluation at the point $\xb + \mu\ub$
is computed via \eqref{eq:gradient_free_oracle_stoch_onept}. 
For this method, and for general convex $f$,
\citeasnoun{Flaxman2005}
demonstrate a WCC bound 
of type \eqref{convex_stoch_stationarity}
in $\bigo{n^2\epsilon^{-4}}$
for general convex $f$.
For strongly convex $f$, \citeasnoun{Flaxman2005}
demonstrate a WCC bound 
of type \eqref{convex_stoch_stationarity}
in $\bigo{n^2\epsilon^{-3}}$.
Various extensions of these results are given in \citeasnoun{Saha2011}, \citeasnoun{Dekel2015} and 
\citeasnoun{Gasnikov2017}.
To the best of our knowledge, the best known upper bound on
the WCC for smooth, strongly
convex problems is in $\bigo{n^2\epsilon^{-2}}$ and the best known upper bound on
the WCC for smooth, convex problems
is in $\bigo{n\epsilon^{-3}}$ \cite{Gasnikov2017}.

For the solution
of \eqref{eq:stoch_prob} when $\Omegab=\Reals^n$,
\citeasnoun{pmlr-v49-bach16} analyse a method resembling
\algref{rs-nesterov} wherein the gradient-free oracle is chosen as
\begin{equation}
 \label{eq:gradient_free_oracle_stoch_bach}
 \gb_{\mu}(\xb;\ub;\xib_{+},\xib_{-}) \defined \ffrac{\tilde{f}(\xb + \mu\ub;\xib_{+})-\tilde{f}(\xb - \mu\ub;\xib_{-})}{\mu}\ub,
\end{equation}
where $\ub$ is again drawn uniformly from the unit $n$-dimensional sphere. 
We remark that $\xib_{+},\xib_{-}\in\Xib$ in \eqref{eq:gradient_free_oracle_stoch_bach}
are \emph{different} realizations; 
this distinction will become particularly relevant in \secref{two-point-feedback}.
For the solution of \eqref{eq:stoch_prob} when $\Omegab\subsetneq\Reals^n$,
\citeasnoun{pmlr-v49-bach16} also consider a gradient-free oracle like \eqref{eq:gradient_free_oracle_stoch_onept}. 
\citeasnoun{pmlr-v49-bach16} are particularly interested in how the smoothness of $f$ can be exploited
to improve complexity bounds in the one-point bandit feedback paradigm; 
they define a parameter $\beta$ and say that a function $f$ is $\beta$-smooth provided
$f$ is almost everywhere $(\beta-1)$-times Lipschitz-continuously differentiable
(a strictly weaker condition than assuming $f\in\cLC^{\beta-1}$).
When $\beta=2$, \citeasnoun{pmlr-v49-bach16} recover similar results to those seen in 
\tabref{one-pt_rates} for when $f\in\cLC^1$. 
When $\beta>2$, however,
in both the constrained and unconstrained case, \citeasnoun{pmlr-v49-bach16} prove a WCC result
of type \eqref{convex_stoch_stationarity} in $\bigo{n^2\epsilon^{2\beta/(1-\beta)}}$ 
when $f$ is convex and $\beta$-smooth.
Further assuming that $f$ is $c$-strongly convex, 
\citeasnoun{pmlr-v49-bach16} prove a WCC result
of type \eqref{convex_stoch_stationarity} 
in $\bigo{n^2\epsilon^{(\beta+1)/(1-\beta)}/c}$ .
Notice that asymptotically, as $\beta\to\infty$, this bound is in $\bigo{n^2\epsilon^{-1}/c}$.
This asymptotic result is particularly important because it attains the lower bound on 
optimization error demonstrated by \citeasnoun{Shamir13} for strongly convex $\infty$-smooth (quadratic) $f$.

\begin{table}
\renewcommand{\tabcolsep}{5.7pt} 
 \caption{Best known WCCs for $N_\epsilon$, the number of evaluations
 required to bound \eqref{convex_stoch_stationarity}, for one-point
 bandit feedback. $N_\epsilon$ is given only in terms of $n$ and $\epsilon$.
 See text for the definition of $\beta$-smooth; here $\beta>2$. Method types 
include random search (RS), mirror descent (MD) and ellipsoidal. 
 \label{table:one-pt_rates}}
{\small 
\begin{tabular*}{\textwidth}{lll}
\hline\hline
Assumption on $f$ & $N_\epsilon$ & Method type (citation) \\ \hline
convex, $f\in\cLC^0$ & $n^2 \epsilon^{-4}$ & RS \cite{Flaxman2005}\\[2pt]
	 & $n^{\gfrac{13}{2}} \epsilon^{-2}$ & ellipsoidal \cite{Belloni15}\\ \hline
$c$-strongly convex, $f\in\cLC^0$ & $n^2 \epsilon^{-3}$ & RS 
\cite{Flaxman2005}\\ \hline
convex, $f\in\cLC^1$ & $n \epsilon^{-3}$ & MD \cite{Gasnikov2017}\\[2pt]
 	 & $n^{\gfrac{13}{2}} \epsilon^{-2}$ & ellipsoidal \cite{Belloni15}\\ \hline
$c$-strongly convex, $f\in\cLC^1$ & $n^2 \epsilon^{-2}$ & MD 
\cite{Gasnikov2017}\\ \hline
convex, $\beta$-smooth & $n^2 \epsilon^{-\gfracp{2\beta}{\beta-1}}$ & RS 
\cite{pmlr-v49-bach16}\\ \hline
$c$-strongly convex, $\beta$-smooth & 
$n^2 \epsilon^{-\gpfracp{\beta+1}{\beta-1}}$ & RS \cite{pmlr-v49-bach16}\\ \hline
 \end{tabular*}}
\end{table}
We also note that \citeasnoun{Bubeck2017} conjecture that a particular kernel method can
achieve a WCC of type \eqref{convex_stoch_stationarity} in $\bigo{n^3\epsilon^{-2}}$
for general convex functions.
In light of well-known results in deterministic convex optimization, the WCCs summarized in
\tabref{one-pt_rates}
may be surprising. 
In particular,
for any $c$-strongly convex function $f\in\cLC^1$, the best known WCC results
are in $\bigo{\epsilon^{-2}}$.
We place particular emphasis 
on this result 
because it illustrates a gap between derivative-free and
derivative-based optimization that is not just a factor of $n$. In this
particular one-point bandit feedback setting, there do not seem to exist methods that achieve the
optimal\footnote{Optimal here is meant in a minimax information-theoretic sense~\cite{AWBR2009}.}
$\bigo{\epsilon^{-1}}$ convergence rate
attainable by gradient-based methods for smooth strongly convex stochastic optimization.
\citeasnoun{Hu2016} partially address this issue concerning one-point bandit feedback, which they refer to as
`uncontrolled noise'. These observations motivated the study of two-point (multi-point) bandit feedback,
which we will discuss in the next section, \secref{two-point-feedback}.

We further remark that every WCC in \tabref{one-pt_rates} has a polynomial dependence on
the dimension $n$, raising natural questions about the applicability of these methods in high-dimensional
settings. 
\citeasnoun{Wang2018AISTATS} consider a mirror descent method employing a special gradient-free oracle computed 
via a compressed sensing technique.
They prove a WCC of type \eqref{convex_stoch_stationarity} in 
$\bigo{\log(d)^{3/2}n_z^2\epsilon^{-3}}$,
under additional assumptions on derivative
sparsity, most importantly, that for all $\xb\in\Omegab$, $\|\nabla f(\xb)\|_0\leq n_z$ for some 
$n_z$. Thus, provided $n_z\ll n$, the polynomial dependence on $n$ becomes a logarithmic dependence on $n$,
at the expense of a WCC with a strictly worse dependence on $\epsilon$ than $\epsilon^{-2}$. 
In \secref{sepspa}, we discuss 
methods that similarly exploit known sparsity of objective function derivatives. 

\anlchange{In an extension of \algref{rs-nesterov}, \citeasnoun[Chapter~3]{RChen2015}}
dynamically updates the difference parameter by exploiting knowledge of the changing variance of
$\xib$.

\subsubsection{Two-point (multi-point) bandit feedback}
\label{sec:two-point-feedback}

We now focus on the stochastic paradigm of two-point (or multi-point) bandit
feedback. In this setting of bandit feedback,
we do not encounter the same gaps in WCC results
between derivative-free and derivative-based optimization
that exist for one-point bandit feedback. 

The majority of methods analysed in the two-point bandit feedback setting
are essentially random-search methods of the form \algref{rs-nesterov}.
The gradient-free oracles from \eqref{gradient_free_oracle} in the two-point setting
takes one of the two forms
\begin{equation} \label{eq:gradient_free_oracle_stoch}
 \begin{aligned}
 \gb_{\mu_k}(\xb;\ub;\xib) & \defined 
\df{\tilde{f}(\cdot,\xib);\xb;\ub;\mu_k}
\Bb\ub , \quad \mbox{or}\\
 \hat{\gb}_{\mu_k}(\xb;\ub;\xib) & \defined 
\dc{\tilde{f}(\cdot,\xib);\xb;\ub;\mu_k}
\Bb\ub.
 \end{aligned}
\end{equation}
The key observation in \eqref{eq:gradient_free_oracle_stoch} is that $\xib$
denotes a single realization used in the computation of both function
values in the definitions of the oracles (see \eqref{eq:fdiff} and \eqref{eq:cdiff}).
 This assumption of
`controllable realizations' separates two-point bandit feedback from the
more pessimistic one-point bandit feedback.
This property of being able to recall a single realization $\xib$ for two (or more)
evaluations of $\tilde{f}$ is precisely why this setting of bandit feedback 
is called `two-point' (or `multi-point'). 
This property will also be exploited in \secref{stoch_recursions}. 

The early work of
\citeasnoun{Agarwal2010} directly addresses the discussed gap in WCC results
and demonstrates that a random-search method resembling
\algref{rs-nesterov}, but applied in the two-point (or multi-point) setting as
opposed to the one-point setting, attains a WCC of type \eqref{convex_stoch_stationarity}
in $\bigo{\epsilon^{-1}}$;
this is a complexity result matching the optimal rate (in terms of
$\epsilon$-dependence) shown by
\citeasnoun{AWBR2009}. 
\anlchange{See
\tabref{two-pt_rates} for a summary of best known WCC results
of type \eqref{convex_stoch_stationarity}
for two-point
bandit feedback methods.}

\citeasnoun{Nesterov2015} provide a WCC result for \algref{rs-nesterov} 
using the stochastic gradient-free oracles
\eqref{eq:gradient_free_oracle_stoch}, but strictly better WCCs have since been established. 
We also note that, in contrast with the
gradient-free oracles of \eqref{gradient_free_oracle} in \secref{rand_det_rs},
the difference parameter $\mu$ is written as $\mu_k$ in \eqref{eq:gradient_free_oracle_stoch},
indicating that a
selection for $\mu_k$ must be made in the $k$th iteration.
In the works that we discuss
here, $\mu_k$ in \eqref{eq:gradient_free_oracle_stoch} is either chosen as a
constant sufficiently small 
or
else $\mu_k\to 0$ at a rate typically of the order of $1/k$. We also remark
that many of the results discussed in this section trivially hold for deterministic
convex problems, and can be seen as an extension of results concerning
methods of the form of \algref{rs-nesterov}. 

Provided
$f\in\cLC^0$,
the
best known WCCs of type \eqref{convex_stoch_stationarity}
(with variations in constants) can be found in
\citeasnoun{Duchi2015}, \citeasnoun{Gasnikov2017} and \citeasnoun{Shamir2017}. These works
consider variants of mirror-descent methods with approximate
gradients given by estimators of the form
\eqref{eq:gradient_free_oracle_stoch}; see \algref{2ptmirror} for a description
of a basic mirror-descent method. 
\algref{2ptmirror} depends on the concept of a \emph{Bregman divergence}
$D_{\psi}$,
used in \lineref{prox_prob} of \algref{2ptmirror} to define a proximal-point
subproblem. 
The Bregman divergence is defined by a function $\psi:\Omegab\to\Reals$, 
which is assumed to be $1$-strongly convex with respect to the norm $\|\cdot\|_p$.
\begin{algorithm}[tb]
 Choose initial point $\xb_0$, sequence of step sizes $\{\alpha_k\}$, sequence of difference
 parameters $\{\mu_k\}$, and distribution of $\xib$

 \For{$k=1,2,\ldots,T-1$}
 {
 Sample $\ub_k$ uniformly from the unit sphere $\B{\zerob}{1}$\\
 Sample a realization $\xib_k$ \\
 $\gb_k\gets\gb_{\mu_k}(\xb_k;\ub_k;\xib_k)$ using an oracle from \eqref{eq:gradient_free_oracle_stoch}\\
 $\xb_{k+1} \gets \argmin_{\yb\in\Omegab} \gb_k^{\T}\yb + \tfrac{1}{\alpha_k}D_{\psi}(\yb,\xb_k)$ \label{line:prox_prob} \\
 }
 \caption{Mirror-descent method with two-point gradient estimate 
\label{alg:2ptmirror}}
\end{algorithm}
To summarize many findings in this area
\cite{Duchi2015,Gasnikov2017,Shamir2017}, if $\|\cdot\|_q$ is the dual norm to $\|\cdot\|_p$ (\ie\
$p^{-1}+q^{-1}=1$) where $p\in\{1,2\}$
and $R_p$ denotes the radius of the feasible set in the $\|\cdot\|_p$-norm, 
then a bound on WCC of type \eqref{convex_stoch_stationarity} 
in $\bigo{n^{2/q}R_p^2\epsilon^{-2}}$
can be established for
a method like \algref{2ptmirror} in the two-point feedback setting where $f\in\cLC^0$.
These WCCs for methods like \algref{2ptmirror} 
are responsible for the popularity of mirror-descent methods 
in machine learning.
For many machine learning problems, solutions are typically sparse, and so,
in some sense, 
$R_1\leq R_2$. 
Thus, using a function $\psi$ that is $1$-strongly convex with respect to the $\|\cdot\|_1$-norm 
(\eg\ simply letting $\psi=\|\cdot\|_1$) may be preferable
to $p=2$,
in both theory and practice. 

\citeasnoun{Duchi2015} also provide an information-theoretic lower bound on convergence rates
for any method in the (non-strongly) convex, $f\in\cLC^0$, two-point feedback setting. This bound matches 
the best known WCCs up to constants,
demonstrating that these results are tight.
This lower bound is still of the order of $\epsilon^{-2}$, 
matching the result of \citeasnoun{AWBR2009} in $\epsilon$-dependence in the case where $f\in\cLC^0$.
It is also remarkable that this result is only a factor of $\sqrt{n}$ worse
than the bounds provided by \citeasnoun{AWBR2009} for the derivative-based case,
as opposed to the factor of $n$ that one may expect.

Additionally assuming $f(\xb)$ 
is strongly convex (but still assuming $f\in\cLC^0$),
\citeasnoun{Agarwal2010} prove, for a method like \algref{rs-nesterov}, a WCC
in $\bigo{n^2\epsilon^{-1}}$.
Using a method like \algref{2ptmirror},
\citeasnoun{Gasnikov2017} improve the dependence on $n$ in this WCC
 to
 $\bigo{n^{2/q}c_p^{-1}\epsilon^{-1}}$, provided $f$ is
$c_p$-strongly convex with respect to $\|\cdot\|_p$.

We now address the case where 
$f$ is convex and $f\in\cLC^1$.
Given an assumption that $\|\gb_{\mu}(\xb;\ub;\xib)\|$ is
uniformly bounded,
 \citeasnoun{Agarwal2010}
demonstrate a WCC of type \eqref{convex_stoch_stationarity}
in $\bigo{n^2\epsilon^{-1}}$. Dropping this somewhat
restrictive assumption
on the gradient-free oracle
and assuming instead that the oracle used
in a method like \algref{rs-nesterov} has bounded variance (\ie\ the oracle
satisfies \anlchange{$\Ea[\xib]{\|\gb_\mu(\xb;\ub;\xib)-\nabla f(\xb)\|^2}\leq
\sigma^2$}), \citeasnoun{Ghadimi2013} 
prove a WCC of a type similar to \eqref{convex_stoch_stationarity} (we avoid a discussion
of randomized stopping) in
$\bigo{\max \{n\Lip{g}\|\xb_0-\xb_*\|\epsilon^{-1},n\sigma^2L_g\|\xb_0-\xb_*\|\epsilon^{-2} \}}$. 
We mention that \citeasnoun{Nesterov2015} hinted at a similar WCC result, but with 
a strictly worse dependence on $n$, and different assumptions on $\xib$. 

\begin{table}
 \caption{Best known WCCs of type \eqref{convex_stoch_stationarity}
 for two-point bandit feedback. $N_\epsilon$ is given only in terms of $n$ and $\epsilon$.
 See text for the definition of $p,q$. $R_p$ denotes the size of the feasible set in the $\|\cdot\|_p$-norm. 
 If $f$ is $c_p$-strongly convex, then $f$ is strongly 
convex with respect to the $\|\cdot\|_p$-norm with constant $c_p$.
 $\sigma$ is the standard deviation on the gradient estimator $\gb_{\mu}(\xb;\ub;\xib)$ 
(\ie\ $\Ea[\xib]{\|\gb_{\mu}(\xb;\ub;\xib)-\nabla 
f(\xb)\|^2}\leq\sigma^2$).
The Lipschitz constant of the gradient $\Lip{g}$ is defined by the
$\|\cdot\|_2$-norm. 
 $\star$ denotes the additional assumption that
 $\Ea[\xib]{\|\gb_{\mu}(\xb;\ub;\xib)\|}<\infty$. 
 Method types 
include random search (RS), mirror descent (MD) and accelerated mirror descent 
(AMD).
 \label{table:two-pt_rates}}
\renewcommand{\tabcolsep}{3pt} 
{\small 
\begin{tabular*}{\linewidth}{lll}
\hline\hline
Assumption on $f$ & $N_\epsilon$ & Method type (citation) \\ \hline 
convex & $n^{\gfrac{2}{q}}R_p \epsilon^{-2}$ & MD \cite{Duchi2015,Shamir2017} \\ 
       &                                     &\hfill    \cite{Gasnikov2017} \\ \hline
${c_p}$-strongly convex & $n^{\gfrac{2}{q}} c_p^{-1} \epsilon^{-1}$ & MD \cite{Gasnikov2017}\\ \hline
convex, smooth & $\max\biggl\{\ffrac{n\Lip{g}R_2}{\epsilon},\ffrac{n\sigma^2}{\epsilon^2}\biggr\}$ & RS \cite{Ghadimi2013}\\[9pt]
& $\max\biggl\{\ffrac{n^{\gfrac{2}{q}}\Lip{g}R_p^2}{\epsilon},\ffrac{n^{\gfrac{2}{q}}
\sigma^2R_p^2}{\epsilon^2}\biggr\}$ & MD \cite{Dvurechensky2018}\\[9pt]
& $\max\biggl\{n^{\gfrac{1}{2}+\gfrac{1}{q}}\sqrt{\ffrac{\Lip{g}R_p^2}{\epsilon}},
\ffrac{n^{\gfrac{2}{q}}\sigma^2R_p^2}{\epsilon^2}\biggr\}$ & AMD \cite{Dvurechensky2018}\\ \hline
convex, smooth, $\star$ & $n^{\gfrac{2}{q}}R_p \epsilon^{-2}$ & MD \cite{Duchi2015}\\ \hline
$c_p$-strongly conv., smooth, $\star$ & $n^{\gfrac{2}{q}} c_p^{-1} \epsilon^{-1}$ & MD \cite{Gasnikov2017}\\ \hline
 \end{tabular*}}
 \end{table}

\section{Methods for structured objectives} %
\label{sec:structured}
The methods discussed in Sections~\ref{sec:det_det} and \ref{sec:rand_det}
assume relatively little about the structure of the objective function $f$ beyond some differentiability
required for analysis. 
\secref{convex} considered the case where $f$ is convex, which resulted, 
for example, in improved worst-case complexity results. 
In this section, we consider a variety of assumptions about additional
known structure in $f$ (including non-linear least squares, sparse, 
composite and minimax-based functional forms) 
and methods designed to exploit this additional structure.
Although problems in this section could be solved by the general-purpose methods 
discussed in Sections~\ref{sec:det_det} and \ref{sec:rand_det}, 
practical gains should be expected by exploiting the additional structure.

\subsection{Non-linear least squares} 
\label{sec:nls}
A frequently encountered objective in many applications 
of computational science, engineering and industry
is
\begin{equation}
 \label{eq:nls}
 f(\xb) = \ffrac{1}{2}\|\Fb(\xb)\|_2^2 = \ffrac{1}{2}\sum_{i=1}^p F_i(\xb)^2.
\end{equation}
For example, data-fitting problems are commonly cast as \eqref{eq:nls};
given data $y_i$ collected at design sites $\theta_i$, one may need to estimate
the parameters $\xb$ of a non-linear model or simulation output that best fit the data. 
In this scenario, $F_i$ is represented by 
$F_i(\xb) = w_i (S_i(\thetab;\xb)-y_i )$, which is a 
weighted residual between the simulation output $S_i$ and target data $y_i$.
In this way, objectives of the form \eqref{eq:nls} (and their correlated 
residual generalizations) encapsulate both the solution of non-linear equations 
and statistical estimation problems.

The methods of \citeasnoun{Zhangdfo10}, \citeasnoun{Zhang2012} and \citeasnoun{SWCHAP14} use
the techniques of \secref{det_det_model} to construct models of the individual
$F_i$ (thereby obtaining a model of the Jacobian $\nabla \Fb$) in
\eqref{eq:nls}. These models are then used to generate search directions in a
trust-region framework resembling the Levenberg--Marquardt
method \cite{Levenberg1944,Marquardt1963,More1978}. 
The analysis by \citeasnoun{Zhangdfo10} and \citeasnoun{Zhang2012} demonstrates -- under
certain assumptions, such as $f(\xb_*)=0$ at an optimal solution $\xb_*$
(\ie\ that the associated data-fitting problem has zero residual) -- that 
the resulting methods achieve the same local quadratic convergence rate does as the
Levenberg--Marquardt method. 
\codes{POUNDERS} is a trust-region-based method 
for minimizing objectives of the form
\eqref{eq:nls} that uses a full 
Newton approach for each residual $F_i$ \cite{SWCHAP14,tao-man_3_10}.
Another model-based method 
\anlchange{(implemented in \codes{DFO-LS} \cite{CFMR2018})}
for minimizing functions of the form
\eqref{eq:nls} -- more closely resembling a Gauss--Newton method -- is 
analysed by \citeasnoun{Cartis2017}. Their method is shown to 
converge to stationary points of \eqref{eq:nls} even when $f(\xb_*) > 0$,
at the expense of slightly weaker theoretical guarantees on the convergence rate. 

\citeasnoun{IFNLS} proposes a hybridization of a Gauss--Newton method with 
implicit filtering (\algref{imfil} from \secref{imfil}) that estimates the
Jacobian of $\Fb$ by building linear models of each component $F_i$
using central differences \eqref{eq:cdiff} with an 
algorithmically updated difference parameter. 
This hybrid method is shown to demonstrate superlinear
convergence for zero-residual (\ie\ $f(\xb_*)=0$) problems. 
 
Earlier methods also used the vector $F$ in order to more efficiently address objectives 
of the form \eqref{eq:nls}.
\citeasnoun{Spendley1969} develops a simplex-based algorithm that employs quadratic
approximations obtained by interpolating the vector $F$ on the current simplex. 
\citeasnoun{Peckham1970} proposes an iterative process that refines models of each
component of $F$ using between $n+1$ and $n+3+n/3$ points.
\citeasnoun{Ralston1978} also develop derivative-free Gauss--Newton methods and 
highlight their performance relative to methods that do not exploit the structure in \eqref{eq:nls}.
\citeasnoun{Brown1971} consider a
variant of the Levenberg--Marquardt method that approximates gradients
using appropriately selected difference parameters. 

\anlchange{\citeasnoun{Li2000}} analyse the convergence of a derivative-free line-search method 
when assuming that the square of the Jacobian (\ie\ $\nabla \Fb(\xb)^{\T} \nabla \Fb(\xb)$)
of \eqref{eq:nls} is \anlchange{positive-definite} everywhere. 
\citeasnoun{Grippo2007} \anlchange{and \citeasnoun{LaCruz2006}} augment
a non-monotone line-search method to incorporate information about $F$ in
the case where $p = n$ in \eqref{eq:nls}.
\anlchange{\citeasnoun{LiLi2011} develop a line-search method that exploits a monotonicity property assumed about $\Fb$.
\citen{LaCruz2014} and \citen{Morini2018} address \eqref{eq:nls} in the case where simple convex constraints are present.}

\subsection{Sparse objective derivatives}
\label{sec:sepspa}
In some applications, it is known that 
\begin{equation}
 \label{eq:sparse_hessian}
 \nabla^2 f(\xb)_{ij}=\nabla^2 f(\xb)_{ji}=0 \quad \gforall\ (i,j) \in S,
\end{equation}
for all $\xb\in\Omegab$, where the index set $S$ defines the sparsity pattern of the Hessian.
Similarly, one can consider \emph{partially separable}
objectives of the form
\begin{equation}
\label{eq:partially_separable}
 f(\xb) = \sum_{i=1}^p F_i(\xb) = \sum_{i=1}^p F_i ( \{\xb_j\}_{j\in S_i} ),
\end{equation}
where each $F_i$ depends only on some subset of
indices $S_i\subset\{1,2,\ldots,n\}$.
The extreme cases of \eqref{eq:partially_separable} are totally separable functions, where
$p=n$ and $S_i=\{i\}$ for $i\in\{1,\ldots,n\}$. 
In this special case, \eqref{eq:det_prob} reduces to the minimization 
of $n$ univariate functions. 

Given the knowledge encoded in
\eqref{eq:sparse_hessian} and \eqref{eq:partially_separable},
derivative-free optimization methods need not consider
interactions between certain components of $\xb$ because they are known to be exactly zero.
In the context of the model-based methods of \secref{det_det_model},
particularly when using quadratic models, 
using this knowledge amounts to dropping monomials in $\phib(\xb)$ in \eqref{eq:quad_monomials}
corresponding to the non-interacting $(i,j)$ pairs from \eqref{eq:sparse_hessian}.
Intuitively, such an action reduces the degrees of freedom in \eqref{eq:interpolation_system} when building models,
necessitating fewer function evaluations in the right-hand side of 
\eqref{eq:interpolation_system}.

\citeasnoun{Colson2005ops} propose a trust-region method for functions of the form
\eqref{eq:partially_separable} that builds and maintains separate
fully linear models for the individual $F_i$ in an effort to use fewer objective 
function evaluations. Similarly, \citeasnoun{Colson2001} propose a method for the case when
$\nabla^2 f$ has a band or block structure that exploits knowledge when building models of the objective; 
the work of \citeasnoun{Colson2002} extends this work to general sparse objective Hessians. 
\citeasnoun{Bagirov2006} develop an algorithm that exploits the fact that efficient discrete gradient
estimates can be obtained for $f$ having the form \eqref{eq:partially_separable}.

In the context of pattern-search methods (discussed in \secref{det_det_DDS}),
\citeasnoun{Price2006} exploit knowledge of $f$ having the form \eqref{eq:partially_separable}
to choose a particular stencil of search directions when forming a positive spanning set. 
In particular, the stencil is chosen as $\{\eb_1,\ldots,\eb_n,\eb_{n+1},\ldots,\eb_{n+p}\}$,
where $\{\eb_1,\ldots,\eb_n\}$ are the elementary basis vectors and
\[
\eb_{n+i} = \sum_{j\in S_i} -\eb_j,
\]
for $i\in\{1,\ldots,p\}$. 
\citeasnoun{Frimannslund2010} also develop a \framework{DDS} method for \eqref{eq:partially_separable}, with the search directions determined based on a smoothed quadratic formed from previously evaluated points.

\citeasnoun{Bandeira11} also assume that $\nabla^2 f(\xb)$ is sparse but do not
assume knowledge of the sparsity structure (\ie\ $S$ in \eqref{eq:sparse_hessian} is not known).
They develop a quadratic model-based trust-region method where
the models are selected by a minimum $1$-norm solution to an underdetermined interpolation system
\eqref{eq:interpolation_system}.
Under certain assumptions on $f$, if $n_z$ is the number of non-zeros in the (unknown) sparsity pattern for
$\nabla^2f$, \citeasnoun{Bandeira11} prove that $\Yb$ in \eqref{eq:interpolation_system} must contain
only $\bigo{(n_z+n)\log(n_z+n)\log(n)}$ (as opposed to $\bigo{n^2}$) randomly
generated points in order to ensure that the constructed interpolation models
are fully quadratic models of $f$ with high probability. 
This work motivated the analysis of randomized trust-region
methods discussed in \secref{rand_det_tr} because the random underdetermined
interpolation models of \citeasnoun{Bandeira11} satisfy the assumptions made in
\citeasnoun{Bandeira2014}. 

In \secref{one-point-feedback}, we noted the work of \citeasnoun{Wang2018AISTATS},
who used assumptions of gradient and Hessian sparsity 
(in particular, $\|\nabla f(\xb)\|_0\leq n_z$)
to improve the
reduce the dependence on $n$ in a WCC of type \eqref{convex_stoch_stationarity}
from polynomial to logarithmic. 
Note that, similar to \citeasnoun{Bandeira11}, 
this is an assumption on knowing a universal bound (\ie\ for all $\xb\in\Omegab$)
on the cardinality 
$\|\nabla f(\xb)\|_0$ rather than the actual non-zero components.
Under a similar sparsity assumption, 
\citeasnoun{Balasubramanian2018} consider the two-point bandit feedback setting 
discussed in \secref{two-point-feedback} 
and show that a 
truncated\footnote{That is, all but the $n_z$ \anlchange{largest values} in $\xb_{k+1}$ are \anlchange{set to 0 in line~5 of \algref{rs-nesterov}}.}
version of the method proposed in \citeasnoun{Ghadimi2013} has a WCC of type \eqref{convex_stoch_stationarity}
in $\bigo{n_z \log(n)^2 / \epsilon^2}$. 
Like the result of \citeasnoun{Wang2018AISTATS}, 
this WCC result also exhibits a logarithmic dependence on $n$, provided $n_z\ll n$. 
As we will discuss in \secref{stoch_bandit}, \citeasnoun{Ghadimi2013} analyse a method
resembling \algref{rs-nesterov} to be applied to non-convex $f$ in \eqref{eq:stoch_prob}. 
\citeasnoun{Balasubramanian2018} prove that the unaltered method of \citeasnoun{Ghadimi2013}
applied to problems in this sparse setting achieves a WCC of type 
 \begin{equation}
 \label{eq:expected_first_order}
 \E{\|\nabla f(\xb_k)\|}\leq\epsilon,
 \end{equation}
in $\bigo{n_z^2 \log(n)^2 \epsilon^{-4}}$;
this WCC result once again eliminates a polynomial dependence on $n$ that would otherwise
exist in a non-sparse setting. 
This WCC result, however, maintains the same $\epsilon$-dependence as the method of 
\citeasnoun{Ghadimi2013} in the non-sparse setting; in this sense, the method of \citeasnoun{Ghadimi2013}
is `automatically' tuned for the sparse setting.

\subsection{Composite non-smooth optimization} \label{sec:nso}
Sometimes, the objective $f$ in \eqref{eq:det_prob}
is known to be non-smooth. Often, one has knowledge about the form of
non-smoothness present in the objective, and we discuss methods that exploit specific forms of non-smoothness 
in this section. For methods that do not access any structural
information when optimizing non-smooth functions $f$, see 
Sections~\ref{sec:det_det_DDS} and \ref{sec:det_det_others}.

We define composite non-smooth functions as those of the form
\begin{equation}
 \label{eq:composite}
 f(\xb) = h(\Fb(\xb)),
\end{equation}
where $h:\Reals^p\to\Reals$ is a non-smooth function (in contrast to smooth $h$ such as the sum of squares in \eqref{eq:nls}), 
and $\Fb:\Reals^n\to\Reals^p$ is continuously differentiable. 
In some of the works we cite, the definition of a composite non-smooth objective may 
include an additional smooth function $g$ so that the objective function has the form 
$f(\xb)+g(\xb)$, but we omit a discussion of this for the sake of focusing on the 
non-smooth aspect in \eqref{eq:composite}.

\subsubsection{Convex $h$}
When $h$ in \eqref{eq:composite} is 
convex (note that $f$ may still be non-convex due to non-convexity in $\Fb$),
one thrust of research extends the techniques in derivative-based composite non-smooth
optimization; see the works of \citeasnoun{Yuan1985} and \citeasnoun[Chapter~14]{Fletcher}.
For example, \citeasnoun{Yuan1985} use derivatives to construct convex 
first-order approximations of $f$ near $\xb$,
\begin{equation}
\label{eq:composite_subprob}
\ell(\xb+\sba) \defined h(\Fb(\xb) + \nabla \Fb(\xb)\sba),
\end{equation}
where $\nabla \Fb$ denotes the Jacobian of $\Fb$; see
\citeasnoun{Hare2017} for properties of such approximations in the derivative-free setting.
By replacing $\nabla \Fb$ in \eqref{eq:composite_subprob} with the matrix
$\Mb(\xb_k)$ containing the gradients of a fully linear approximation to $\Fb$ at
$\xb_k$, \citeasnoun{Grapiglia2016} and \citeasnoun{Garmanjani2016} independently analyse a
model-based trust-region method similar to \algref{tr} from \secref{det_det_TRM} that uses the non-smooth
trust-region subproblem
\begin{equation}\label{eq:trsp_psi}
 \minimize_{\sba: \| \sba \| \le \Delta_k} \ell(\xb_k+\sba) \defined h(\Fb(\xb_k) + \Mb(\xb_k)\sba).
\end{equation}
Note that only $\Fb$ is assumed to be a black-box; these methods exploit the fact that
$h$ is convex with a known form in order to appropriately solve \eqref{eq:composite_subprob}. 
Both \citeasnoun{Grapiglia2016} and \citeasnoun{Garmanjani2016} use the stationarity
measure of \citeasnoun{Yuan1985} in their analysis,
\begin{equation}
 \label{eq:psi_measure}
 \Psi(\xb) \defined \ell(\xb) - \min_{\|\sba\|\leq 1} \ell(\xb+\sba),
\end{equation}
for which it is known that $\Psi(\xb_*)=0$ if and only if $\xb_*$ is a 
critical point of $f$ in the sense that $\ell(\xb_*) \leq \ell(\xb_*+\sba)$
for all $\sba\in\Reals^n$. Worst-case complexity results that bound the effort required to attain $\| \Psi(\xb_k) \| \le \epsilon$ 
are included in \tabref{rates}.
 
Methods using the local approximation \eqref{eq:composite_subprob} require convexity of $h$ 
in order for $\Psi$ in \eqref{eq:psi_measure} to be interpreted
as a stationarity measure. From a practical perspective, the form of $h$
directly affects the difficulty of solving the trust-region subproblem. 
\citeasnoun{Grapiglia2016} demonstrate this
approach on a collection of problems of the form
\begin{equation}
 \label{eq:finite_minimax}
 f(\xb) = \max_{i=1,\ldots,p} F_i(\xb),
\end{equation}
where each $F_i$ is assumed smooth. 
\citeasnoun{Garmanjani2016} test their method on a collection 
of 
problems of the form
\begin{equation}
 \label{eq:l1_norm}
 f(\xb) = \|\Fb(\xb)\|_1 = \sum_{i=1}^p |F_i(\xb)|.
\end{equation}
For objectives of the form \eqref{eq:finite_minimax} or \eqref{eq:l1_norm}, the
associated trust-region subproblems \eqref{eq:trsp_psi} can be cast
as linear programs when the $\infty$-norm defines the trust region.
(An early example of such an approach appears in \citeasnoun{Madsen1975}, where
linear approximations to each $F_i$ in \eqref{eq:finite_minimax} are constructed.)
Although more general convex $h$ could fit into this framework, one must be wary
of the difficulty of the resulting subproblems.

Direct-search methods have also been adapted
for composite non-smooth functions of specific forms. In these variants,
knowledge of $h$ in \eqref{eq:composite} informs the selection of search
directions in a manner similar to that described in \secref{sepspa}. 
\citeasnoun{Ma2009} and \citeasnoun{Bogani2009} consider the
cases of $f$ of the form \eqref{eq:finite_minimax} and \eqref{eq:l1_norm},
respectively.

Objectives of the form \eqref{eq:finite_minimax} are
also addressed by \citeasnoun{Hare2017b}, who develop an algorithm that decomposes 
such problems into orthogonal subspaces associated with directions of non-smoothness 
and directions of smoothness. The resulting derivative-free $VU$-algorithm 
employs model-based estimates of gradients to form and update this decomposition \cite{Hare2014}.
\citeasnoun{Liuzzi2006} address finite minimax problems by converting the
original problem into a smooth problem using an exponential penalty function.
Their \framework{DDS} method adjusts the penalty parameter via a rule that
depends on the current step size in order to guarantee convergence to a Clarke
stationary point.

Approximate gradient-sampling methods are developed and analysed by
\citeasnoun{Hare2013b} and \citeasnoun{Hare2013} for the finite mini\-max problem
\eqref{eq:finite_minimax}. These methods effectively exploit the
subdifferential structure of $h(\yb) = \max_{i=1,\ldots,p} y_i$ and employ 
derivative-free approximations of each $\nabla F_i(\xb)$. \citeasnoun{LMW16} propose
a variant of gradient sampling, called \emph{manifold sampling},
for objectives of the form \eqref{eq:l1_norm}.
Unlike (approximate) gradient sampling, manifold sampling does not depend on a random sample
of points to \anlchange{estimate the $\epsilon$-subdifferential}.

\subsubsection{Non-convex $h$}
\label{sec:cno_noncon}

When $h$ is non-convex, minimization of \eqref{eq:composite} is considerably more challenging than when $h$ is convex.
Few methods exist that exploit the structure of non-convex $h$.
One of the many challenges is that
 the model in \eqref{eq:composite_subprob} 
may no longer be an underestimator of $h$.
\citeasnoun{KLW18} propose a manifold sampling method for piecewise linear $h$; in
contrast to the previously discussed methods, this method does not require that $h$ be convex. 
Other methods applicable for non-convex $h$ employ \emph{smoothing functions}.

As mentioned in \secref{det_det_DDS}, the worst-case complexity of \framework{DDS} methods applied to
non-smooth (Lipschitz-continuous) objective functions is 
difficult to analyse. The reason that \framework{DDS} methods 
generate an asymptotically dense set of polling directions is to ensure that no
descent directions exist.
An exception to this generality, however, is functions for which an appropriate
smoothing function exists. 
Given a locally Lipschitz-continuous $f$, we say that 
$f_{\mu}:\Reals^n \rightarrow \Reals$
is a smoothing function for $f$ provided that for any $\mu\in(0,\infty)$,
$f_{\mu}$ is continuously differentiable 
and that 
\[
\lim_{\zb\to \xb,\mu\to 0^{+}} f_{\mu}(\zb) = f(\xb),
\]
for all $\xb\in\Reals^n$.

Thus, if a smoothing function $f_{\mu}$ exists for $f$, 
it is natural to iteratively apply a method for smooth unconstrained optimization
to obtain approximate solutions $\xb_k$ 
to $\min_{\xb} f_{\mu_k}(\xb)$ while decreasing $\mu_k$. We roughly prescribe such a smoothing
method in \algref{smoothing}.

\begin{algorithm}[tb]
 Set initial smoothing parameter $\mu_1>0$, terminal smoothing parameter $\mu_*<\mu_1$, and decrease parameter $\gamma\in(0,1)$\\
 Choose initial point $\xb_0$ and smooth optimization method $\mathfrak{M}$\\
 $k \gets 1$\\
 \While{$\mu_k < \mu_*$}{
 Apply $\mathfrak{M}$ to $f_{\mu_k}$, supplying
 $\xb_{k-1}$ as an initial point to $\mathfrak{M}$, until
 a termination criterion is satisfied and $\xb_k$ is returned\\
 $\mu_{k+1} \gets \gamma\mu_k$\\
 $k\gets k+1$\\
 }
 \caption{Smoothing method \label{alg:smoothing}}
\end{algorithm}

\citeasnoun{Garmanjani2012} consider the 
\framework{DDS} framework analysed by \citeasnoun{Vicente2013} as the method $\mathfrak{M}$ in
\algref{smoothing}. They terminate $\mathfrak{M}$ when the
step-size parameter $\alpha$ of \algref{ds} 
is sufficiently small, where the notion of sufficiently small scales with $\mu_k$
in \algref{smoothing}.
\citeasnoun{Garmanjani2012} prove a first-order stationarity result of the form
\begin{equation}\label{eq:nablaftilde}
 \liminf_{k\to\infty} \|\nabla f_{\mu_k}(\xb_k)\| = 0.
\end{equation}
Under certain assumptions
(for instance, that $h$ satisfies some regularity conditions at $\Fb(\xb_*)$) 
this first-order stationarity result is equivalent to $0\in \partial f(\xb_*)$;
that is, $\xb_*$ is Clarke stationary. 

\citeasnoun{Garmanjani2012} consider the decrease to be sufficient in 
\lineref{descent} of \algref{descent} if
$f(\pb_i) < f(\xb) - c_1\alpha^{3/2}$ for some $c_1>0$.
If, furthermore, $\mathfrak{M}$ terminates in each iteration of \algref{smoothing}
when $\alpha < c_2\mu_k^2$ for some $c_2>0$, 
then an upper bound on the number of function evaluations needed to obtain
\begin{equation}
\|\nabla f_{\mu_*}(\xb_k)\|\leq\epsilon,
\label{eq:1st_order_eps_stat}
\end{equation}
for $\epsilon\in(0,1)$ and $\mu_*\in\bigo{n^{-1/2}\epsilon}$,
is in 
$\bigo{\epsilon^{-3}}$; see \tabref{rates}. 
We note that while the sequence of smoothing parameters $\mu_k$ induces a type
of limiting behaviour of the gradients (as seen in \eqref{eq:1st_order_eps_stat}) 
returned by the method $\mathfrak{M}$ used in
\algref{smoothing}, this still does not necessarily recover 
elements of the Clarke subdifferential of $f$. The smoothing functions $f_{\mu_k}$
must satisfy an additional \emph{gradient consistency property} in order for
\algref{smoothing} to produce a sequence of \anlchange{points $\xb_k$} converging to Clarke
stationary points \cite[Theorem~9.67]{RockWets98}. 

\citeasnoun{Garmanjani2016} consider
the use of a model-based trust-region method $\mathfrak{M}$ in \algref{smoothing}. 
The authors demonstrate the first-order convergence result \eqref{eq:nablaftilde};
they also prove the same WCC as is proved by \citeasnoun{Garmanjani2012}.

\subsection{Bilevel and general minimax problems}\label{sec:bp}
Bilevel optimization addresses problems where a
\emph{lower-level} objective is embedded within an \emph{upper-level}
problem. Bilevel problems take the form
\begin{equation}\label{eq:bilevel}
 \begin{aligned}
 &\anlchange{\underset{\xb \in \Omegab}{\minimize}} && f^u(\xb,\xb^l) \\
 &\mbox{\,subject to} && \xb^l \in \argmin_{\zb \in \Omegab^l} \{ f^l(\xb,\zb)\},
 \end{aligned}
\end{equation}
where $f^u:\Omegab \subseteq \Reals^{n}\to\Reals$ and $f^l:\Omegab^l \subseteq \Reals^{n^l}\to\Reals$.
\citeasnoun{ConnVicente12} propose a model-based trust-region method for solving 
\eqref{eq:robust} 
in the absence of derivative information.
They show how to obtain approximations of the upper-level objective by solving
the lower-level problem to sufficient accuracy.
\citeasnoun{Mersha2011} and \citeasnoun{Zhang2014bl} develop \framework{DDS}-based algorithms for \eqref{eq:bilevel} under particular assumptions (\eg\ strict convexity of the lower-level problem).

A special case of \eqref{eq:bilevel} is when $f^l = -f^u$, which results in the
minimax problem \eqref{eq:det_prob}, where the objective is given by a
maximization:
\begin{equation}
 \label{eq:robust}
f(\xb) = \max_{\xb^l\in \Omegab^l} f^l(\xb,\xb^l).
\end{equation}
In contrast to the finite minimax problem \eqref{eq:finite_minimax}, the
objective in \eqref{eq:robust} involves a potentially infinite set $\Omegab^l$.

\citeasnoun{bertsimas2010robust} and \citeasnoun{BertsSimAnn10} consider \eqref{eq:robust}
when exact gradients of $f^l$ may not be available. The authors 
assume that approximate gradients of $f^l$ 
are available
and propose methods with convergence analysis restricted to functions $f$ in
\eqref{eq:robust} that are convex in $\xb$.
\citeasnoun{Ciccazzo13} and \citeasnoun{Latorre2019} develop derivative-free 
methods that employ approximate solutions of the inner 
problem in \eqref{eq:robust}. 
\citeasnoun{MMSMW2017} also consider \eqref{eq:robust} 
and develop a derivative-free method of
outer approximations for more general $f$. Their
analysis shows that the resulting limit points 
are Clarke stationary for $f$.

\vspace{5pt}
\section{Methods for stochastic optimization} %
\label{sec:stoch}
We now turn our attention to methods for solving the stochastic optimization problem
\eqref{eq:stoch_prob}. In \secref{stoch_convex_bandit}, we
considered the case where $f(\xb)=\Ea[\xib]{\tilde{f}(\xb;\xib)}$ is convex. In
this section, we lift the assumption of convexity to consider a more general
class of stochastic functions $\tilde{f}$. 

In general, the analysis of methods for stochastic optimization requires assumptions
on the random variable $\xib$. In this section, we use the convention
that $\xib\sim \Xib$ denotes that the random variable $\xib$ is from a
distribution $\Xib$ and that $\xib \in \Xib$ refers to a random variable in the
support of this distribution. 
Frequently, realizations $\xib\sim\Xib$ are assumed to be independent and
identically distributed (i.i.d.).
Throughout this section, we assume that $\Ea[\xib]{\tilde{f}(\xb;\xib)}$ exists
for each $\xb\in\Omegab$ and $f(\xb)=\Ea[\xib]{\tilde{f}(\xb;\xib)}$; that is,
the objective of \eqref{eq:stoch_prob} is \anlchange{well-defined}. Another common
assumption in the stochastic optimization literature 
is that some bound on the variance of
$\tilde{f}(\xb;\xib)$ is assumed, that is,
\begin{equation}
\label{eq:SAassumption}
\E[\xib]{(\tilde{f}(\xb;\xib)-f(\xb))^2} < \sigma^2 <\infty \quad 
\mbox{for all } \xb\in\Omegab.
\end{equation}

If, for a given $\xb$, $\nabla_{\xb} \tilde{f}(\xb;\xib)$ exists for each
$\xib\in\Xib$, then under certain regularity conditions 
it follows that 
$\nabla f(\xb) =
\Ea[\xib]{\nabla_{\xb} \tilde{f}(\xb;\xib)}$;
one such regularity condition is that 
\[\tilde{f}(\cdot ;\xib) \mbox{ is } 
L_{\tilde{f}(\cdot ;\xib)}\mbox{-Lipschitz-continuous and } \E[\xib]{L_{\tilde{f}(\cdot ;\xib)}}<\infty.\] 
We note that when first-order information is available, the assumption 
\eqref{eq:SAassumption} is often replaced by an assumption on the variance of 
the expected gradient norm; see \eg\ \citeasnoun[Assumption~4.3]{BottCurtNoce16}.
 In this setting, a key class of methods
for \eqref{eq:stoch_prob} are \emph{stochastic approximation} (\framework{SA})
methods; see the paper proposing \framework{SA} methods by
\citeasnoun{RobbinsMonro1951} and a survey of modern \framework{SA} 
methods
(often also referred to as `stochastic gradient' methods when first-order 
information is available)
by \citeasnoun{BottCurtNoce16}. 
Here we focus on situations where no objective derivative information is
available; that is, stochastic gradient methods are not directly applicable.
That said,
some of the work we discuss attempts to approximate
stochastic gradients, which are then used 
in an \framework{SA} framework.
As discussed in \secref{intro}, we will not address global optimization
methods, such as Bayesian optimization.%
\footnote{We recommend \citeasnoun{Shahriari2016} and \citeasnoun{Frazier2018} to readers
interested in recent surveys of Bayesian optimization.}

\secref{stoch_recursions} discusses stochastic approximation methods, and
\secref{stoch_ds} presents direct-search methods for stochastic optimization.
In \secref{stoch_model} we highlight
modifications to derivative-free model-based 
methods to address~\eqref{eq:stoch_prob}, and
in \secref{stoch_bandit} we discuss bandit methods for (non-convex) stochastic 
optimization. 

\subsection{Stochastic and sample-average approximation}
\label{sec:stoch_recursions}
One of the first analysed approaches for solving \eqref{eq:stoch_prob} is the method of
\citeasnoun{KieferWolfowitz}, inspired by the \framework{SA} method of
\citeasnoun{RobbinsMonro1951}. 
We state the basic Kiefer--Wolfowitz framework in 
\algref{kw}.
Since \citeasnoun{KieferWolfowitz} consider only univariate problems,
\algref{kw} is in fact the multivariate extension first of
\citeasnoun{Blum1954}. In \algref{kw}, $\nabla f(\xb_k)$ is approximated by
observing realizations of $\tilde{f}$ using central differences. That is, 
$\nabla f(\xb_k)$ is approximated by 
\begin{equation}
\label{eq:fdg}
\gK(\xb_k;\mu_k;\xib_k) = 
\begin{bmatrix}
 \ffrac{\tilde{f}(\xb_k + \mu_k \eb_1;\xib_1^+) - \tilde{f}(\xb_k - \mu_k \eb_1;\xib_1^-)}{2\mu_k}\\
 \vdots\\
 \ffrac{\tilde{f}(\xb_k + \mu_k \eb_n;\xib_n^+) - \tilde{f}(\xb_k - \mu_k \eb_n;\xib_n^-)}{2\mu_k}
\end{bmatrix}\!,
\end{equation}
where $\mu_k>0$ is a difference parameter, $\eb_i$ is the $i$th
elementary basis vector, and $2n$ realizations $\xib\sim\Xib$ are employed. 
The 
next point
$\xb_{k+1}$ is then set to be $\xb_k - \alpha_k \gK(\xb_k;\mu_k;\xib_k)$, where
$\alpha_k>0$ is a step-size parameter. 
As in \secref{rand_det}, we note that $\xb_{k+1}$ is a random variable that 
depends on the filtration generated by the method before $\xb_{k+1}$ is 
realized; this will be the case throughout this section.
In the \framework{SA} literature, the
sequences $\{ \alpha_k \}$ and $\{ \mu_k \}$ are often
referred to as \emph{gain} sequences. 

\begin{algorithm}[t]
 Choose initial point $\xb_0$, sequence of step sizes
 $\{\alpha_k\}$ and sequence of difference parameters
 $\{\mu_k\}$ 

 \For{$k=0,1,2,\ldots$}
 {
 Generate $\xib_k = (\xib_1^+,\xib_1^-,\dots,\xib_n^+,\xib_n^-)$\\
 Compute gradient estimate $\gK_k(\xb_k;\mu_k;\xib_k)$ via \eqref{eq:fdg}\\
 $\xb_{k+1}\gets \xb_k - \alpha_k \gK_k(\xb_k;\mu_k;\xib_k)$\\
 }
 \caption{Kiefer--Wolfowitz method \label{alg:kw}}
\end{algorithm}

Because evaluation of the function $f$ 
requires computing an expectation (and in contrast to the primarily monotone 
algorithms in \secref{det_det}), stochastic optimization methods generally do 
not monotonically decrease $f$. 
This is exemplified by \algref{kw}, which updates $\xb_{k+1}$ 
without considering the value of $\tilde{f}(\xb_{k+1}; \xib)$ 
for \emph{any} realization of $\xib$.

Historically, \algref{kw} has been analysed by a community more concerned
with stochastic processes than with optimization. Hence, convergence results
differ from those commonly found in the optimization literature. For
example, many results in the \framework{SA} literature consider a continuation
of the dynamics of \algref{kw} applied to the deterministic $f$ as an ordinary
differential equation (ODE) in terms of $\xb(t):\Reals\to\Reals^n$. That is, 
they
consider
\[
\ffrac{\mathrm{d}\xb}{\mathrm{d}t} = -\nabla f(\xb), \quad \xb=\xb(t).
\]
and define the set of fixed points of the ODE,
$ \Sb = \{\xb: \nabla f(\xb) = 0 \}.$
Many convergence results then demonstrate that the continuation $\xb(t)$
satisfies $\xb(t)\to \Sb$ with probability one as the continuous iteration 
counter $t\to \infty$; see \citeasnoun{Kushner2003} for a
complete treatment of such ODE results.

In order to prove that the sequence of \anlchange{points $\xb_k$} generated by \algref{kw} converges almost 
surely (\ie\ with probability one), 
conditions must be placed on the objective function, step sizes and
difference parameters.
In %
the \framework{SA} literature there is no single
consistent set of conditions, but there are nearly always conditions on the sequence of
step sizes $\{\alpha_k\}$ 
requiring $\alpha_k\to 0$ and
$\sum_k \alpha_k = \infty$. Intuitively, this divergence
condition ensures that any point in
the domain $\Omegab$ can be reached, independent of the history of iterations.
As one
example of convergence conditions,
\citeasnoun{Bhatnagar2013} prove almost sure convergence of \algref{kw} under
the following assumptions (simplified for presentation).
\begin{enumerate}\renewcommand{\theenumi}{(\arabic{enumi})}
\setlength\itemsep{4pt}
 \item The sequences of step sizes and difference parameters satisfy 
$\alpha_k > 0$, $\mu_k > 0$, 
 $\alpha_k \to 0$,
 $\mu_k \to 0$,
 $\sum_{k} \alpha_k = \infty$
 and $\sum_{k} \alpha_k^2\mu_k^{-2} < \infty$.
 \item The realizations $\xib\sim\Xib$ are i.i.d.~and the distribution 
$\Xib$ has a finite second moment. 
 \item The function $f$ is in $\cLC^1$.
 \item $\sup_{k}\{\|\xb_k\|\}<\infty$ with probability one. 
\end{enumerate}
Similar assumptions on 
algorithms of the form \algref{kw} appear throughout the \framework{SA} literature
\cite{Blum1954b,Derman1956,Sacks1958,Fabian1971,Kushner1979,Ruppert1991SA,SpallBook}. 
Convergence of \algref{kw} under similar assumptions to those above, but with
the modification that $\mu_k$ is fixed in every iteration to a sufficiently
small constant (that scales inversely with $\Lip{g}$),
is additionally demonstrated by \citeasnoun{Bhatnagar2013}.

In terms of WCCs, the convergence rates that have been historically derived for
\algref{kw} are also non-standard for optimization. In
particular, results concerning convergence rates are typically shown as
a convergence in distribution \cite[Chapter~3.2]{Durrett2010}: given a fixed $\xb_*\in\Sb$, 
\begin{equation}
\label{eq:kw_rate}
 \ffrac{1}{k^{\gamma}}(\xb_k-\xb_*)\to\cN(\zerob,\Bb),
\end{equation}
where $\gamma>0$ and $\Bb$ is a covariance matrix, the entries of which depend 
on algorithmic
parameters and $\nabla^2 f(\xb_*)$ (provided it exists). 
With few 
assumptions on $\xib$, it has been shown that
\eqref{eq:kw_rate} holds with $\gamma=1/3$ 
\cite{SPSA1992,LEcuyerYin1998}.
Observe that these convergence rates are distinct from WCC results like those
in \eqref{eq:nonconvex_stationarity}. 

Later, the use of \emph{common random numbers} (CRNs) was considered. In 
contrast to \eqref{eq:fdg}, which employs a
realization $\xib_k=(\xib_1^+,\xib_1^-,\dots,\xib_n^+,\xib_n^-)$, a gradient 
estimator in
the CRN regime uses a single realization $\xib_k$ and
has the form
\begin{equation}
\label{eq:fdg_crn}
\gK(\xb_k;\mu_k;\xib_k) = 
\begin{bmatrix}
 \dc{\tilde{f}(\cdot;\xib_k);\xb_k;\eb_1;\mu_k}\\
 \vdots\\
 \dc{\tilde{f}(\cdot;\xib_k);\xb_k;\eb_n;\mu_k}
\end{bmatrix}\!,
\end{equation}
where $\dc{\cdot}$ is defined in \eqref{eq:cdiff}

The difference between \eqref{eq:fdg} and \eqref{eq:fdg_crn} is analogous to
the difference between one-point and two-point bandit feedback in the context
of bandit problems (see \secref{stoch_convex_bandit}). In the CRN regime, we 
can recall a single realization
$\xib_k$ to compute a finite-difference approximation in each
coordinate direction. By using \eqref{eq:fdg_crn} as the gradient estimator in
\algref{kw}, the rate \eqref{eq:kw_rate} holds with $\gamma=1/2$
\cite{LEcuyerYin1998,Kleinman1999}.
Thus, as in the analysis of bandit methods, the use of 
CRNs allows for strictly better convergence rate results. 

\citename{Dai2016} \citeyear{Dai2016,Dai2016b} studies the complexity of \algref{kw}, 
as well as a method
that uses the estimator \eqref{eq:fdg} in \algref{2ptmirror}, under varying
assumptions on $\Xib$.
\anlchange{Dai} considers a gradient estimate of the form \eqref{eq:fdg_crn}
with $\df{\tilde{f}(\cdot;\xib_k);\xb_k;\eb_i;\mu_k}$
replacing each central difference; 
recall the definition of $\df{\cdot}$ in \eqref{eq:fdiff}.
\anlchange{Dai} %
 demonstrates
that the best rate of the form 
\eqref{eq:kw_rate}
achievable by \algref{kw} with forward differences has $\gamma=1/3$, even when common random numbers are used. 
However, a rate of the form \eqref{eq:kw_rate} with $\gamma=1/2$ can be achieved using forward
differences in \algref{2ptmirror};
\citename{Dai2016} \citeyear{Dai2016,Dai2016b} draws
parallels between this
result and the WCC of \citeasnoun{Duchi2015}, discussed in
\secref{two-point-feedback}.

We remark that the gradient estimate \eqref{eq:fdg} used in
\algref{kw} requires $2n$ evaluations of $\tilde{f}$ per iteration. Although 
replacing $\dc{\tilde{f}(\cdot;\xib_k);\xb_k;\eb_i;\mu_k}$ with
$\df{\tilde{f}(\cdot;\xib_k);\xb_k;\eb_i;\mu_k}$ could 
reduce this cost to $n+1$ evaluations of 
$\tilde{f}$ per iteration, it
is still desirable to reduce this per-iteration cost from $\bigo{n}$ to $\bigo{1}$
evaluations. The \codes{SPSA} method of \citeasnoun{SPSA1992} achieves this goal by
using 
the gradient estimator
\begin{equation}
\label{eq:spsa_est}
\gS(\xb_k;\mu_k;\xib_k;\ub_k) = \dc{\tilde{f}(\cdot,\xib_k);\xb_k;\ub_k;\mu_k}
\begin{bmatrix}
\ffrac{1}{[\ub_k]_1}\\
\vdots\\
\ffrac{1}{[\ub_k]_n}
\end{bmatrix}\!,
\end{equation}
where $\ub_k\in\Reals^n$ is randomly generated from some distribution in each
iteration.

The construction of \eqref{eq:spsa_est} requires evaluations of
$\tilde{f}(\cdot;\xib_k)$ at exactly two points. \algref{kw} is then modified by
replacing the gradient estimator $\gK_k(\xb_k;\mu_k;\xib_k)$ with 
$\gS_k(\xb_k;\mu_k;\xib_k;\ub_k)$.
Informally, the conditions on the distribution governing $\ub_k$ originally
proposed by \citeasnoun{SPSA1992} cause each entry of $\ub_k$ to be bounded away 
from 0 with high
probability (intuitively, to avoid taking huge steps). A simple
example distribution satisfying these properties is to let each entry of
$\ub_k$ independently follow a Bernoulli distribution with support $\{1,-1\}$,
both events occurring with probability $1/2$. Under appropriate assumptions
resembling those for \algref{kw}, the sequence $\{\xb_k\}$ generated by 
\codes{SPSA} can be shown to
converge in the same sense as \algref{kw}. Convergence rates of the form
\eqref{eq:kw_rate} matching those obtained for \algref{kw} have also been
established \cite{Gerencser1997,Kleinman1999}.

The performance of \framework{SA} methods is highly sensitive 
to the chosen sequence of step sizes $\{\alpha_k\}$ \cite{Hutchison2013}.
This mirrors the situation in gradient-based \framework{SA} methods where the 
tuning of algorithmic parameters is an active area of research 
\cite{Diaz2017,Ilievski2017,DeepHyper18}.

The \framework{SA} methods above consider only a single evaluation of
the stochastic function $\tilde{f}$ at any point. Other methods more accurately estimate
$f(\xb_k)$ by querying $\tilde{f}(\xb_k;\xib)$ for multiple, different realizations 
(`samples') of
$\xib$. These methods belong to the framework of sample average 
approximation, wherein 
the original problem \eqref{eq:stoch_prob} is replaced with a (sequence of) 
deterministic
\emph{sample-path problem(s)}:
\begin{equation}
 \minimize_{\xb \in \Omegab} \ffrac{1}{p} \sum_{i = 1}^{p}\tilde{f}(\xb; 
\xib_i).
\label{eq:SAA}
\end{equation}
Retrospective approximation methods \cite{Chen2001} vary the number of samples,
$p$, in a predetermined sequence $\{p_0, p_1, \ldots\}$; the accuracy to which 
each instance of
\eqref{eq:SAA} subproblem is solved can also vary as a sample average 
approximation 
method progresses. Naturally, the performance of such a method depends
critically on the sequence of sample sizes and accuracies used at each
iteration; \citeasnoun{Pasupathy2010} characterizes a class of sequences of
predetermined sample sizes and accuracies for which 
derivative-free retrospective approximation
methods can be shown to converge for smooth objectives.

Other approaches dynamically adjust the number of samples $p$ from iteration 
to the next. For example, the
method of \citeasnoun{Pasupathy2018} adjusts the number of samples $p_k$ to balance 
the contributions
from deterministic and stochastic errors in 
iteration $k$. 
The stochastic error at $\xb_k$ is then 
\[\biggl| f(\xb_k) - \ffrac{1}{p_k} \sum_{i = 1}^{p_k}\tilde{f}(\xb_k; 
\xib_{k,i})\biggr|.\]
The deterministic error 
is the difference between the objective $f$ and a specified approximation; for example, the deterministic error at $\xb_k-\alpha_k \nabla f (\xb_k)$ using a first-order Taylor approximation is
\[| f(\xb_k - \alpha_k \nabla f(\xb_k))- (
f(\xb_k) - \alpha_k\|\nabla f(\xb_k)\|^2 ) |.\]
\citeasnoun{Pasupathy2018} establish convergence rates for a variant 
of \algref{kw} drawing independent samples $\{\xib_{k,1}, \ldots, 
\xib_{k,p_k}\}$ in each iteration.

\subsection{Direct-search methods for stochastic optimization}
\label{sec:stoch_ds}

Unsurprisingly, researchers have modified methods for
deterministic objectives in order to produce methods appropriate for stochastic
optimization. For example, in the paper inspiring \citeasnoun{NelderMead}, 
\citeasnoun{Spendley1962} propose re-evaluating the point corresponding to the best simplex vertex if it hasn't
changed in $n+1$ iterations, saying that if the vertex is best `only by reason
of errors of observation, it is unlikely that the repeat observation will [be
the best observed point], and the point will be eliminated in due course'.
\citeasnoun{Barton1996} propose modifications to the Nelder--Mead method
in order to avoid premature termination
due to repeated shrinking. To alleviate this problem, they suggest
reducing the amount the simplex is shrunk, re-evaluating the best point after
each shrink operation, and re-evaluating each reflected point before performing
a contraction. 
\citeasnoun{Chang2012} proposes a Nelder--Mead variant that samples candidate points
and all other points in the simplex an increasing number of times; this method
ultimately ensures that stochasticity in the function evaluations will not
affect the correct ranking of simplex vertices. 

\citeasnoun{Sriver2009} augment a \framework{GPS} method with a ranking and
selection procedure and dynamically determine the number of samples 
performed for each polling point. The ranking and selection procedure allows 
the method to also address cases where $\xb$ contains discrete variables.
For the case of additive unbiased, Gaussian 
noise (\ie\ $\tilde{f}(\xb; \xib) = f(\xb) + \sigma \xib$ with $\xib$ from a 
standard normal distribution and $\sigma>0$ finite),
they prove that the resulting method converges 
almost surely to a stationary point of $f$.
For problems involving more general distributions, 
\citeasnoun{KimZhang2010} consider a \framework{DDS} method that employs 
the sample mean
\begin{equation}
\ffrac{1}{p_k} \sum_{i = 1}^{p_k}\tilde{f} (\xb; \xib_{k,i}),
 \label{eq:samplemean}
\end{equation}
with a dynamically increasing sample size $p_k$. 
They establish a consistency result and 
appeal to the convergence properties of \framework{DDS} methods.
\citeasnoun{Sankaran2010} propose a surrogate-assisted method for stochastic 
optimization inspired by stochastic collocation techniques (see \eg\ 
\citeb{Gunzburger2014}). Convergence for the method is established by appealing 
to the \framework{GPS} and \framework{MADS} mechanisms underlying the method.

\citeasnoun{Chen2016a} consider an \anlchange{implicit-filtering} method in which values of $f$
are observable only through the sample average \eqref{eq:samplemean}. 
\citeasnoun{Chen2016a} demonstrate that the sequence of points generated by the method converges
(\ie\ $\{\nabla f(\xb_k)\}$ admits a subsequence that converges to zero)
with probability one if the sample size $p_k$ increases to infinity.
Algorithmically, $p_k$ is adjusted to scale with the square of the inverse of the stencil step size 
($\Delta_k$ in \algref{imfil}).

\citeasnoun{Chen2018} consider the bound-constrained minimization of a composite
non-smooth function of the form \eqref{eq:composite}, where $h$ is
Lipschitz-continuous
(but non-smooth) and $\Fb$ is continuously differentiable. However,
they assume that values of $\Fb$ are observable only through sample averages 
 and that a smoothing function $h_{\mu}$ of $h$ (as discussed in
\secref{cno_noncon}) is available. 
They show that with probability one, the sequence of points from a smoothed
implicit-filtering method converges to a first-order stationary point, where the
stationarity measure is appropriate for non-smooth optimization.

\subsection{Model-based methods for stochastic optimization}
\label{sec:stoch_model}
Analysis of the model-based trust-region 
methods in \secref{det_det_model} generally depends on 
the construction of fully linear models of a
deterministic function $f$; see \eqref{eq:fullylinear}. 
In particular, methods of the form of \algref{tr} typically require that a model $m_k$ satisfy
\[
|f(\xb_k+\sba) - m_k(\xb_k+\sba)| \leq \kappaef \Delta_k^2 
\quad \gforall\ \sba \in \B{\zerob}{\Delta_k}.
\]
A natural model-based trust-region approach to stochastic optimization is
to build a model $m_k$ of the function $f$ 
by fitting the model to observed values of the stochastic function $\tilde{f}$. 
Intuitively, if such models satisfy \eqref{eq:fullylinear}, then 
an extension of the
analysis described in \secref{det_det_TRM} should also apply to the minimization of $f$ in
\eqref{eq:stoch_prob}.
The methods described here formalize the approximation properties of such models 
(which are stochastic because of their dependence on $\xib$) and employ the models in a trust-region framework.
For example,
by employing an estimator $\bar{f}_p$ of $f$ at each
interpolation point $\xb$ used in model construction, we can replace 
each function value $f(\xb)$ with $\bar{f}_p(\xb)$ in the interpolation system 
\eqref{eq:interpolation_system}. 
One example of such an estimator $\bar{f}_p$ is the sample average 
\eqref{eq:samplemean}.

Early work in applying derivative-free trust-region methods for stochastic optimization includes
that of \citeasnoun{Deng2006aua}, \anlchange{which modifies} 
the \codes{UOBYQA} method of
\citeasnoun{Powell2002a}. The $k$th iteration of the method of \citeasnoun{Deng2006aua} uses Bayesian
techniques to dynamically update a budget of $p_k$ new $\tilde{f}$ evaluations. This budget is then apportioned among the current set
of interpolation points $\yb\in\Yb$ in order to reduce the variance in each
value of $\bar{f}_{p_k}(\yb)$, with the authors using the sample mean for the estimator $\bar{f}_{p_k}$. \citeasnoun{Deng09} 
\anlchange{show} that, given assumptions on the sequence of evaluated $\xib$ (\ie\ the sample path), every limit point $\xb_*$ produced by this method is
stationary with probability 1. 

Another method in this vein, \codes{STRONG}, was proposed by
\citeasnoun{ChangHongWan2013} 
and combines response surface methodology
\cite{Box1987}
with a trust-region mechanism. In
the analysis of \codes{STRONG}, it is assumed that model gradients $\nabla
m_k(\xb_k)$ almost surely equal the true gradients $\nabla f(\xb_k)$ as
$k\to\infty$, which is algorithmically encouraged by monotonically increasing
the sample size $p_k$ in an inner loop. 
\codes{QNSTOP} by
\citeasnoun{CastlePhD} presents a similar approach using response surface models in
a trust-region framework, but its convergence analysis and assumptions mirror
those of stochastic approximation methods. 

Both \citeasnoun{Larson2016} and \citeasnoun{Chen2017} build on the idea of
probabilistically fully linear models in \eqref{eq:tr_prob}, which
essentially says that the condition \eqref{eq:fullylinear} needs to hold
on a given iteration only with some probability \cite{Bandeira2014}.
In contrast to the usage of such models in randomized methods for deterministic objectives
(the subject of \secref{rand_det_tr}), in stochastic optimization the filtration in 
\eqref{eq:tr_prob} also includes the realizations of the stochastic evaluations of $\tilde{f}$.
This probabilistic notion of uniform local model quality is powerful.
For example, although the connection is not made by \citeasnoun{Reiger2017}, 
this notion of model quality implies probabilistic descent properties such as those required by
\citeasnoun{Reiger2017}. 
This implication is an example of a setting in which stochastic gradient estimators can be
replaced by gradients of probabilistically fully linear models.

One way to satisfy \eqref{eq:tr_prob} is to build a 
regression model using randomly sampled points. For example,
\citeasnoun[Theorem~4.2.6]{MMenickelly2017}
shows that evaluating $\tilde{f}$ on a sufficiently 
large set of points uniformly sampled from $\B{\xb_k}{\Delta_k}$ can be used to 
construct a probabilistically fully linear regression model.

\citeasnoun{Larson2016} prove convergence of a probabilistic variant of
\algref{tr} in the sense that, for any $\epsilon>0$,
\[
 \lim_{k\to\infty} \P{\|\nabla f(\xb_k)\|>\epsilon} = 0.
\]
Under similar assumptions, \citeasnoun{Chen2017} prove 
almost sure convergence to a stationary point, that is, 
\begin{equation}
 \lim_{k\to\infty} \|\nabla f(\xb_k)\| = 0 \quad \mbox{with probability one}.
\label{eq:asc}
\end{equation}

\citeasnoun{BCMS2018} provide a WCC result for the variant of \algref{tr} presented
by \citeasnoun{Chen2017}.
\citeasnoun{BCMS2018} extend the analysis of \citeasnoun{Cartis2015} to study the
stopping time of
the stochastic process generated by the method of \citeasnoun{Chen2017}.
In contrast to previous WCC
results discussed in this survey, which bound the number of function
evaluations $N_{\epsilon}$ needed to attain some form of \emph{expected
$\epsilon$-optimality} (\eg\ \eqref{eq:gauss_nonconvex_complexity} 
or \eqref{eq:expected_first_order}),
\citeasnoun{BCMS2018} prove that 
the \emph{expected number of iterations}, $\E{T_\epsilon}$,
needed to achieve \eqref{eq:nonconvex_stationarity}
is in $\bigo{\epsilon^{-2}}$.
\citeasnoun{PaquetteScheinberg2018} apply similar analysis
to a derivative-free
stochastic line-search method, where they demonstrate that for
non-convex $f$, 
$\E{T_{\epsilon}}\in\bigo{\epsilon^{-2}}$, while for convex and strongly convex $f$,
$\E{T_{\epsilon}}\in\bigo{\epsilon^{-1}}$
and 
$\E{T_{\epsilon}}\in\bigo{\log(\epsilon^{-1})}$, respectively.
Since the number of function evaluations per iteration of the derivative-free
methods 
of \citeasnoun{BCMS2018} and \citeasnoun{PaquetteScheinberg2018}
is highly variable across iterations, the total work (in terms of 
function evaluations) is not readily apparent from such WCC results.

\citeasnoun{Larson2016} and \citeasnoun{Chen2017} 
demonstrate that
sampling $\tilde{f}$ on $\B{\xb_k}{\Delta_k}$
of the order of $\Delta_k^{-4}$ times will ensure that \eqref{eq:fullylinear} holds 
(\ie\ one can obtain a fully linear model) with high
probability. \citeasnoun{Shashaani2016} and \citeasnoun{Shashaani2018} take a related but
distinct approach.
As opposed to requiring that models be probabilistically fully linear,
their derivative-free trust-region method performs adaptive Monte Carlo sampling 
both at current points $\xb_k$ and interpolation points\anlchange{;}
 the number of samples $p_k$ is 
 chosen to balance a measure of statistical error with 
the optimality gap at $\xb_k$. 
\citeasnoun{Shashaani2018} prove that their
method achieves almost sure convergence of the form \eqref{eq:asc}.

A model-based trust-region method for constrained stochastic optimization, \codes{SNOWPAC},
is developed by \citeasnoun{Augustin2017}. 
Their method
addresses the stochasticity by employing Gaussian process-based models of robustness
measures such as expectation and conditional value at risk.
The approach used is an extension of the 
constrained deterministic method \codes{NOWPAC} of \citeasnoun{NOWPAC2014}, which we
discuss in \secref{constrained}.

\subsection{Bandit feedback methods}
\label{sec:stoch_bandit}

While much of the literature on bandit methods for stochastic optimization focuses on convex
objectives $f$ (as discussed in \secref{stoch_convex_bandit}), here we discuss treatment of non-convex objectives $f$.
We recall our notation and discussion from \secref{stoch_convex_bandit},
in particular the notion of regret minimization shown in~\eqref{expected_regret_convex}.

In the absence of convexity, regret bounds do not translate into bounds on optimization error
as easily as in \eqref{eq:regret_implies_opt}. 
Some works address the case where each $\tilde{f}(\cdot;\xib_k)$ in
\eqref{expected_regret_convex} is
Lipschitz-continuous
and employ a partitioning of the feasible region $\Omegab$ 
\cite{Kleinberg2008,Bubeck2011a,Bubeck2011b,Valko2013,Zhang2015}. 
These methods employ global optimization strategies that
we do not discuss further here. 

In another line of work, \citeasnoun{Ghadimi2013} consider the application of an algorithm like 
\algref{rs-nesterov} with the choice of gradient estimator
$\gb_\mu(\xb;\ub;\xib)$ from \eqref{eq:gradient_free_oracle_stoch}. 
Under an assumption of bounded variance of the estimator (\ie\
\anlchange{$\Ea[\xib]{\|\gb_\mu(\xb;\ub;\xib)-\nabla f(\xb)\|^2}\leq
\sigma^2$}), \citeasnoun{Ghadimi2013} prove a WCC result similar to the
one they obtained in the convex case; see
\secref{two-point-feedback}.
They show 
that an upper bound on the (randomized) number of iterations needed to attain
\begin{equation}
 \label{eq:ghadimi_lan_stationarity}
 \E{\|\nabla f(\xb_k)\|^2}\leq \epsilon
\end{equation}
is in $\bigo{\max\{n L_g \anlchange{\Rx} \epsilon^{-1},n L_g \anlchange{\Rx} \sigma^2\epsilon^{-2}\}}$.
Notice that the stationarity condition given in \eqref{eq:ghadimi_lan_stationarity}
involves a square on the gradient norm, making it distinct from
a result like
\eqref{eq:gauss_nonconvex_complexity} or \eqref{eq:expected_first_order}. 
Thus, assuming $\sigma^2$ is sufficiently large, the result of
\citeasnoun{Ghadimi2013} translates to a WCC of type \eqref{eq:expected_first_order}
in $\bigo{n^2\epsilon^{-4}}$.

\citename{Balasubramanian2018} \citeyear{Balasubramanian2018,Balasubramanian2019} propose a method that uses two-point bandit feedback (\ie\
a gradient estimator from \eqref{eq:gradient_free_oracle_stoch})
within a derivative-free conditional
gradient method \cite{Ghadimi2019}. The gradient estimator is used to define a linear
model, which is minimized over $\Omegab$ to produce a trial step. 
If $\Omegab$ is bounded, they show a WCC of type \eqref{eq:expected_first_order} that again 
grows like~$\epsilon^{-4}$.

By replacing gradients with estimators of the form \eqref{eq:gradient_free_oracle_stoch} in 
the stochastic variance-reduced gradient framework of machine learning 
\cite{Reddi2016a},
\citeasnoun{SLiu2018} prove a WCC of type \eqref{eq:ghadimi_lan_stationarity} in $\bigo{n\epsilon^{-1} + b^{-1}}$,
where $b$ is the size of a minibatch drawn with replacement in each iteration.
\citeasnoun{BGu2016} prove a similar WCC result in an asynchronous parallel computing environment 
for a distinct method using minibatches for variance reduction.

\vspace{5pt}
\section{Methods for constrained optimization} %
\label{sec:constrained}

In this section, we discuss derivative-free methods for problems where 
the feasible region $\Omegab$ is a proper subset of $\Reals^n$.
In the derivative-free setting, such constrained 
optimization problems can take many forms since an additional distinction is
associated with the derivative-free nature of objective and constraint 
functions. For example, and in contrast to the preceding sections, a 
derivative-free constrained optimization problem may involve an objective 
function $f$ for which a gradient is made available to the optimization method. 
The problem is still derivative-free if there is a constraint function defining 
the feasible region $\Omegab$ for which a (sub)gradient is not available to the 
optimization method.

As is common in many application domains where derivative-free methods are applied,
the feasible region $\Omegab$ may also involve discrete 
choices. In particular, these choices can include categorical variables that are 
either ordinal (\eg\ letter grades in $\{$A, B, C, D, F$\}$) or non-ordinal 
(\eg\ 
compiler type in $\{$flang, gfortran, ifort$\}$). Although ordinal categorical 
variables can be mapped to a subset of the reals, the same cannot be done for 
non-ordinal variables. Therefore, we generalize the formulations of
\eqref{eq:det_prob} and \eqref{eq:stoch_prob} to the problem 
\begin{equation}\label{eq:con_prob}
\begin{aligned}
 &\minimize_{\xb,\yb} && f(\xb,\yb) \\
 &\mbox{\,subject to} && \xb \in \Omegab \subset \Reals^n \\
 		 &&& \yb \in \Nonord,
\end{aligned}
\tag{CON}\end{equation}
where $\yb$ represents a vector of non-ordinal 
variables and $\Nonord$ is a finite set of feasible values.
Here we assume that discrete-valued ordinal variables are 
included in $\xb$. Furthermore, most of the 
methods we discuss do not explicitly treat non-ordinal variables $\yb$; hence, 
except where indicated, we will drop the use of $\yb$.

Similar to \secref{structured}, here we distinguish methods based on the 
assumptions made about the problem structure. We organize 
these assumptions based on the black-box optimization constraint taxonomy of 
\citeasnoun{taxonomy15}, which characterizes the type of constraint functions that 
occur in a particular specification of a derivative-free optimization problem. 
When constraints are explicitly stated (\ie\ `known' to the method), this 
taxonomy takes the 
form of the tree in \fo{Figure~\ref{fig:taxtree}}. 

\begin{figure}[!h] %
\centering
\begin{tikzpicture}[scale=0.87, every node/.style={scale=0.9},
    fact/.style={rectangle, draw=none, rounded corners=1mm, fill=blue!90, 
        text centered, anchor=north, text=white},
    state/.style={circle, draw=none, fill=orange!80, 
        text centered, anchor=north, text=white},
    leaf/.style={rectangle, draw=none, rounded corners=1mm, fill=red!30, 
        text centered, anchor=north, text=black, text depth=0pt},
    level distance=0.5cm, growth parent anchor=south
]
    \tikzstyle{every node}=[font=\footnotesize]
        \node (Fact01) [fact] {Known constraints (K)}
        child{ [sibling distance=7.2cm]
            node (State01) [state] {A or S?}
            child{
                node (Fact02) [fact] {Algebraic (A)}
                child{ [sibling distance=3.6cm]
                    node (State02) [state] {R or U?}
                    child{
                        node (Fact03) [fact] {Relaxable (R)}
                        child{ [sibling distance=1.8cm]
                            node (State041) [state] {Q or N?}
                            child{
                                node (Fact051) [fact] {\shortstack[c]{Quant.\\(Q)}}
                                child{
                                    node (State051) [leaf] {KARQ}
                                }
                            }
			  child{
                                node (Fact051) [fact] {\shortstack[c]{Nonquant.\\(N)}}
                                child{
                                    node (State051) [leaf] {KARN}
                                }
                            }
                        }}
                    child{
                        node (Fact04) [fact] {Unrelaxable (U)}
                        child{ [sibling distance=1.8cm]
                            node (State04) [state] {Q or N?}
                            child{
                                node (Fact05) [fact] {\shortstack[c]{Quant.\\(Q)}}
                                child{
                                    node (State05) [leaf] {KAUQ}
                                }
                            }
			  child{
                                node (Fact05) [fact] {\shortstack[c]{Nonquant.\\(N)}}
                                child{
                                    node (State05) [leaf] {KAUN}
                                }
                            }
		      }
                    }
                }
            }
            child{ [sibling distance=4cm]
                node (Fact10) [fact] {Black-box simulation-based (S)}
                child{ [sibling distance=3.6cm]
                    node (State02) [state] {R or U?}
                    child{
                        node (Fact03) [fact] {Relaxable (R)}
                        child{ [sibling distance=1.8cm]
                            node (State041) [state] {Q or N?}
                            child{
                                node (Fact051) [fact] {\shortstack[c]{Quant.\\(Q)}}
                                child{
                                    node (State051) [leaf] {KSRQ}
                                }
                            }
			  child{
                                node (Fact051) [fact] {\shortstack[c]{Nonquant.\\(N)}}
                                child{
                                    node (State051) [leaf] {KSRN}
                                }
                            }
                        }}
                    child{
                        node (Fact04) [fact] {Unrelaxable (U)}
                        child{ [sibling distance=1.8cm]
                            node (State04) [state] {Q or N?}
                            child{
                                node (Fact05) [fact] {\shortstack[c]{Quant.\\(Q)}}
                                child{
                                    node (State05) [leaf] {KSUQ}
                                }
                            }
			  child{
                                node (Fact05) [fact] {\shortstack[c]{Nonquant.\\(N)}}
                                child{
                                    node (State05) [leaf] {KSUN}
                                }
                            }
		      }
                    }
                }
            }
        }
;
        
\end{tikzpicture} %
\caption{Tree-based taxonomy of known (\ie\ non-hidden) constraints from
\cite{taxonomy15}.\label{fig:taxtree}}
\end{figure}

The first distinction in \fo{Figure~\ref{fig:taxtree}} is whether a constraint is 
algebraically available to the optimization method or whether it depends on a 
black-box simulation. In the context of derivative-free optimization, we will 
assume that it is these latter constraint functions for which a (sub)gradient is 
not made available to the optimization method. Algebraic constraints are those 
for which a functional form or simple projection operator is provided to the 
optimization method. \secref{const_alg} discusses methods that exclusively 
handle algebraic constraints. 
Examples of such algebraic constraints have been discussed earlier in this 
paper (\eg\ Sections~\ref{sec:RS} and~\ref{sec:convex}),
wherein it is assumed that satisfaction of the constraints (\eg\ through a 
simple projection) is trivial relative to evaluation of the objective. 
This imbalance between the ease of the constraint and objective functions is 
also the subject of recent WCC analysis \cite{CGT2018}. 

\secref{const_sim} discusses methods that target situations where one or 
more constraints do not have available derivatives.

The next distinction in \fo{Figure~\ref{fig:taxtree}} is whether a constraint can be relaxed 
or whether the constraint must be satisfied in order to obtain meaningful 
information for the objective $f$ and/or other constraint functions. 
Unrelaxable constraints are a relatively common occurrence in derivative-free 
optimization. 
In contrast to classic optimization, constraints are sometimes introduced 
solely to prevent errors in the evaluation of, for example, a simulation-based 
objective function. Methods for addressing relaxable algebraic constraints are
discussed in \secref{relax_alg}, and unrelaxable algebraic constraints are the focus of 
\secref{unrelax_alg}.

Hidden constraints are not represented in \fo{Figure~\ref{fig:taxtree}}. Hidden constraints 
are constraints that are not explicitly stated in a problem specification. 
Violating these constraints is detected only when attempting to 
evaluate the objective or constraint functions; for example, a simulation may 
fail to return output, thus leaving one of these functions undefined. Some 
derivative-free methods directly account for the possibility that such 
failures may be present despite not being explicitly stated. Hidden constraints 
have been addressed in works including those of
\citeasnoun{Avriel1967}, \citeasnoun{Kelley2000a}, \citeasnoun{Choi2000oav}, \citeasnoun{Carter2001}, \citeasnoun{Conn2001}, \citeasnoun{Huyer2008}, \citeasnoun{Leehidden11}\anlchange{, \citeasnoun{Chen2016a}, \citeasnoun{Porcelli2017} and \citeasnoun{JMMD2019}}.

\subsection{Algebraic constraints}\label{sec:const_alg}

When all constraints are algebraically available,
we can characterize the ordinal feasible region by a collection of inequality constraints:
\begin{equation}
\Omegab = \{\xb \in \Reals^n : c_i(\xb) \leq 0, \mbox{ for all } i\in \Iset \}, 
\label{eq:constraints}
\end{equation}
where each $c_i: \Reals^n \to \Reals \cup \{\infty\}$ and the set $\Iset$ is finite for all of the methods discussed. 
Problems with semi-infinite constraints can be addressed by using structured approaches as in \secref{bp}.
In this setting,
we define the constraint function $\cb:\Reals^n\to ( \Reals \cup 
\{\infty\})^{|\Iset|}$,
where the $i$th entry of the vector $\cb(\xb)$ is given by $c_i(\xb)$. 
Equality constraints can be represented \anlchange{in \eqref{eq:constraints}} 
by including both $c_i(\xb)$ and 
$-c_i(\xb)$\anlchange{; however, 
this practice should be avoided since it 
can hamper both theoretical and empirical performance.}

\subsubsection{Relaxable algebraic constraints}
\label{sec:relax_alg}

Relaxable algebraic constraints are the constraints that are typically treated in 
derivative-based non-linear optimization.
We will organize our discussion into three primary types of methods:
penalty approaches, filter approaches, and approaches with subproblems that 
employ models of the constraint functions.

\paragraph{Penalty approaches.}

Given constraints defined by \eqref{eq:constraints}, it is natural in the setting of relaxable constraints
to quantify the violation of the $i$th constraint via the value of $\max\{0,c_i(\xb)\}$. 
In fact, given a penalty parameter $\rho>0$, 
a common approach in relaxable constrained optimization is to replace the minimization of $f(\xb)$ with 
the minimization of a merit function such as
\begin{equation}
f(\xb) + \ffrac{\rho}{2} \ds\sum_{i\in\Iset} \max\{0,c_i(\xb)\}.
 \label{eq:exact_penalty}
\end{equation}

The merit function in \eqref{eq:exact_penalty} is typically called an 
\emph{exact penalty function},
because for a sufficiently large (but finite) value of $\rho>0$, every local 
minimum $\xb_*$ of \eqref{eq:con_prob} is also a
local minimum of the merit function in \eqref{eq:exact_penalty}. 
We note that each summand $\max\{0,c_i(\xb)\}$ is generally non-smooth; the 
summand is still convex provided $c_i(\xb)$ is convex.
Through the mapping
\begin{equation*}
 \Fb(\xb) = \begin{bmatrix}
 f(\xb) \\ \cb(\xb)
 \end{bmatrix}\!,
\end{equation*}
functions of the form \eqref{eq:exact_penalty} can be seen 
as cases of the composite non-smooth function \eqref{eq:composite}
and are hence amenable to the methods discussed in \secref{nso}.
In contrast to this non-smooth approach, 
a more popular merit function historically has been the \emph{quadratic 
penalty function},
\begin{equation}
 f(\xb) + \rho\ds\sum_{i\in\Iset} \max\{0,c_i(\xb)\}^2.
 \label{eq:quad_penalty}
\end{equation}
However, the merit function in \eqref{eq:quad_penalty} lacks the same 
exactness guarantees 
that come with \eqref{eq:exact_penalty}; even as $\rho$ grows arbitrarily large, 
local minima of \eqref{eq:con_prob}
need not correspond in any way with minima of \eqref{eq:quad_penalty}. 

Another popular means of maintaining the smoothness (and convexity, when 
applicable) of \eqref{eq:quad_penalty} but regaining
the exactness of \eqref{eq:exact_penalty} is to consider Lagrangian-based merit 
functions.
Associating multipliers $\lambda_i$ with each of the constraints in 
\eqref{eq:constraints}, the Lagrangian of \eqref{eq:con_prob} is
\begin{equation}
\label{eq:lag}
L(\xb;\lambdab) = f(\xb) + \sum_{i\in \Iset} \lambda_i \, c_i(\xb).
\end{equation}
Combining \eqref{eq:lag} with \eqref{eq:quad_penalty} yields the \emph{augmented 
Lagrangian} merit function
\begin{equation}
\label{eq:auglag}
L_A(\xb;\lambdab;\rho)
= f(\xb) + \sum_{i\in \Iset} \lambda_i \, c_i(\xb) 
+ \ffrac{\rho}{2} \sum_{i\in \Iset} \max\{0,c_i(\xb)\}^2 
\end{equation}
with the desired properties; that is, for non-negative $\lambdab$ and $\rho$, 
$L_A(\xb;\lambdab;\rho)$ is smooth (convex) provided that $\cb$ is.

In all varieties of these methods, which we broadly refer to as penalty approaches, 
the parameter $\rho$ is dynamically updated between iterations.
Methods typically increase $\rho$ in order to promote feasibility; penalty 
methods tend to approach solutions from outside of $\Omegab$ and hence 
typically assume that the penalized constraints are relaxable.
\anlchange{For a review on general penalty approaches, see \cite[Chapter~12]{Fletcher}.}

\citeasnoun{Lewis2002GCA} adapt the augmented Lagrangian approach of \citeasnoun{CGT91} 
in one
of the first proofs that \framework{DDS} methods can be globally 
convergent
for non-linear optimization. They utilize pattern search (see the 
discussion in 
\secref{det_det_DDS}) 
to approximately minimize the augmented Lagrangian function \eqref{eq:auglag} 
in each iteration
of their method. 
That is, each iteration of their method solves a subproblem
\begin{equation}
\label{eq:bc_auglag}
\minimize_{\xb} \{L_A(\xb;\lambdab;\rho)
: \lb\leq\xb\leq\ub \}.
\end{equation}
\citeasnoun{Lewis2002GCA} prove global convergence of their method to first-order 
Karush--Kuhn--Tucker (KKT) points. We note that the algebraic availability of 
bound constraints is explicitly used in \eqref{eq:bc_auglag}. Other constraints 
could be algebraic or simulation-based because the method used to approximately 
solve \eqref{eq:bc_auglag} does not require availability of the derivative 
$\nabla_{\xb} L_A(\xb;\lambdab;\rho)$.
The approach of \citeasnoun{Lewis2002GCA} is expanded by \citeasnoun{Lewis2010a}, who 
demonstrate the benefits of treating linear constraints (including bound 
constraints) outside of the 
augmented Lagrangian merit function.
That is, they consider subproblems of the form
\begin{equation}
\label{eq:linear_auglag}
\minimize_{\xb} \{L_A(\xb;\lambdab;\rho)
: \Ab \xb \leq \bb \}.
\end{equation}

\citeasnoun{Bueno2013} propose an inexact restoration method for problems 
\eqref{eq:con_prob} where $\Omegab$
is given by equality constraints. 
The inexact restoration method alternates between
improving feasibility (measured through the constraint violation 
$\|\cb(\xb)\|_2$ in this equality-con\-strained case)
and then approximately minimizing a $\|\cdot\|_2$-based exact penalty function 
before dynamically adjusting the penalty parameter 
$\rho$. 
Because of the separation of the feasibility and optimality phases of the inexact restoration method, 
the feasibility phase requires no evaluations of $f$. This feasibility 
phase is easier when constraint functions are available 
algebraically because (sub)derivative-based methods can be employed.
\citeasnoun{Bueno2013} prove global convergence to first-order KKT points of this 
method under appropriate assumptions. 

\citeasnoun{Amaioua2018} study the performance of a search step in 
\codes{MADS} when solving \eqref{eq:con_prob}. One of their approaches uses 
the exact penalty \eqref{eq:exact_penalty}, a second approach uses the 
augmented Lagrangian \eqref{eq:auglag} and a third combines
these two.

\citeasnoun{Audet2014b} show that the convergence properties of \codes{MADS} extend
to problems with linear equality constraints. They explicitly address these 
algebraic constraints by reformulating the original problem into a new problem 
without equality constraints (and possibly fewer variables); other constraints 
are treated as will be discussed in \secref{const_sim}.

\paragraph{Filter approaches.}
Whereas a penalty approach combines an objective function $f$ and a
measure of constraint
violation into a single merit function to be minimized approximately, a 
\emph{filter method} can be understood as a biobjective method 
minimizing the objective and the constraint violation simultaneously.
For this general discussion, we will refer to the measure of constraint 
violation as $h(\xb)$.
For example, in \eqref{eq:exact_penalty}, 
\[h(\xb) = \ds\sum_{i\in\Iset} \max\{0,c_i(\xb)\}.\]
From the perspective of biobjective optimization, a \emph{filter} can be 
understood as a subset of non-dominated points in the ($f,h$) space. 
A two-dimensional point 
$(f(\xb_l),h(\xb_l))$ 
is non-dominated, in the finite set of points $\{\xb_j:j\in J\}$ evaluated 
by a method, provided there is no $j\in J\setminus\{l\}$ with
\[f(\xb_j) \leq f(\xb_l) \quad \mbox{and} \quad h(\xb_j) \leq h(\xb_l).\]
Unlike biobjective optimization, however, filter methods adaptively vary the 
subset of non-dominated points considered in order to identify feasible points 
(\ie\ points where $h$ vanishes).
Different filter methods employ different 
mechanisms for managing the filter and generating new points.

\anlchange{\citeasnoun{Brekelmans2005} employ a filter for handling relaxable algebraic constraints.
Their model-based method attempts to have model-improving points satisfy the constraints.} 
\citeasnoun{Ferreira2017} extend the inexact-restoration method of \citeasnoun{Bueno2013} 
by replacing
the penalty formulation with a filter mechanism and again prove global 
convergence to first-order KKT points.

\paragraph{Approaches with subproblems using modelled constraints.}
Another means of constraint handling is to construct 
local models $m^{c_i}$ of each constraint $c_i$ in \eqref{eq:constraints}.
Given a local model $m^f$ of the objective function $f$, such methods 
generally employ a sequence of subproblems of the form
\begin{equation}
\label{eq:con_algsubprob}
 \minimize_{\sba} \{ m^f(\sba) : c_i(\xb+\sba)\leq 0, \mbox{ for all 
} i\in \Iset \}.
\end{equation}
As an example approach, sequential quadratic programming (SQP) 
methods 
are popular derivative-based methods that employ a quadratic model of the 
objective function and linear models of the constraint functions.
Several derivative-free approaches of this form exist, which we detail 
in this section.
We mention that many of these approaches will generally impose an additional 
trust-region constraint (\ie\ $\|\sba\|\leq\Delta$)
on \eqref{eq:con_algsubprob}.
As in \secref{det_det_TRM}, this trust-region constraint often has the 
additional role of monitoring the quality of the model $m^f$.
Furthermore, such a trust-region constraint ensures that whenever $\sba = 
\zerob$ is feasible for 
\eqref{eq:con_algsubprob}, the feasible region of \eqref{eq:con_algsubprob} 
is 
compact.

\citeasnoun{CST98} consider an adaptation of a model-based trust-region method to 
constrained problems with differentiable algebraic constraints treated via 
the trust-region subproblem \eqref{eq:con_algsubprob}. 
They target problems where they deem the algebraic 
constraints to be `easy', meaning that the resulting trust-region subproblem
is not too difficult to solve.
This method is implemented in the solver \codes{DFO} \cite{Conn2001}.

The \codes{CONDOR} method of \citeasnoun{VandenBerghen2004} and \citeasnoun{Condor} extends the
unconstrained \codes{UOBYQA} method of \citeasnoun{Powell2002a} to address algebraic
constraints. The trust-region subproblem considered takes the form
\begin{equation}
\label{eq:condorsqp}
\minimize_{\xb} \{ m^f(\xb) : 
c_i(\sba)\leq 0, \mbox{ for all 
} i\in \Iset'; \Ab \xb \leq \bb; \|\sba\|\leq\Delta \},
\end{equation}
where $m^f$ is a quadratic model 
and $\Iset'\subseteq \Iset$ captures the non-linear constraints in 
\eqref{eq:constraints}. 
In solving \eqref{eq:condorsqp}, the linear constraints are enforced explicitly 
and the non-linear constraints are addressed via an SQP approach. As will be 
discussed in \secref{unrelax_alg}, this corresponds to the linear constraints 
being treated as unrelaxable.

The \codes{LINCOA} model-based method of \citeasnoun{Powell2015} addresses linear 
inequality constraints. The \codes{LINCOA} trust-region subproblem, which can be 
seen 
as \eqref{eq:condorsqp} with $\Iset'=\emptyset$, enforces the linear 
constraints 
via an active set approach. The active set 
decreases the degrees of freedom in the 
variables by restricting $\xb$ to an affine subspace. 
Numerically efficient conjugate gradient and Krylov methods are
proposed for working in the resulting subspace. Although considerable care is 
taken to have most points satisfy the linear constraints $\Ab \xb \leq \bb$, 
these constraints are ultimately treated as relaxable, since the method does 
not enforce these constraints when attempting to improve the quality of the 
model $m^f$.

\citeasnoun{Conejo2013} propose a trust-region algorithm when $\Omegab$ is
closed and convex. They assume that it is easy to compute the 
projection onto $\Omegab$, which facilitates enforcement of the constraints
via the trust-region subproblem \eqref{eq:con_algsubprob}.
This approach is extended to include more general forms of $\Omegab$
by \citeasnoun{Conejo2015}. As with \codes{LINCOA}, although subproblem solutions 
are feasible, the constraints are treated as relaxable since they may be 
violated in the course of improving the model $m^f$.

\citeasnoun{Martinez2012} propose a feasibility restoration method intended for
problems with inequality constraints where the feasible region is `thin': for
example, if $\Omegab$ is defined by both $c_i(\xb) \le 0$ and $-c_i(\xb) \le 0$
for some $i$. Each iteration contains two steps: one that seeks to minimize the
objective and one that seeks to decrease infeasibility using many evaluation of
the constraint functions (without evaluating the objective). 
Similar to the progressive-barrier method discussed in \secref{const_sim}, the
method by \citeasnoun{Martinez2012} dynamically updates a tolerable level of
infeasibility.

\subsubsection{Unrelaxable algebraic constraints}
\label{sec:unrelax_alg}

We now address the case when all of the constraints are available 
algebraically but an unrelaxable constraint also exists. In this 
setting, such unrelaxable constraints are typically necessary to ensure 
meaningful output of a black-box objective function. Consequently, methods must 
always maintain feasibility (or at least establish feasibility 
and then maintain it) with respect to the unrelaxable constraints.

An early example of a method for unrelaxable constraints is the `complex' 
method of \citeasnoun{Box1965}. This extension of the simplex method of 
\citeasnoun{Spendley1962} treats 
unrelaxable bound constraints by modifying the simplex operations to 
project into the interior of any potentially violated bound constraint.
\citename{May1974} \citeyear{May1974,May1979} extends
the unconstrained 
derivative-free method of \citeasnoun{Mifflin1975} to address unrelaxable linear 
constraints. The method of \citeasnoun{May1979} uses finite-difference estimates, 
but care is taken to ensure that the perturbed points never violate the 
constraints.

As seen in \secref{relax_alg}, several approaches treat 
non-linear algebraic constraints via a merit function and 
enforce unrelaxable linear constraints via a constrained subproblem. These 
include the works of
\citeasnoun{Lewis2002GCA} for bound constraints in \eqref{eq:bc_auglag},
\citeasnoun{Lewis2010a} for inequality constraints in \eqref{eq:linear_auglag}, and
\citeasnoun{VandenBerghen2004} for inequality constraints in 
\eqref{eq:condorsqp}. 
Another merit function relevant for unrelaxable constraints is the 
extended-value merit function 
\begin{equation}
h(\xb) = f(\xb)+\infty \, \delta_{\Omegab^C}(\xb),
\label{eq:extvalue}
\end{equation}
where $\delta_{\Omegab^C}$ is the indicator function of $\Omegab^C$.
Such an \emph{extreme-barrier} approach (see \eg\ the 
discussion by \citeb{Lewis1999}) is particularly relevant for simulation-based 
constraints. Hence, with the exception of explicit treatment of unrelaxable 
algebraic constraints, we postpone significant discussion of extreme-barrier 
methods until \secref{const_sim}.

\paragraph{DDS methods for unrelaxable algebraic constraints.}
Within \framework{DDS} methods, an intuitive approach to handling unrelaxable constraints is to
limit poll directions $\Db_k$ so that $\xb_k+\db_k$ is feasible with respect to the unrelaxable constraints. 
\citeasnoun{Lewis1999} and \citeasnoun{LuSc02}, respectively, develop pattern-search and 
coordinate-search methods for unrelaxable bound-constrained problems.
By modifying the polling 
directions \citeasnoun{Lewis2000a} show that pattern-search methods are also 
convergent in the presence of unrelaxable linear constraints. 
\citeasnoun{Chandramouli2019} address unrelaxable bound 
constraints 
within a \framework{DDS} method that employs a model-based method in the search step
in addition to a bound-constrained line-search step.
\citeasnoun{Kolda2006} develop and analyse a new condition, related to the tangent 
cone of nearby active constraints, on the
sets of directions used within a generating set search method when solving 
linearly constrained problems. The condition ensures that evaluated points are 
guaranteed to satisfy the linear constraints.
 \citeasnoun{lucidi99objectivederivativefree} propose feasible descent methods
that sample the objective over a finite set of search directions. Each iteration
considers a set of $\epsilon$-active constraints (\ie\ those constraints for 
which $c_i(\xb_k)\geq - \epsilon$) for general algebraic inequality constraints. 
Poll steps are projected in
order to ensure they are feasible with respect to these $\epsilon$-active
constraints. The analysis of \citeasnoun{lucidi99objectivederivativefree} extends 
that of \citeasnoun{Lewis2000a} and establishes convergence
to a first-order KKT point under standard assumptions. 

As introduced in \secref{relax_alg}, \citeasnoun{Audet2014b} reformulate 
optimization problems with unrelaxable linear equality constraints in the 
context of \framework{MADS}. 

\citeasnoun{Gratton2017direct} extend the randomized \framework{DDS} method of 
\citeasnoun{Gratton2015} 
to linearly constrained problems; candidate points are accepted only
if they are feasible. 
\citeasnoun{Gratton2017direct} establish probabilistic convergence and complexity results using a 
stationary measure appropriate for linearly constrained problems.

\paragraph{Model-based methods for unrelaxable algebraic constraints.}
\anlchange{Model-based}\linebreak[4]
methods are more challenging to design in the presence of 
unrelaxable constraints because enforcing guarantees of model 
quality such as those in \eqref{eq:fullylinear} can be difficult. 
For a fixed value of $\kappab$ in \eqref{eq:fullylinear}, it 
may be impossible to obtain a $\kappab$-fully linear model using only 
feasible points. As an example, consider two linear constraints for which 
the angle between the constraints is too small to allow for $\kappab$-fully 
linear model construction; avoiding interpolation points drawn from
such thin regions motivated development of 
the wedge-based method of \citeasnoun{Marazzi2002a} from \secref{det_det_TRM}.

\citeasnoun{Powell2009a} proposes \codes{BOBYQA}, a model-based trust-region method 
for bound-constrained optimization without derivatives that extends the 
unconstrained method \codes{NEWUOA} in \citeasnoun{MJDP06}. \codes{BOBYQA}
 ensures that all points at which $f$ is evaluated satisfy the bound 
constraints.
\citeasnoun{Arouxet2011} modify \codes{BOBYQA} to use 
an active-set strategy in solving the bound-constrained trust-region 
subproblems; an $\|\cdot\|_{\infty}$-trust region is employed so that these 
subproblems correspond to minimization of a quadratic over a 
compact, bound-constrained domain. 
\citeasnoun[Section~6.3]{WildPhD} develops an RBF-model-based method for 
unrelaxable bound constraints by enforcing the bounds during both model 
improvement and $\|\cdot\|_{\infty}$-trust-region subproblems. 
\citeasnoun{Gumma2014} extend the \codes{NEWUOA} method to address 
linearly constrained 
problems. The linear constraints are enforced both when solving the trust-region subproblem
and when seeking to improve the geometry of the interpolation points.
\citeasnoun{Gratton2011} propose a model-based method for 
unrelaxable bound-constrained
optimization, which restricts the construction of fully linear models
to subspaces defined by nearly active constraints. 
Working in such a reduced space means that the machinery for unconstrained 
models in \secref{det_det_TRM} again applies.

\paragraph{Methods for problems with unrelaxable discrete constraints.}
Constraints that certain variables take discrete values are often unrelaxable in
derivative-free optimization. For example, a black-box simulation may be unable 
to assign meaningful output when input variables take non-integer values. 
That such integer constraints are unrelaxable presents challenges distinct from 
those typically arising in mixed-integer non-linear 
optimization \cite{Belotti2013}.

Naturally, researchers have modified derivative-free methods for continuous 
optimization to address integer constraints.
\citeasnoun{Audet:2000} and \citeasnoun{Abramson:2009}, respectively, propose
integer-constrained pattern-search and \framework{MADS} methods to 
ensure that evaluated points respect integer constraints.
\citeasnoun{Abramson2007} develop a pattern-search method that employs a filter 
that handles general inequality constraints and ensures that 
integer-constrained variables are always integer.

\citeasnoun{Porcelli2017} propose the `brute-force optimizer' \codes{BFO}, a 
\framework{DDS} method for mixed-variable problems (including those with 
ordinal categorical variables) that aligns the poll points to respect the 
discrete constraints. A recursive call of the method reduces the number of 
discrete variables by fixing a subset of these variables.

\citeasnoun{Liuzzi2011} \anlchange{solve} 
mixed-integer bound-constrained problems by using 
both a local discrete search (to address integer variables) and a line
search (for continuous variables).
This approach is extended by \citeasnoun{Liuzzi2015} to also address mixed-integer
problems with general constraints using the SQP approach from 
\citeasnoun{Liuzzi2010}.
\citeasnoun{Liuzzi2018} solve \anlchange{constrained} integer optimization problems by
performing non-monotone line searches along feasible \emph{primitive directions} 
$\Db$ in a neighbourhood of the current point $\xb_k$. Feasible primitive 
directions are those
$\db\in \Ints^n \cap \Omegab$ satisfying 
$\operatorname{GCD}(\db_1,\ldots,\db_{| \Db |})=1$, that is, 
directions in a bounded neighbourhood that are not integer multiples of one 
another.

The method by \citeasnoun{Rashid12} for mixed-integer problems builds multiquadric
RBF models. 
Candidate points are produced by using 
gradient-based mixed-integer optimization techniques; the authors' relaxation-based 
approach employs a `proxy model' that coincides with 
function values from points satisfying the unrelaxable integer 
constraints. 
\anlchange{The methods of 
\citename{Mueller2012} \citeyear{Mueller2013,Mueller2012} and \citen{muller2016miso} similarly employ 
a global RBF model over the integer lattice, with various strategies for generating
trial points based on this model.}
\citeasnoun{Newby2015} build on \codes{BOBYQA} to address bound-constrained
mixed-integer problems. 
They outline an approach for building interpolation models of objectives using
only points that are feasible. Their trust-region subproblems consist of
minimizing a quadratic objective subject to bound and integer constraints.

Many of the discussed methods have been shown to converge to points that are
mesh-isolated local solutions;
see 
\citeasnoun{Newby2015} for discussion of such `local minimizers'. When an 
objective is convex, one can do better. 
\citeasnoun{Palkar2019} propose a method for certifying a global minimum of a
convex objective $f$ subject to unrelaxable integer constraints. They form a
piecewise linear underestimator by interpolating $f$ through subsets of $n+1$ 
affinely independent points. The resulting underestimator is then used to 
generate new candidate points until global optimality has been certified.

\subsection{Simulation-based constraints} \label{sec:const_sim} 

As opposed to the preceding section, methods in this section are not limited to
constraints that have closed-form solutions but also address constraints that depend
on the output from some calculated function. 
Many methods address such simulation-based constraints by using approaches similar to
those used for algebraic constraints. 

\paragraph{Filter approaches.}
Filter methods for simulation-based constraints, as with algebraic constraints, 
seek to simultaneously decrease the objective and constraint violation.
For example, \citeasnoun{Audet2004} develop a pattern-search method for general
constrained optimization that accepts steps that improve either the objective or
some measure of violation of simulation-based constraints. Their hybrid approach 
applies an extreme barrier to points that violate linear or bound constraints. 
\citeasnoun{Audet2004b} provides examples where the method by \citeasnoun{Audet2004}
does not converge to stationary points.

\citeasnoun{PourmohamadPhD} models objective and constraint functions using Gaussian 
process models in a filter-based method. Because these models are stochastic, 
point acceptability is determined by criteria such as probability of filter
acceptability or expected area of dominated region (in the filter space).

\citeasnoun{Echebest2015} develop a derivative-free method in the inexact 
feasibility restoration filter method framework of \citeasnoun{Gonzaga2004}. 
\citeasnoun{Echebest2015} employ fully linear models of the objective and 
constraint functions and show that the resulting limit points are first-order 
KKT points.

\paragraph{Penalty approaches.}
The original \codes{MADS} method \cite{Audet06mads} converts constrained
problems to unconstrained problems by using the extreme-barrier approach 
mentioned above;
that is, the merit function \eqref{eq:extvalue}
effectively assigns a value of infinity to points that violate any constraint.
A similar approach to general constraints is used by the `complex' method of 
\citeasnoun{Box1965}; the simplex
(complex) is updated to maintain the feasibility of the  \anlchange{vertices} of the simplex.
As a consequence of its generality, the extreme-barrier approach is applicable 
for algebraic constraints, simulation-based constraints and even hidden 
constraints. Furthermore, because \eqref{eq:extvalue} is independent of the 
degree of both constraint satisfaction and constraint violation, the 
extreme barrier is able to address non-quantifiable constraints.

In contrast, the progressive-barrier method by \citeasnoun{Audet2009} 
employs a quadratic constraint penalty similar to \eqref{eq:quad_penalty} for
relaxable simulation-based constraints $\{c_i:i\in\Iset_r\}$ and an
extreme-barrier penalty for unrelaxable simulation-based constraints
$\{c_i:i\in\Iset_u\}$.
Their progressive-barrier method maintains a non-increasing 
threshold value $\epsilon_k$ that quantifies the 
allowable relaxable constraint violation in each iteration. 
Their approach effectively uses the merit function
\begin{equation}
h_k(\xb) = \begin{cases}
 f(\xb) & \mbox{if $ \xb\in\Omegab_u $ and 
 $\sum_{i\in\Iset_r} \max\{0,c_i(\xb)\}^2 < \eps_k $}, \\
	 \infty \, & \mbox{otherwise},\\
 \end{cases}
\label{eq:progvalue}
\end{equation}
where $\Omegab_u = \{\xb: c_i(\xb) \le 0, \; \forall i \in \Iset_u\}$ denotes 
the feasible domain with respect to the unrelaxable constraints. 
The progressive-barrier method maintains a set of feasible and infeasible 
incumbent
points and seeks to decrease the threshold $\epsilon_k$ to 0
based on the value
of \eqref{eq:progvalue} at infeasible incumbent points. 
Trial steps are accepted as incumbents based on criteria resembling, but
distinct from, those used by filter methods.
Convergence to Clarke stationary points is obtained for particular sequences
of the incumbent points.
The \codes{NOMAD} \cite{LeDigabel:2011} implementation of \framework{MADS} 
allows users to
choose to address inequality constraints handled via extreme-barrier,
progressive-barrier or filter approaches. 

Also within the \framework{DDS} framework, \citeasnoun{Gratton2014} use an 
extreme-barrier
approach to handle unrelaxable constraints and an exact 
penalty function to handle the relaxable constraints. That is, step 
acceptability is based on satisfaction of the unrelaxable constraints as well 
as sufficient decrease in the merit function \eqref{eq:exact_penalty}, with 
the set $\Iset$ containing only those constraints that are relaxable. As the 
algorithm progresses, relaxable constraints are transferred to the set of 
constraints treated by the extreme barrier; this approach is similar to that 
underlying the progressive-barrier approach.

\citeasnoun{Liuzzi2009} and \citeasnoun{Liuzzi2010} consider line-search methods that
apply a penalty to simulation-based constraints;
\citeasnoun{Liuzzi2009} employ an exact penalty function (a smoothed version of 
$\|\cdot\|_{\infty}$), whereas \citeasnoun{Liuzzi2010} employ a sequence 
of
quadratic penalty functions of the form \eqref{eq:quad_penalty}. 
\anlchange{\citen{Fasano2014} propose a similar line-search approach to address 
constraint and objective functions that are not differentiable.}
 
Primarily concerned with equality constraints,
\citename{Sampaio2015derivative} \citeyear{Sampaio2015derivative,Sampaio2016numerical}
propose a
derivative-free variant of trust-funnel methods, a class
of methods proposed by \citeasnoun{Gould2010} that avoid the use of both merit
functions and filters.

\citeasnoun{Diniz-Ehrhardt2011} propose a method that models objective and 
constraint functions in an augmented Lagrangian framework. 
Similarly, \citeasnoun{ALSlack16} use an augmented Lagrangian framework, 
wherein the merit function in \eqref{eq:auglag} uses 
Gaussian process models of the objective and 
constraint functions in place of the actual objective and constraint 
functions.

\paragraph{Approaches with subproblems using modelled constraints.}

In early work, \citeasnoun{Glass1965} \anlchange{develop} a coordinate-search method that also 
uses linear models of the objective and constraint functions. On each iteration, 
after the coordinate directions are polled, the models are used in a linear 
program to generate new points; points are accepted only if they are feasible.
Extending this idea,
\citeasnoun{Powell1994} \anlchange{develops} the constrained optimization by 
linear approximation (\codes{COBYLA}) method, which builds linear 
interpolation models of
 the objective and constraint functions on a common set of $n+1$ affinely 
independent points. Care is taken to maintain the non-degeneracy of this 
simplex. The method can handle both inequality and equality constraints, with 
candidate points obtained from a linearly constrained subproblem 
and then accepted based on a merit function of the form 
\eqref{eq:exact_penalty}.

\citeasnoun{Brmen2015} propose a variant of \codes{MADS} with a specialized 
model-based search step
\begin{equation}
\label{eq:brmensp}
\minimize_{\xb} \{ m^f(\xb) : \Ab \xb \leq \bb \},
\end{equation}
where $m^f$ is a strongly convex quadratic model of $f$ and ($\Ab, \bb$) are 
determined 
from linear regression models of the constraint functions. Both the search and 
poll 
steps are accepted only if they are feasible; this corresponds to the 
method effectively treating the constraints with an extreme-barrier approach.

\citeasnoun{Gramacy2015} extend the \framework{MADS} framework by using treed 
Gaussian
processes to model both the objective and simulation-based
constraint functions. The resulting models are used both within the search step 
and to order the poll points (within an opportunistic polling paradigm) using a 
filter-based approach.

A number of methods work with restrictions of the domain $\Omegab$ in order to 
promote feasibility (typically with respect to the simulation-based 
constraints) of the generated points. Such strategies are often motivated by 
a desire to avoid the situation where feasibility is established only 
asymptotically. An example of such a restricted domain is the set
\begin{equation}
\label{eq:restricted_omega}
 \Omegab_{\rm res}(\epsb) = \{\xb \in 
\Reals^n : c_i (\xb) \leq 0 \; \forall i\in 
\Iset_a, \;
 m^{c_i}(\xb) + \epsilon_i(\xb) \leq 0 \; \forall i\in 
\Iset_s \},
\end{equation}
where algebraic constraints (corresponding to $i\in 
\Iset_a$) are explicitly enforced and a parameter (or function of $\xb$) 
$\epsb$ controls the degree of restriction for the modelled simulation-based 
constraints
(corresponding to $i\in \Iset_s$).

The methods of \citeasnoun{Regis2013b} utilize interpolating radial basis 
function surrogates
of the objective and constraint functions. Acceptance of infeasible 
points is allowed and is followed by a constraint restoration phase that 
minimizes a quadratic penalty based on the modelled constraint violation. When 
the current point is feasible, a subproblem is solved with a feasible set defined by
\eqref{eq:restricted_omega} in addition to a constraint that lower-bounds the distance 
between the trial point and the current point. Each parameter $\epsilon_i$ is 
adjusted based on the feasibility of constraint $i\in \Iset_s$ in recent 
iterations.

\citeasnoun{NOWPAC2014} develop a trust-region
 method employing fully linear models of both constraint and objective 
functions.
 They introduce a \emph{path augmentation} scheme intended to locally
 convexify the simulation-based constraints.
Their trust-region subproblem at the current point $\xb_k$ minimizes the model 
of the objective function subject to a trust-region constraint and the 
restricted feasible set \eqref{eq:restricted_omega}, where
$\epsilon_i(\xb) = \epsilon_0 \|\xb-\xb_k\|^{2/(1+p)}$ and where 
$\epsilon_0>0$ and $p\in (0,1)$ are fixed constants.
 \citeasnoun{NOWPAC2014} establish convergence of their method from feasible
 starting points; that is, they show a first-order criticality measure
 asymptotically tends to 0. \citeasnoun{NOWPAC2014} produce a code,
 \codes{NOWPAC}, that employs minimum-Frobenius norm quadratic models of
 both the objective and constraint functions. This work is extended 
by
 \citeasnoun{Augustin2017} to the stochastic optimization problem 
\eqref{eq:stoch_prob}.

Whereas \citeasnoun{NOWPAC2014} consider a local
 convexification of inequality constraints through the addition of a convex
 function to the constraint models, \citeasnoun{CONORBIT15} consider a similar
 model-based approach but define an envelope around models of nearly active
 constraints. In particular, at the current point $\xb_k$, the restricted 
feasible set 
 \eqref{eq:restricted_omega} uses the parameter
 \[ 
 \epsilon_i(\xb) = \begin{cases}
 0 & \mbox{if } c_i(\xb_k) >-\xi_i, \\
 \xi_i & \mbox{if } c_i(\xb_k) \leq -\xi_i,
 \end{cases}
 \]
 where $\{\xi_i : i\in \Iset_s\}$ is fixed. This form of $\epsb$ ensures that 
trust-region subproblems remain non-empty and avoids applying a restriction 
when the algorithm is sufficiently close to the level set $\{\xb : c_i(\xb) = 
0\}$.
 
\citeasnoun{Troltzsch2016} considers an SQP method in the style of 
\citeasnoun{Omojokun1990}, which applies a two-step process that first seeks to 
improve a measure of constraint violation and then solves a subproblem 
restricted to the null space of modelled constraint gradients.
\citeasnoun{Troltzsch2016} uses linear models of the constraint functions and 
quadratic models of the objective function, with these models replacing $\cb$ 
and $f$ in the augmented Lagrangian merit function in \eqref{eq:auglag}.
Step acceptance uses a merit function (an exact penalty function). 

\anlchange{\citeasnoun{Mueller2017} develop a method for addressing
computationally inexpensive objectives while satisfying computationally expensive
constraints. Their two-phase method first seeks feasibility 
by solving a multi-objective optimization problem \anlchange{(a problem class that is the subject of \secref{sec84})} in which the constraint violations are
minimized simultaneously; the second phase seeks to reduce the objective subject to 
constraints derived from cubic RBF models of the constraint functions.}

\citeasnoun{Bajaj2018} propose a two-phase method. In the feasibility phase, a
trust-region method is applied to a quadratic penalty function 
that employs models of the simulation-based constraints. The trust-region 
subproblem at iteration $k$ takes the form
\begin{equation}
\label{eq:bsp1}
\minimize_{\xb} \Biggl\{
\sum_{i\in \Iset_s} \max\{0, m^{c_i}(\xb)\}^2 : c_i(\xb) \leq 0 
\; \forall i\in \Iset_a, \; \xb \in \B{\xb_k}{\Delta_k} \Biggr\}
\end{equation}
and thus explicitly enforces the algebraic constraints ($i\in 
\Iset_a$) and penalizes violation of the modelled simulation-based constraints
($i\in \Iset_s$).
In the optimality phase, a trust-region method is applied to a model of the 
objective function, and the modelled constraint violation is 
bounded by that achieved in the feasibility phase; that is, the 
trust-region subproblem is
\begin{equation*} %
\begin{aligned}
 &\minimize_{\xb} && m^f(\xb) \\
 &\mbox{\,subject to} && c_i(\xb) \leq 0 \quad \gforall\ i\in 
\Iset_a \\
 &&& m^{c_i}(\xb) \leq c_i(\xb_{\rm pen}) \quad \gforall\ i\in 
\Iset_s \\
&&& \xb \in \B{\xb_k}{\Delta_k},
\end{aligned}
\end{equation*}
where $\xb_{\rm pen}$ is the point returned from the feasibility phase.

\anlchange{\citeasnoun{Hare2005} present an approach for approximating the normal
and tangent cones; their approach is quite general and applies to the 
case when the domain is defined by non-quantifiable black-box constraints. 
\citeasnoun{Davis2013} consider a simplex-gradient-based approach for 
approximating normal cones when the black-box constraints are quantifiable. 
Naturally, such approximate cones could be used to determine if a method's
candidate solution approximately satisfies a stationarity condition.}

\vspace{5pt}
\section{Other extensions and practical considerations} %
\label{sec:extensions}
We conclude with a cursory look at extensions of the methods 
presented, especially highlighting active areas of development.

\subsection{Methods allowing for concurrent function evaluations}
\label{sec:concurrent}
A number of the methods presented in this survey readily allow for the
concurrent evaluation of the objective function at 
\anlchange{multiple points $\xb\in \Reals^n$.} 
\anlchange{Performing function evaluations concurrently through the} 
use of parallel computing resources should decrease the wall-clock time
required by a given method\anlchange{.}
\anlchange{Depending on the method, there is a natural limit} to the amount of
concurrency that can be utilized efficiently.
\anlchange{Below we summarize such methods and their limits for concurrency.}

The simplex methods discussed in \secref{det_det_DS} benefit from
performing $n$ concurrent evaluations of the objective when a shrink operation
is performed. Also, the points corresponding to the expansion and reflection
operations could be evaluated \anlchange{in parallel}. 
\anlchange{Non-opportunistic d}irectional direct-search method\anlchange{s} are especially amenable to parallelization
\cite{DennisJr1991} because the $| \Db_k |$ poll directions 
can be evaluated concurrently. 

Model-based methods from \secref{det_det_model} can use concurrent evaluations
during model building when, for example, evaluating \anlchange{up to} $\dim{\cP^{d,n}}$ \anlchange{additional points for use} in
\eqref{eq:interpolation_system}.
In another example, \codes{CONDOR} \cite{Condor,VandenBerghen2004}
utilizes concurrent
evaluations of the objective to replace points far away from the current
trust-region centre by maximizing the associated Lagrange polynomial.
The thesis by \citeasnoun{Olsson2014} considers three ways of using concurrent
resources within a model-based algorithm: using multiple starting points,
evaluating different models in order to better predict a point's value, and 
generating multiple points with each model (\eg\ solving with the trust-region
subproblem with different radii).
\anlchange{A similar approach of generating multiple trial points concurrently 
is employed in the parallel direct-search, trust-region method of \citeasnoun{HoughMeza02}.}

Finite\anlchange{-difference-based approaches (\eg\ }\secref{det_det_others}) 
allow for $n$ concurrent
evaluations \anlchange{with} forward differences \eqref{eq:fdiff} or $2n$ concurrent
evaluations with central differences \eqref{eq:cdiff}. Implicit filtering also
performs such a central-difference calculation that can utilize $2n$ concurrent
evaluations (\lineref{stencil} of \algref{imfil}). 
\anlchange{L}ine-search methods 
\anlchange{can evaluate multiple
points concurrently} during their line-search procedure. 
The methods of \citeasnoun{Garcia-Palomares2002} and \citeasnoun{Garcia-Palomares2013}
also consider using parallel resources to concurrently evaluate points in 
\anlchange{a neighbourhood of interest.}

\anlchange{When} using a set of independently generated points\anlchange{, pure random search exhibits}
perfect scaling as the level of available concurrency increases. 
Other\-wise, the randomized methods for deterministic objectives from
\secref{rand_det} can utilize concurrent \anlchange{evaluations} in a 
manner similar to that of their deterministic counterparts. Nesterov random search can 
use $n$ or $2n$ concurrent objective evaluations when computing \anlchange{an} approximate
gradient in \eqref{gradient_free_oracle}. Randomized DDS methods can
concurrently evaluate $| \Db_k |$ poll points, and randomized
trust-region methods can concurrently evaluate points needed for building and
improving models.

In addition to the above approaches for using parallel resources, methods from
\secref{structured}
for structured problems can use concurrent evaluations to calculate \anlchange{parts} of
the objective. For example, methods for optimizing \anlchange{separable objectives such as} 
\eqref{eq:nls} or \eqref{eq:partially_separable} can evaluate the $p$ \anlchange{component} 
functions $F_i$ concurrently. 

The various gradient approximations used by methods in \secref{stoch} are
amenable to parallelization in the same manner as previously discussed\anlchange{, but with the additional possibility of also evaluating at multiple $\xib$ values.
\framework{SA} methods} can use $2n$ concurrent evaluations \anlchange{of $\tilde{f}$} in calculating \eqref{eq:fdg} or
\eqref{eq:fdg_crn} and \codes{SPSA} can use \anlchange{two} concurrent evaluations when
calculating \eqref{eq:spsa_est}. Methods employing the sample mean estimator
\eqref{eq:samplemean} can utilize $p_k$ evaluations concurrently.

\subsection{Multistart methods}
\label{sec:multistart}

\anlchange{A natural approach for addressing non-convex objectives for which}
it is not known whether multiple local 
minima exist \anlchange{is to start} a local optimization method
from different points in the domain in the hope of identifying different local
minima. Such multistart approaches also allow for the use of methods that are
specialized for optimizing problems with known structure.

Multistart methods allow for the use of concurrent objective
evaluations if two or more local \anlchange{optimization run}s are being performed at the
same time. \anlchange{Multistart methods also allow one} to utilize additional computational
resources\anlchange{; this ability is especially useful when} an objective evaluation does not become faster with \anlchange{additional} resources or
when the local optimization method is inherently sequential. 

\anlchange{\citeasnoun{Boender1982}} derive confidence intervals on the
\anlchange{objective value of a global minimizer} when starting a local optimization method at uniformly
drawn points. Their 
\anlchange{analysis} 
 gives rise
to the multilevel single linkage (\codes{MLSL}) method \cite{RinnooyKan1987,RinnooyKan1987a}.
Iteration $k$ of the method draws $N$ points uniformly over the domain and
starts a local optimization method from sampled points that do not have any
other point within a \anlchange{specific distance, depending on $k$ and $N$,} 
with a smaller objective value. With this rule, and under assumptions on the
distance between minimizers in $\Omegab$ and properties of the local
optimization method used, \codes{MLSL} is shown to \anlchange{almost surely} identify all local minima
while starting the local optimization method from only finitely many points.
\citename{LarWild14} \citeyear{LarWild14,LW16}
generalize \codes{MLSL} by showing similar
theoretical results when 
\anlchange{starting-point selection utilizes points both from the random sampling and from those generated by local optimization runs.}

If a meaningful variance exists in the objective evaluation times, batched
evaluation of points may result in an inefficient use of computational
resources. Such concerns have motivated the development of a number of methods
including the \codes{HOPSPACK} framework \cite{Planteng2009}, which supports the
sharing of information between different local optimization methods.
\anlchange{\citeasnoun{CASRGR03} also use information from multiple 
optimization methods to determine points at which to evaluate the objective function.}
Similarly, \anlchange{the \codes{SNOBFIT}}
method by \citeasnoun{Huyer2008} uses concurrent objective evaluations while combining
local searches in a global framework. The software focuses on robustness in
addressing many practical concerns including soft constraints, hidden
constraints, and a problem domain that is modified by the user as the method
progresses.

Instead of coordinating concurrent instances of a pattern-search method,
\citeasnoun{Audet2008} propose an implementation of \codes{MADS} that decomposes
the domain into subspaces to be optimized over in parallel. \citeasnoun{Alarie2018}
study different approaches for selecting subsets of variables to define
subproblems in such an \anlchange{approach. 
\citeasnoun{Custodio2015}} maintain concurrent instances of a pattern-search method,
and merge those instances that \anlchange{become} sufficiently close.
\citeasnoun{Taddy2009} use a global treed-Gaussian process to guide \anlchange{a 
local pattern-search method to encourage the identification of} better local minima.

\subsection{Other global optimization methods} 
\label{sec:global} 

Guarantees of global optimality for \anlchange{general continuous} functions rely on   \anlchange{candidate points being generated densely} in the domain \cite[Theorem~1.3]{ATAZ89}; such  \anlchange{candidate points} can
be generated in either a deterministic or randomized fashion.
\anlchange{When} $f$ is Lipschitz-continuous on $\Omegab$ and the Lipschitz constant $\Lip{f}$
is available to the optimization method, one need not generate 
 \anlchange{points densely} in the domain. \anlchange{In particular, if} $\hat{\xb}$ is an approximate minimizer of $f$
and 
\anlchange{$\xb$}
is a point satisfying 
\anlchange{$f(\xb)>f(\hat{\xb})$,}
no global
minimizer can lie in -- and therefore no point needs to be sampled
from\anlchange{~--}
\anlchange{$\B{\xb}{\gpfrac{f(\xb)-f(\hat{\xb})}{\Lip{f}}}$.}
Naturally, the \anlchange{benefit of exploiting this fact}
requires accurate knowledge of the
Lipschitz constant. One can empirically observe a lower bound on $\Lip{f}$, but
\anlchange{obtaining useful upper bounds on}
$\Lip{f}$ may not be possible. \anlchange{Methods that exploit this Lipschitz knowledge} 
may suffer a considerable performance decrease when overestimating $\Lip{f}$
\cite{Hansen1991}.

Motivated by situations where the Lipschitz constant of $f$ is
unavailable, \citeasnoun{jones1993lipschitzian} develop the \codes{DIRECT} (DIviding
RECTangles) method. 
\codes{DIRECT} partitions
a bound-constrained $\Omegab$ into $2n+1$ hyper-rectangles (hence the method's 
name) with an
evaluated point at the centre of each. 
\anlchange{Each hyper-rectangle is scored via}
a combination of the length of its longest side and the function value at its \anlchange{centre}.
\anlchange{This scoring favours hyper-rectangles exhibiting both long sides and small function values;
the best-scoring hyper-rectangles}
are further divided. (As such, \codes{DIRECT}'s performance can be \anlchange{significantly} affected
by adding a constant value to the objective \cite{Finkel2006}.)
\codes{DIRECT} \anlchange{generates centres} that are dense in 
\anlchange{$\Omegab$}
and will therefore identify the global \anlchange{minimizers} of $f$ over 
\anlchange{$\Omegab$}, even when $f$ is 
\anlchange{non-smooth}
\cite{jones1993lipschitzian,Finkel2004,Finkel2009}.
\anlchange{Several versions of \codes{DIRECT} that perform concurrent function evaluations} take significant care to
ensure the sequence of points generated is the same as that produced by
\codes{DIRECT} \cite{He2009,He2007,He2009a,He2009b}.
Similar \anlchange{hyper-rectangle} partitioning strategies are used by the \anlchange{methods of \citeasnoun{Munos2011}.}
The multilevel coordinate-search (\codes{MCS}) method by \citeasnoun{Huyer1999} is
inspired by \codes{DIRECT} in many ways. \codes{MCS} maintains a 
\anlchange{partitioning} of
the domain and subdivides \anlchange{hyper-rectangles} based on their size and value.
\codes{MCS} uses \anlchange{the function values at} boundary points\anlchange{, rather than the centre points,} to determine the value of a \anlchange{hyper-rectangle}; such 
\anlchange{boundary points}
can be \anlchange{shared by} more than one \anlchange{hyper-rectangle}. \citeasnoun{Huyer2008} show that a
\anlchange{version}
of \codes{MCS} needs
to consider only finitely many \anlchange{hyper-rectangles} before identifying a global minimizer.

Many randomized approaches for generating points densely in 
\anlchange{a domain $\Omegab$}
have been developed. These include Bayesian optimization methods and related variants 
\anlchange{\cite{Mockus1989,Jones1998,Frazier2018}}, some of which have established
complexity rates \cite{Bull2011}. 
Such 
\anlchange{randomized samplings of $\Omegab$}
can be
used to produce a global surrogate; similar to other model-based methods, this
global model can be minimized to produce points where the objective 
\anlchange{should}
be evaluated. Although minimizing such a global surrogate may be difficult, such a
subproblem 
\anlchange{may}
be easier than the original problem, which typically entails 
a computationally expensive objective function for which derivatives are 
unavailable. \citeasnoun{Vu2016} provide a recent survey of such
surrogate-based methods for global optimization.

\subsection{Methods for multi-objective optimization}\label{sec:sec84}  
 
Multi-objective optimization problems are typically stated as
\begin{equation}
\begin{aligned}\label{eq:moo_prob}
 &\underset{\xb}{\minimize} && \Fb(\xb) \\
 &\mbox{\,subject to} && \xb \in \Omegab \subset \Reals^n,
\end{aligned}
\tag{MOO}\end{equation}
where $p>1$ objective functions $f_i:\Reals^n\to\Reals$ for $i=1,\dots,p$\anlchange{}
define the vector-valued mapping $\Fb$ via 
$\Fb(\xb)=[f_1(\xb), \, \ldots, \, f_p(\xb)]$.
Given potentially conflicting 
objectives $f_1,\dots,f_p$, \anlchange{the problem \eqref{eq:moo_prob} is well-defined only when given} an 
ordering on the vector of objective values $\Fb(\xb)$. 
Given distinct points $\xb_1,\xb_2\in\Reals^n$, $\xb_1$ \emph{\anlchange{Pareto dominates}} $\xb_2$ 
provided
\[
f_i(\xb_1)\leq f_i(\xb_2) \quad \mbox{for all } i=1,\dots,\anlchange{p \qquad \mbox{and} \quad f_j(\xb_1) < f_j(\xb_2) \mbox{ for some } j}.
\]
The set of all \anlchange{feasible} points that are not \anlchange{Pareto-}dominated by any other \anlchange{feasible} point
is referred to
as the \emph{Pareto\anlchange{(-optimal)} set of} \eqref{eq:moo_prob}. \anlchange{An in-depth treatment of such problems is provided by \citeasnoun{Ehrgo00b}.}

Ideally, a method designed for the solution of \eqref{eq:moo_prob} should return \anlchange{an}  approximation of the
Pareto \anlchange{set.} If at least one objective $f_1,\dots,\anlchange{f_p}$
is non-convex, however, \anlchange{approximating} the Pareto set 
\anlchange{can be challenging. Consequently,}  
methods for multi-objective optimization \anlchange{typically pursue 
\emph{Pareto stationarity}, which is a form of} local optimality characterized by a first-order stationarity condition.
\anlchange{If} $\Omegab=\Reals^n$, 
a point $\xb_*$ is a \anlchange{Pareto} stationary point of $\Fb$ provided \anlchange{that for each $\db\in\Reals^n$, there exists
 $j\in\{1,\dots,p\}$} such that $f'_j(\xb_*;\db)\geq 0$.
\anlchange{This notion of stationarity is an extension of the one given for single-objective optimization in \eqref{eq:uncon_Clarke_stationary}.}

\anlchange{Typical methods for \eqref{eq:moo_prob}} return a collection 
\anlchange{of points that are not known to be Pareto-dominated and thus}
serve as an approximation to the 
\anlchange{set of Pareto points.}
\anlchange{From a theoretical point of view, most methods endeavour only to demonstrate that all accumulation \anlchange{points} 
are Pareto stationary, and rarely prove the existence of more than one such point. From a practical point of view, comparing 
the approximate Pareto sets returned by a method for multi-objective optimization}
 is not straightforward.
\anlchange{For discussions of some comparators used in multi-objective
optimization, see \citeasnoun{Knowles2002} and
\citeasnoun{Audet2018b}.}

Various derivative-free methods discussed in this survey have been extended to \anlchange{address
\eqref{eq:moo_prob}.}
\anlchange{The method of \citeasnoun{audet2008multiobjective} solves biobjective optimization problems} by iteratively combining the two objectives
into a single objective (for instance, by considering a weighted sum of the two objectives)\anlchange{; 
\framework{MADS} is then applied to this single-objective problem}. 
\citeasnoun{audet2010mesh} extend \anlchange{the method of  \citeasnoun{audet2008multiobjective}} to multi-objective \anlchange{problems} with more than two objectives. 
\citename{audet2008multiobjective} \citeyear{audet2008multiobjective,audet2010mesh}
demonstrate that all refining points of the sequence \anlchange{of candidate points} produced
by these methods are \anlchange{Pareto stationary}. 

\citeasnoun{Custodio2011} propose \emph{direct-multisearch} methods, a multi-objective analogue of direct-search methods. 
Like direct-search methods, direct-multisearch methods involve both a search step and \anlchange{a poll} step. 
Direct-multisearch methods \anlchange{maintain} a list of non-dominated points;
at the start of an iteration, one non-dominated point must be selected to serve as the centre for a poll step. 
\citeasnoun{Custodio2011} demonstrate that at any accumulation point $\xb_*$ \anlchange{of the maintained sequence of non-dominated points from} a direct-multisearch method,
it holds that 
\anlchange{for any direction $\db$ that appears in a poll step infinitely often},
$f'_j(\xb_*,\db)\geq 0$ for at least one $j$\anlchange{.
In} other words, accumulation points of the method are
\anlchange{Pareto stationary} when restricted to these directions $\db$. 
\citeasnoun{Custodio2017} \anlchange{incorporate} these direct-multisearch
methods \anlchange{within a multistart framework in an effort to find multiple Pareto stationary points and thus to better approximate the Pareto set.} 

\anlchange{For stochastic
  biobjective problems, \citeasnoun{Kim2011} employ sample average approximation to estimate $\Fb(\xb)=\Ea[\xib]{\tilde{\Fb}(\xb;\xib)}$ 
 and propose a model-based trust-region method.}
\citeasnoun{Ryu2014} \anlchange{adapt the approach of \citeasnoun{Kim2011} to the deterministic  biobjective setting.} 
At the start of each iteration, \anlchange{these methods construct} 
fully linear models of both objectives around a (currently) non-dominated point.
\anlchange{These methods solve} three trust-region \anlchange{subproblems} -- one for each of the \anlchange{two
objectives}, and a third that weights the two objectives as in
\citeasnoun{audet2008multiobjective} -- and \anlchange{accept} all non-dominated trial points. 
If both objectives are in $\mathcal{LC}^1$, \citeasnoun{Ryu2014} prove that one of the three 
\anlchange{objectives}
satisfies
a lim-inf convergence result of the form \eqref{eq:liminf_first_order_convergence}, 
implying the existence of a Pareto-stationary accumulation point. 

\citeasnoun{Liuzzi2016} propose a method for constrained multi-objective non-smooth
optimization that separately handles each objective and constraint via an
exact penalty (see \eqref{eq:exact_penalty}) in order to determine whether 
a point is non-dominated.
Given the non-sequential nature of how non-dominated points are selected,
\citeasnoun{Liuzzi2016} 
\anlchange{identify and link the} 
subsequences implied by a lim-inf convergence result.
They show that limit points of \anlchange{these} linked sequences are \anlchange{Pareto stationary}
provided the search directions used in each
linked sequence are asymptotically dense in the unit sphere. 

\citeasnoun{cocchi2018implicit} extend implicit filtering to the multi-objective
case. They approximate \anlchange{each objective gradient
separately using implicit-filtering techniques; they} 
combine these approximate gradients in a disciplined way to 
generate search directions\anlchange{.} 
\citeasnoun{cocchi2018implicit} demonstrate that their method generates at least one
accumulation point and that every such accumulation point is \anlchange{Pareto stationary}.

\subsection{Methods for multifidelity optimization}
Multifidelity optimization concerns the minimization of a high-fidelity 
objective function $f=f_0$ 
\anlchange{in situations where}
a lower-fidelity version $f_\epsilon$ (for 
$\epsilon>0$) also exists. Evaluations of the lower-fidelity function \anlchange{$f_\epsilon$} are 
less computationally expensive than are evaluations of \anlchange{$f_0$}; hence, a goal in 
multifidelity optimization is to 
\anlchange{exploit the existence of}
the lower-fidelity $f_\epsilon$ in order to 
perform as few evaluations of \anlchange{$f_0$} as possible. \anlchange{An} example of such a
\anlchange{setting occurs} when there exist multiple grid resolutions defining 
discretizations for the numerical solution of partial differential equations
\anlchange{that defines $f_0$ and $f_\epsilon$}.

\citeasnoun{polak:650} develop a pattern-search method that exploits the 
existence of multiple levels of fidelity. The method begins at the coarsest 
available level and then monotonically refines the level of fidelity (\ie\ 
decreases $\epsilon$) after a sufficient number of 
\anlchange{consecutive unsuccessful iterations}
occur. 

A method that both decreases $\epsilon$ and increases $\epsilon$ (akin to the 
V- and W-cycles of multigrid \anlchange{methods} \cite{Xu2017}), 
is the multilevel method of 
\citeasnoun{Frandi2013}. The method follows the MG/Opt framework of 
\citeasnoun{Nash2000a} and 
\anlchange{instantiates runs}
of a coordinate-search method at
specified fidelity and solution accuracy levels. Another multigrid-inspired 
method is developed by \citeasnoun{Liu2014}, wherein a hierarchy of \codes{DIRECT} 
\anlchange{runs}
are performed at varying fidelity and budget levels.

Model-based methods have also been extended to the multifidelity setting. For 
example, \citeasnoun{AM2012} employ a fully linear RBF model to interpolate the 
error between \anlchange{two} different fidelity levels. Their method then employs this model 
within a trust-region \anlchange{framework, but uses $f_0$}
to determine whether to accept a given step.

Another model-based approach for multifidelity optimization is
co-kriging; see, for example, \citeasnoun{Xiong2013} and
\citeasnoun{LeGratiet2015}. 
In such approaches, \anlchange{a statistical surrogate (typically a Gaussian process model) 
is constructed for each fidelity level} 
with the aim of modelling the 
relationships among the fidelity levels in areas of the domain relevant \anlchange{to} 
optimization. 

\medskip 

Derivative-free methods for multifidelity, \anlchange{multi-objective and} 
concurrent/parallel optimization remain an especially open avenue of future 
research. 

\vspace{5pt}
\section*{Acknowledgements}
We are grateful to referees and colleagues whose comments greatly improved the
manuscript; these include
Mark Abramson,
Charles Audet,
Ana Lu\'isa Cust\'odio,
S\'ebastien Le Digabel,
Warren Hare,
Paul Hovland,
Jorge Mor\'e,
\anlchange{Raghu Pasupathy,
Margherita Porcelli,}
Francesco Rinaldi,
Cl\'ement Royer,
Katya Scheinberg,
Sara Shashaani,
Aekaansh Verma and
Zaikun Zhang.
We are especially indebted to Gail Pieper and Glennis Starling for their
invaluable editing.
This material is based upon work supported by the applied mathematics and
SciDAC activities of the Office of Advanced Scientific Computing Research,
Office of Science, US Department of Energy, under Contract DE-AC02-06CH11357.

\newpage %
\section*{Appendix: Collection of WCC results}
\label{sec:massive_table}
\addcontentsline{toc}{section}{Appendix: Collection of WCC results}

\tabref{rates} contains select WCC bounds for methods appearing in the
literature.
Given $\epsilon>0$, all WCC bounds in this appendix are given in terms of 
$N_\epsilon$, an upper
bound on the number of \emph{function evaluations} of a method to guarantee 
that the specified condition is met. 
We present results in this form because function evaluation complexity of
derivative-free methods is often of greater interest than is iteration
complexity.
We present $N_\epsilon$ in terms of four parameters:
\begin{itemize}\setlength\itemsep{3pt}
 \item the accuracy $\epsilon$;
 \item the dimension $n$;
 \item the Lipschitz constant of the function $\Lip{f}$, the Lipschitz constant
 of the function gradient $\Lip{g}$ or the Lipschitz constant of the
 function Hessian $\Lip{H}$ (provided these constants are \anlchange{well-defined}); and
 \item a measure of how far the starting point $\xb_0$ is from a stationary 
point $\xb_*$. In this appendix, this measure is either 
 $f(\xb_0)-f(\xb_*)$,
 \begin{equation}
 \label{eq:supR}
 \anlchange{\Rlevel} \defined \sup_{\xb\in\Reals^n} \{\|\xb-\xb_*\|: f(\xb)\leq f(\xb_0)\}
 \end{equation}
 or
 \begin{equation}
 \label{eq:supR2}
 \anlchange{\Rx} \geq \|\xb_0-\xb_*\|. 
 \end{equation}
\end{itemize}
We present additional constants in $N_\epsilon$ when particularly 
informative. 

Naturally, each method in \tabref{rates} has additional algorithmic parameters that 
influence algorithmic behaviour.
We have omitted the dependence of each method's WCC on the selection of 
algorithmic parameters to allow for an easier comparison of methods. 

\anlchange{We recall that, with the exception of the methods from \citeasnoun{Nesterov2015} and \citeasnoun{KP2014}, 
the methods referenced in Table~\ref{table:rates} do not require knowledge of the value of the relevant
Lipschitz constants.}

\afterpage{
 \begin{landscape}
\begin{table}
{\small
 \caption{Known WCC bounds on the number of function evaluations
 needed 
 to achieve a given stationarity measure.\label{table:rates}}
 \begin{tabular*}{46.5pc}{p{175pt}p{175pt}l}
\hline\hline
Rate type & Method type (citation)$^{\textsc{notes}}$ & $N_\epsilon$ \\ \hline 
 ${f\in \cLC^1}$ & \\ \hline
$\|\nabla f(\xb_k)\|\leq\epsilon$ & \framework{DDS} \cite{KP2014} & $\ffrac{n^2 \Lip{g} (f(\xb_0)-f(\xb_*))}{\epsilon^{2}}$ 
 \\[9pt]
 & \framework{TR} \cite{Garmanjani2016} & $\ffrac{n^2\Lip{g}^2(f(\xb_0)-f(\xb_*))}{\epsilon^{2}}$ \\[9pt]
 & \codes{ARC-DFO} \cite{Cartis2012}\trem{A} & $\ffrac{n^2\max\{\Lip{H},\Lip{g}\}^{3/2}(f(\xb_0)-f(\xb_*))}{\epsilon^{3/2}}$\\[9pt]
 \anlchange{$\Ea[\Ub_{k-1}]{\|\nabla f(\hat\xb_k)\|}\leq\epsilon$}
 & \framework{RS} \cite{Nesterov2015}\trem{B} & $\ffrac{n\Lip{g}(f(\xb_0)-f(\xb_*))}{\epsilon^2}$ \\[9pt]
 $\|\nabla f(\xb_k)\|\leq\epsilon$ w.p.\ $1-p_1$ &
 \framework{DDS} \cite{Gratton2015}\trem{C} & 
$\ffrac{mn\Lip{g}^2(f(\xb_0)-f(\xb_*))}{\epsilon^{2}}$
 \\[9pt]
 $\|\nabla f(\xb_k)\|\leq\epsilon$ w.p.\ $1-p_2$ &
 \framework{TR} \cite{Gratton2017complexity}\trem{{\rm C,D}} & $\ffrac{m\max\{\kappaef,\kappaeg\}^2(f(\xb_0)-f(\xb_*))}{\epsilon^{2}}$
 \\ \hline
 ${f\in \cLC^2}$ & \\ \hline
$\max\{\|\nabla f(\xb_k)\|,-\lambda_k\}\leq\epsilon$ & \framework{DDS} 
\cite{Gratton2016} & 
$\ffrac{n^5\max\{\Lip{H},\Lip{g}\}^3(f(\xb_0)-f(\xb_*))}{\epsilon^{3}}$ \\[9pt]
& \framework{TR} \cite{Gratton2017b} & $\ffrac{n^5 \max\{\Lip{H}^3,\Lip{g}^2\}(f(\xb_0)-f(\xb_*))}{\epsilon^3}$ \\[9pt]
 $\max\{\|\nabla f(\xb_k)\|,-\lambda_k\}\leq\epsilon$ w.p.\ $1-p_3$ &
 \framework{TR} \cite{Gratton2017complexity}\trem{{\rm C,D}} & 
$\ffrac{m\max\{\kappaeg,\kappaeH\}^3(f(\xb_0)-f(\xb_*))}{\epsilon^{3}}$
 \\ \hline
 \end{tabular*}
 }
 \end{table}

 \begin{table}
 {\small
 \begin{tabular*}{46.5pc}{p{175pt}p{175pt}l}
 \hspace{-6pt}Table 8.1 continued.\\
\hline\hline
Rate type & Method type (citation)$^{\textsc{notes}}$ & $N_\epsilon$ \\ \hline
\textit{${f\in \cLC^1}$, ${f}$ is ${\lambda}$-strongly convex} & \\ \hline
 $f(\xb_k)-f(\xb_*)\leq\epsilon$
 & \framework{DDS} \cite{KP2014} & $\ffrac{n^2\Lip{g}}{\lambda}\log\biggl(\ffrac{1}{\epsilon}\biggr)$ \\[9pt]
 \anlchange{$\Ea[\Ub_{k-1}]{f(\hat\xb_k)}-f(\xb_*)\leq\epsilon$}
 & \framework{RS} \cite{Nesterov2015}\trem{B} & $\ffrac{n\Lip{g}}{\lambda}\log\biggl(\ffrac{\Lip{g}\anlchange{\Rx}^2}{\epsilon}\biggr)$ \\ \hline
\textit{${f\in \cLC^1}$, ${f}$ is convex} & \\ \hline
 $f(\xb_k)-f(\xb_*)\leq\epsilon$
 & \framework{DDS} \cite{KP2014}\trem{E} & $\ffrac{n^2\Lip{g}\anlchange{\Rlevel}}{\epsilon}$ \\[9pt]
 \anlchange{$\Ea[\Ub_{k-1}]{f(\hat\xb_k)}-f(\xb_*)\leq\epsilon$}
 & \framework{RS} \cite{Nesterov2015}\trem{B} & $\ffrac{n\Lip{g}\anlchange{\Rx}^2}{\epsilon}$ \\ \hline
\textit{${f\in \cLC^0}$, ${f}$ is convex} & \\ \hline
 \anlchange{$\Ea[\Ub_{k-1}]{f(\hat\xb_k)}-f(\xb_*)\leq\epsilon$}
 & \framework{RS} \cite{Nesterov2015}\trem{B} & $\ffrac{n^2\Lip{f}^2\anlchange{\Rx}^2}{\epsilon^2}$ \\ \hline
 ${f\in \cLC^0}$ & \\ \hline
\anlchange{$\Ea[\Ub_{k-1}]{\|\nabla f_{\bar\mu}(\hat\xb_k)\|}\leq\epsilon, 
\bar\mu=\ffrac{\epsilon}{\Lip{f}\sqrt{n}}$} 
& \framework{RS} \cite{Nesterov2015}\trem{B} & $\ffrac{n^3 \Lip{f}^5(f(\xb_0)-f(\xb_*))}{\epsilon^3}$ \\ \hline
 \end{tabular*}
 }
 \end{table}

 \begin{table}
\label{tab:footnotes} {\small
 \begin{tabular*}{46.5pc}{p{175pt}p{175pt}l}
 \hspace{-6pt}Table 8.1 continued.\\
\hline\hline
Rate type & Method type (citation)$^{\textsc{notes}}$ & $N_\epsilon$ \\ \hline
\textit{${f = h\circ \Fb}$, convex ${h\in\cLC^0}$, ${\Fb\in\cLC^{1}}$} & \\ \hline
 $\Psi(\xb_k)\leq\epsilon$ & 
 \framework{TR} \cite{Garmanjani2016}\trem{F} & $\ffrac{p n^2\Lip{g}(\Fb)^2\Lip{f}(h)^2(f(\xb_0)-f(\xb_*))}{\epsilon^2}$ 
 \\ \hline
 \textit{A smoothed ${f_{\mu}(\xb)}$ for ${f}$} & \\ \hline
 $\|\nabla f_{\mu_k}(\xb_k)\|\leq\epsilon$ where 
$\mu_k\in \Bigoh \biggl( \ffrac{\epsilon}{\sqrt{n}}\biggr) $& 
 \framework{DDS} \cite{Garmanjani2012}\trem{G} & $\ffrac{n^{5/2}\left[-\log(\epsilon)+\log(n)\right](f(\xb_0)-f(\xb_*))}{\epsilon^3}$ \\[7pt]
 & \framework{TR} \cite{Garmanjani2016}\trem{G} & $\ffrac{n^{5/2}\left[|\log(\epsilon)|+\log(n)\right](f(\xb_0)-f(\xb_*))}{\epsilon^3}$ \\ \hline
 \end{tabular*}
 }
{\small 
\renewcommand{\tabcolsep}{2.5pt} 
 \begin{tabular*}{46.5pc}{lp{45pc}}
A & We omit an additional $|\log(\epsilon)|$ dependence.\\[3pt]

B & $\hat\xb_k = \argmin_{j=1,\ldots,k} f(\xb_j)$.\\[3pt]

C & $m$ is the number of function evaluations performed in each iteration, 
independent of $n$.\\[3pt]

D & \citeasnoun{Gratton2017complexity} prove results for an arbitrary model-building
scheme that assumes the ability to yield $p$-probabilistically $\kappab_Q$-fully quadratic models 
(where $\kappab_Q=(\kappaef,\kappaeg,\kappaeH)$) when $f\in\cLC^2$
and $p$-probabilistically $\kappab_L$-fully linear models (where $\kappab_L=(\kappaef,\kappaeg)$) 
when $f\in\cLC^1$.
The construction of probabilistically fully quadratic models or probabilistically fully linear models
when $m \ll(n+1)(n+2)/2$ remains an open question. Note that when $p=1$, it is known
that by using $m\in\bigo{n^2}$ %
points, one can guarantee $\kappab_Q$-fully quadratic models with 
$\kappaef,\kappaeg,\kappaeH\in\bigo{n\Lip{H}}$
\cite[Theorem 3]{Conn2006a}.
In this case, the result of \citeasnoun{Gratton2017complexity}
yields a rate weaker than that %
obtained by \citeasnoun{Gratton2016} by a factor of $\Lip{g}$. 
Similarly, when $p=1$, it is known
that by using $m\in\bigo{n}$ %
points, one can guarantee $\kappab_L$-fully linear models with
$\kappaef,\kappaeg\in\bigo{n^{\gfrac{1}{2}}\Lip{g}}$ \cite[Theorem~2]{Conn2006a}.
In this case, the result of \citeasnoun{Gratton2017complexity}
yields a rate comparable to that %
obtained by \citeasnoun{Garmanjani2016}. \\[3pt]

E & \citeasnoun{Vicente2013} derives the same bound but with $\Lip{g}^2$ instead of 
$\Lip{g}$\anlchange{; however, the method of \citeasnoun{Vicente2013} does not require the value $\Lip{g}$.}  \\[3pt]

F & $\Lip{g}(\Fb)$ is the Lipschitz constant of the Jacobian $J(F)$, $\Lip{f}(h)$ 
is the Lipschitz 
constant of $h$, and $p$ is the dimension of the domain of $h$. A bound for a
similar method with an additional $|\log(\epsilon)|$ dependence appears in
\citeasnoun{Grapiglia2016}.\\[3pt]

G & Lipschitz constants do not appear because they are `cancelled' by 
choosing the rate at which smoothing parameter $\mu_k\to 0$. \\ \hline
\end{tabular*}
 }
 \end{table}
 \end{landscape}
 \clearpage
 }

In \tabref{rates}, we employ the constants
\begin{align*}
p_1 &=
\exp\biggl(-\ffrac{n\Lip{g}^2}{\epsilon^{2}} (f(\xb_0)-f(\xb_*)) 
\biggr),
\\*
p_2 &= 
\exp\biggl(-\ffrac{\max\{\kappaef,\kappaeg\}^2}{\epsilon^{2}
} (f(\xb_0)-f(\xb_*)) \biggr),
\\*
p_3 &=
\exp\biggl(-\ffrac{\max\{\kappaeg,\kappaeH\}^3}{\epsilon^{3}
} (f(\xb_0)-f(\xb_*))\biggr). 
\end{align*}

\clearpage
\bibliographystyle{abbrvnat} 
\bibliography{extracted.bib}

\label{lastpage}
\vspace*{\fill}
\begin{flushright}
\scriptsize
\framebox{\parbox{0.95\textwidth}{
The submitted manuscript has been created by UChicago Argonne, LLC, Operator of Argonne National Laboratory (“Argonne”). 
Argonne, a U.S. Department of Energy Office of Science laboratory, is operated under Contract No. DE-AC02-06CH11357. 
The U.S. Government retains for itself, and others acting on its behalf, a paid-up nonexclusive, irrevocable worldwide 
license in said article to reproduce, prepare derivative works, distribute copies to the public, and perform publicly 
and display publicly, by or on behalf of the Government.  The Department of Energy will provide public access to these 
results of federally sponsored research in accordance with the DOE Public Access Plan. 
\url{http://energy.gov/downloads/doe-public-access-plan}.
}}
\normalsize
\end{flushright}

\end{document}